\documentclass{article}

\usepackage{amsmath}
\usepackage{amsthm}
\usepackage{amssymb}
\usepackage{amscd}
\usepackage{xspace}
\usepackage{verbatim}
\newtheorem{theorem}{Theorem}[section]

\newtheorem{lemma}{Lemma}[section]
\newtheorem{proposition}{Proposition}[section]
\theoremstyle{definition}
\newtheorem{definition}{Definition}[section]

\theoremstyle{remark}
\newtheorem{remark}{Remark}[section]
\numberwithin{equation}{section}

\newcommand{\ov}{\overline}

\newcommand{\e}{\varepsilon}

\newcommand{\G}{\Gamma}

\renewcommand{\O}{\Omega}

\renewcommand{\liminf}{\varliminf}
\renewcommand{\limsup}{\varlimsup}

\renewcommand{\vec}[1]{\mathbf{#1}}
\newcommand{\field}[1]{\mathbb{#1}}

\newcommand{\R}{\field{R}}

\newcommand{\er}{\eqref}

\DeclareMathOperator{\Div}{div} 
\DeclareMathOperator{\supp}{supp}

\renewcommand{\O}{\Omega}

\newcommand{\f}{\varphi}
\renewcommand{\vec}[1]{\boldsymbol{#1}}

\date{}

\begin{document}
\title{On the $\Gamma$-limit of singular perturbation problems with optimal profiles which are not one-dimensional. Part III: The
energies with non local terms} \maketitle
\begin{center}
\textsc{Arkady Poliakovsky \footnote{E-mail:
poliakov@math.bgu.ac.il}
}\\[3mm]
Department of Mathematics, Ben Gurion University of the Negev,\\
P.O.B. 653, Be'er Sheva 84105, Israel
\\[2mm]
\end{center}
\begin{abstract}
We use the technique developed in \cite{PI}--\cite{PII} to construct
the upper and the lower bounds for classes of problems containing
non-local terms, including problems in micromagnetics and problems
arising in the variational study of the Method of Vanishing
Viscosity for systems of conservation laws. We reduced these
problems to the problems considered in \cite{PI}--\cite{PII}, with
the appropriate prescribed differential constraint.
\end{abstract}

\section{Introduction}
Consider a family $\{I_\varepsilon\}_{\varepsilon>0}$ of functionals
$I_\varepsilon(\phi):U\to[0,+\infty]$, where $U$ is a given metric
space. The lower and upper $\Gamma$-limits of $I_\varepsilon$ are
defined by:
\begin{align*}
(\Gamma-\liminf_{\varepsilon\to 0^+} I_\varepsilon)(\phi)
:=\inf\left\{\liminf_{\varepsilon\to
0^+}I_\varepsilon(\phi_\varepsilon):\;\,\{\phi_\varepsilon\}_{\varepsilon>0}\subset
U,\; \phi_\varepsilon\to\phi\text{ in }U\;
\text{as}\;\varepsilon\to 0^+\right\},\\
(\Gamma-\limsup_{\varepsilon\to 0^+} I_\varepsilon)(\phi)
:=\inf\left\{\limsup_{\varepsilon\to
0^+}I_\varepsilon(\phi_\varepsilon):\;\,\{\phi_\varepsilon\}_{\varepsilon>0}\subset
U,\; \phi_\varepsilon\to\phi\text{ in }U\; \text{as}\;\varepsilon\to
0^+\right\}.
\end{align*}
In the case where the lower and upper $\Gamma$-limits coincide we
define:
\begin{align*}
(\Gamma-\lim_{\varepsilon\to 0^+}
I_\varepsilon)(\phi):=(\Gamma-\liminf_{\varepsilon\to 0^+}
I_\varepsilon\big)(\phi)=(\Gamma-\limsup_{\varepsilon\to 0^+}
I_\varepsilon)(\phi).
\end{align*}
It is useful to know the $\Gamma$-limit of $I_\varepsilon$, because
it describes the asymptotic behavior as $\varepsilon\downarrow 0$ of
minimizers of $I_\varepsilon$, as follows from the following simple
well known result:
\begin{proposition}[De-Giorgi]\label{propdj}
Assume that $\phi_\varepsilon$ is a minimizer of $I_\varepsilon$ for
every $\varepsilon>0$. Then:
If $I_0(\phi)=(\Gamma-\liminf_{\varepsilon\to 0^+}
I_\varepsilon)(\phi)$ and $\phi_\varepsilon\to\phi_0$ as
$\varepsilon\to 0^+$ then $\phi_0$ is a minimizer of $I_0$.
If $I_0(\phi)=(\Gamma-\lim_{\varepsilon\to 0^+}
I_\varepsilon)(\phi)$ (i.e. it is a full $\Gamma$-limit of
$I_\varepsilon(\phi)$) and for some subsequence $\varepsilon_n\to
0^+$
we have $\phi_{\varepsilon_n}\to\phi_0$, then $\phi_0$ is a
minimizer of $I_0$.
\end{proposition}
Usually, for finding the $\Gamma$-limit of $I_\varepsilon(\phi)$, we
need to find two bounds.
\begin{itemize}
\item[{\bf(I)}] Firstly, we find a lower bound, i.e. a functional
$\underline{I}(\phi)$ such that for every family
$\{\phi_\varepsilon\}_{\varepsilon>0}$, satisfying
$\phi_\varepsilon\to \phi$ as $\varepsilon\to 0^+$, we have
$\liminf_{\varepsilon\to 0^+}I_\varepsilon(\phi_\varepsilon)\geq
\underline{I}(\phi)$.
\item[{\bf(II)}] Secondly, we find an upper
bound, i.e. a functional $\overline{I}(\phi)$, such that for every
$\phi\in U$ there exists a family
$\{\psi_\varepsilon\}_{\varepsilon>0}$, satisfying
$\psi_\varepsilon\to \phi$ as $\varepsilon\to 0^+$ and
$\limsup_{\varepsilon\to 0^+}I_\varepsilon(\psi_\varepsilon)\leq
\overline{I}(\phi)$.
\item[{\bf(III)}] If we find that
$\underline{I}(\phi)=\overline{I}(\phi):=I(\phi)$, then $I(\phi)$ is
the $\Gamma$-limit of $I_\varepsilon(\phi)$.
\end{itemize}

Let $G\in C^1\big(\R^{m\times N^n}\times\R^{m\times
N^{(n-1)}}\times\ldots\times\R^{m\times N}\times \R^m,\R\big)$ and
$W\in C^1(\R^m,\R)$ be nonnegative functions such that
$G\big(0,0,\ldots,0,b)=0$ and let $\Psi\in C^1(\R^m,\R^{l\times
N})$. Consider the energy functional with nonlocal term defined for
every $\e>0$ by
\begin{equation}\label{general1}
I_\e(\phi)=\int_\O\frac{1}{\e}G\Big(\e^n\nabla\phi^n,\ldots,\e\nabla\phi,\phi\Big)\,dx
+\int_\O\frac{1}{\e}W\big(\phi\big)\,dx+\frac{1}{\e}\int_{\R^N}\big|\nabla
\bar H_{\Psi(\phi)}\big|^2
dx\quad\;\;\text{for}\;\;\phi:\O\to\mathcal{M}\subset\R^m\,.
\end{equation}
Here given $u:\O\to\R^{l\times N}$, $\bar H_{u}:\R^N\to\R^l$ is
defined by
\begin{equation}\label{defH12}
\begin{cases}\Delta \bar H_u=div\,\{\chi_\O u\}\quad\text{in the sense of distributions in }\R^N\,,\\
\nabla \bar H_u\in L^2(\R^N,\R^{l\times N})\,,
\end{cases}
\end{equation}
where $\chi_\O$ is the characteristic function of $\O$. One of the
fields where functionals of type \er{general1} are relevant is
Micromagnetics (see \cite{ARS}, \cite{DKMO}, \cite{RS1}, \cite{RS2}
and other). The full 3-dimensional model of ferromagnetic materials
deals with an energy functional, which, up to a rescaling, has the
form
\begin{equation}\label{micro2}
E_\e(m):=\e\int_\O|\nabla m|^2dx+\frac{1}{\delta_\e}\int_\O
W(m)dx+\frac{1}{\e}\int_{\R^3}|\nabla \bar H_{m}|^2dx\,,
\end{equation}
where $\O\subset\R^3$ is a bounded domain, $m:\O\to S^2$ stands for
the magnetization, $\delta_\e>0$ is a material parameter and $\bar
H_{m}:\R^3\to\R$ is defined, as before, by
\begin{equation}\label{defH123}
\begin{cases}\Delta \bar H_m=div\,\{\chi_\O m\}\quad\text{in }\R^3\,,\\
\nabla \bar H_m\in L^2(\R^3,\R^3)\,,
\end{cases}
\end{equation}
The first term in \er{micro2} is usually called the exchange energy
while the second is called the anisotropy energy and the third is
called the demagnetization energy.
One can consider the infinite cylindrical domain $\O=G\times\R$ and
configurations which don't depend on the last coordinate. These
reduce the original model to a 2-dimensional one, where the energy,
up to a rescaling, has the form
\begin{equation}\label{micro3}
E_\e(m):=\e\int_G|\nabla m|^2dx+\frac{1}{\delta_\e}\int_G
W(m)dx+\frac{1}{\e}\int_{\R^2}|\nabla \bar H_{m'}|^2dx\,,
\end{equation}
where $G\subset\R^2$ is a bounded domain, $m=(m_1,m_2,m_3):G\to S^2$
stands for the magnetization, $m':=(m_1,m_2)\in\R^2$ denotes the
first two components of $m$, $\delta_\e>0$ and $\bar
H_{m'}:\R^2\to\R$ is defined, as before, by
\begin{equation}\label{defH124}
\begin{cases}\Delta \bar H_{m'}=div\,\{\chi_G m'\}\quad\text{in }\R^2\,,\\
\nabla \bar H_{m'}\in L^2(\R^2,\R^2)\,.
\end{cases}
\end{equation}
Note that in the case $\delta_\e=\e$ (i.e. the anisotropy and the
demagnetization energies have the same order as $\e\to 0$) the
energy-functionals in \er{micro2} and \er{micro3} are special cases
of the energy in \er{general1}.

 In this work, using the technique developed in \cite{PI} and \cite{PII}, we construct the upper and the lower bounds as
$\e\downarrow 0$ for the general energy of the form \er{general1}
under certain conditions on set $\mathcal{M}$ for functions $\phi\in
BV$. In particular our upper bound improves, in general, one
obtained in \cite{polmag}.

  In order to reduce problem \er{general1} to the local problems studied in \cite{PI} and \cite{PII},
the following trivial observation was made for problem
\er{general1}. For $\phi:\O\to\mathcal{M}$, such that $\Psi(\phi)\in
L^2$, consider the variational problem
\begin{equation}\label{defH12sldfgvrftghbjhgh}
J_{\Psi,\phi}(L):=\inf\bigg\{
P_{\Psi,\phi}(L):=\int_{\R^{N}}\Big|L(x)+\chi_\O(x)
\Psi\big(\phi(x)\big)\Big|^2dx:\;\;L\in L^2(\R^N,\R^{l\times N}),\;
div \,L\equiv 0\bigg\}\,.
\end{equation}
Then
\begin{equation}\label{defH12sldfgvrftghbjhghjh}
J_{\Psi,\phi}(L)=\int_{\R^N}\big|\nabla \bar
H_{\Psi(\phi)}(x)\big|^2dx\,,
\end{equation}
where
given $u:\O\to\R^{l\times N}$, $\bar H_{u}:\R^N\to\R^l$ is defined
by \er{defH12}.
Moreover
\begin{equation}\label{defH12sldfgvrftghbjhghjkk}
L_0(x):=\nabla \bar H_{\Psi(\phi)}(x)-\chi_\O(x)
\Psi\big(\phi(x)\big)\quad\quad\text{is a minimizer to
\er{defH12sldfgvrftghbjhgh}}\,.
\end{equation}
Therefore the $\Gamma$-limit of the family of functionals
\er{general1}, as $\phi_\e\to\phi$, where
\begin{equation}\label{hggnwsl88ijii9nh}
\begin{cases}
W(\phi)=0\\ \Div \Psi(\phi)=0\quad\text{in}\;\O\\
\Psi(\phi)\cdot\vec
n=0\quad \text{on}\;\partial\O,
\end{cases}
\end{equation} is the same as the
$\Gamma$-limit of the family of functionals
\begin{multline}\label{general1sljgjkf}
\bar I_\e(\phi,L)=:
\int_\O\frac{1}{\e}G\Big(\e^n\nabla\phi^n,\ldots,\e\nabla\phi,\phi\Big)\,dx
+\int_\O\frac{1}{\e}W\big(\phi\big)\,dx+\frac{1}{\e}\int_{\R^{N}}\Big|L(x)+\chi_\O(x)
\Psi\big(\phi(x)\big)\Big|^2dx\\
\text{where}\quad\phi:\O\to\mathcal{M},\;\;\text{and}\;\; div\,
L=0\,,
\end{multline}
as $(\phi_\e,L_\e)\to\big(\phi,-\chi_\O \Psi(\phi)\big)$. More
generally we have the following simple result (see also Lemma
\er{vfgcfgjfvhgkjkj} as a generalization).
\begin{proposition}\label{vfgcfgjfvhgkjkjghjgjkjhkjgk}
Let $\Psi(\psi)\in C(\R^m,\R^{l\times N})$  which satisfies
\begin{equation}\label{ghhghjghjhbfggffggjgj}
\big|\Psi(\psi)\big|\leq C_0|\psi|^{p/2}\quad\forall\psi\in\R^m\,,
\end{equation}
for some constant $C_0>0$ and $p\geq 1$. Furthermore, for every
$\e>0$ consider the functional
$E_\e\big(\phi(x)\big):L^p(\O,\R^m)\to[0,+\infty)\cup\{+\infty\}$
which (possibly) can attain the infinite values. Next for every
$\e>0$ consider the functional
$P_\e\big(\phi(x)\big):L^p(\O,\R^m)\to[0,+\infty)\cup\{+\infty\}$,
defined by
\begin{equation}\label{fhfhjffhjhffhjh}
P_\e\big(\phi(x)\big):=E_\e\big(\phi(x)\big)+\frac{1}{\delta_\e}\int_{\R^N}\Big|\nabla_x
\bar H_{\Psi(\phi)}(x)\Big|^2 dx\,,
\end{equation}
where $\delta_\e>0$ satisfies $\lim_{\e\to 0^+}\delta_\e=0$ and
given $\phi(x)\in L^p(\O,\R^m)$, $\bar
H_{\Psi(\phi)}(x):\R^N\to\R^l$ is defined by
\begin{equation}\label{defH12slvvvbghhjk}
\begin{cases}\Delta_x \bar
H_{\Psi(\phi)}(x)=\Div_x\Big\{\chi_\O(x)\cdot \Psi\big(\phi(x)\big)\Big\}\quad\text{in }\R^N\,,\\
\nabla_x \bar H_{\Psi(\phi)}(x)\in L^2(\R^N,\R^{l\times N})\,.
\end{cases}
\end{equation}
Furthermore, for every $\e>0$ consider the
functional $Q_\e\big(\phi(x),L(x)\big):L^p(\O,\R^m)\times
L^2(\R^N,\R^{l\times N}) \to[0,+\infty)\cup\{+\infty\}$ defined by
\begin{equation}\label{fhfhjffhjhfgyukuoijkkjbmjg}
Q_\e\big(\phi(x),L(x)\big):=\begin{cases}E_\e\big(\phi(x)\big)+\frac{1}{\delta_\e}\int_{\R^N}\Big|L(x)+\chi_\O(x)\cdot
\Psi\big(\phi(x)\big)\Big|^2 dx\quad\quad\text{if}\quad \Div_x
L(x)\equiv 0\,,\\+\infty\quad\quad\quad\text{otherwise}\,.
\end{cases}
\end{equation}
Then for every $\phi(x)\in L^p(\O,\R^m)$, such that
$\Div_x\big\{\chi_\O(x)\cdot \Psi\big(\phi(x)\big)\big\}=0$ in
$\R^N$ we have the following equalities
\begin{multline}\label{hfghjfhfhfhhnkh}
\inf\Bigg\{\liminf_{\e\to
0^+}P_\e\big(\phi_\e(x)\big):\;\;\phi_\e(x)\to\phi(x)\;\;\text{in}\;\;L^p(\O,\R^m)\Bigg\}=
\inf\Bigg\{\liminf_{\e\to
0^+}Q_\e\big(\phi_\e(x),L_\e(x)\big):\\\;\;\phi_\e(x)\to\phi(x)\;\;\text{in}\;\;L^p(\O,\R^m)\;\;
\text{and}\;\;\;L_\e(x)\to\Big(-\chi_\O(x)\cdot
\Psi\big(\phi(x)\big)\Big)\;\;\text{in}\;\;L^2(\R^N,\R^{l\times
N})\Bigg\}\,,
\end{multline}
and
\begin{multline}\label{hfghjfhfhfhj;k;jkhkhhkkkjgk}
\inf\Bigg\{\limsup_{\e\to
0^+}P_\e\big(\phi_\e(x)\big):\;\;\phi_\e(x)\to\phi(x)\;\;\text{in}\;\;L^p(\O,\R^m)\Bigg\}=
\inf\Bigg\{\limsup_{\e\to 0^+}Q_\e\big(\phi_\e(x),L_\e(x)\big):\\
\phi_\e(x)\to\phi(x)\;\;\text{in}\;\;L^p(\O,\R^m)\;\;
\text{and}\;\;\;L_\e(x)\to\Big(-\chi_\O(x)\cdot
F\big(\phi(x)\big)\Big)\;\;\text{in}\;\;L^2(\R^N,\R^{l\times
N})\Bigg\}\,.
\end{multline}
%
%
%
%
%
%
%
%
\end{proposition}
Next since the energy \er{general1sljgjkf} with $\mathcal{M}\equiv
\R^m$ is a particular case of the functionals studied in \cite{PI},
where we get its upper bound  and in \cite{PII}, where we get its
lower bound, we can apply this results to problem \er{general1}.
Then we get the following Theorems providing the upper and the lower
bound (see Theorems
\ref{dehgfrygfrgygenbgggggggggggggkgkgthtjtfjkjkkhhk} and
\ref{dehgfrygfrgygenbgggggggggggggkgkgthtjtfjkjkkhhkjfjgfgghnew} for
the proof).
\begin{theorem}\label{dehgfrygfrgygenbgggggggggggggkgkgthtjtfjkjkkhhkggkjj}
Let $\O\subset\R^N$ be an open set with locally Lipschitz's
boundary, let $G\in C^1\big(\R^{m\times N^n}\times\R^{m\times
N^{(n-1)}}\times\ldots\times\R^{m\times N}\times \R^m,\R\big)$ and
$W\in C^1(\R^m,\R)$ be nonnegative functions such that
$G\big(0,0,\ldots,0,b)=0$ and let $\Psi\in C^1(\R^m,\R^{l\times
N})$. Furthermore, let $\varphi\in BV(\R^N,\R^{m})\cap L^\infty$ be
such that $\|D \varphi\|(\partial\Omega)=0$, $W\big(
\varphi(x)\big)=0$ for a.e. $x\in\O$, $\Div_x
\Psi\big(\varphi(x)\big)=0$ in $\O$ and
$\Psi\big(\varphi(x)\big)\cdot\vec n(x)=0$ on $\partial\O$.
Then there exists a sequence $\{\psi_\e\}_{\e>0}\subset
C^\infty(\R^N,\R^m)$ such that
$\int_\O\psi_\e(x)dx=\int_\O\varphi(x)dx$, for every $q\geq 1$ we
have $\lim_{\e\to 0^+}\psi_\e=\varphi$ in $L^q$ and
\begin{multline}\label{a1a2a3a4a5a6a7s8hhjhjjhjjjjjjkkkkgenhjhhhhjtjurtgfhfhfjfjfjhjjjnjhjjjgjkjgj}
\limsup_{\e\to
0^+}I_\e(\psi_\e):=\int_\O\frac{1}{\e}G\Big(\e^n\nabla\psi_\e^n,\ldots,\e\nabla\psi_\e,\psi_\e\Big)\,dx
+\int_\O\frac{1}{\e}W\big(\psi_\e\big)\,dx+\frac{1}{\e}\int_{\R^N}\big|\nabla
\bar H_{\Psi(\psi_\e)}\big|^2 dx\\ \leq\int_{\O\cap J_\varphi}\hat
E_{per}\Big(\varphi^+(x),\varphi^-(x),\vec \nu(x)\Big)d \mathcal
H^{N-1}(x)\,,
\end{multline}
where $\bar H_{u}$ is defined by \er{defH12},
\begin{multline}\label{L2009hhffff12kkkhjhjghghgvgvggcjhggghtgjuthjhhhjjhlhhj}
\hat E_{per}\Big(\varphi^+,\varphi^-,\vec \nu\Big):=
\inf\Bigg\{\int_{I_{\vec \nu}}\frac{1}{L} \bigg(G\big(L^n\nabla^n
\zeta(y),\ldots,L\nabla
\zeta(y),\zeta(y)\big)\,dy+W\big(\zeta(y)\big)+ \big|\nabla
H_{\Psi,\zeta,\vec\nu}(y)\big|^2\bigg)\,dy:
\\
L\in(0,+\infty)\,,\;\zeta\in
\mathcal{S}(\varphi^+,\varphi^-,I_{\vec\nu})\Bigg\}\,,
\end{multline}
where $H_{\Psi,\zeta,\vec\nu}\in W^{2,2}_{loc}(\R^N,\R^l)$ satisfies
\begin{equation}\label{defH12slbgghfjfghjhjhjhjjkjkjhhjjjkj}
\begin{cases}\Delta_y H_{\Psi,\zeta,\vec\nu}(y)=\Div_y \Psi\big(\zeta(y)\big)\quad\text{in }I_{\vec\nu}\,,\\
H_{\Psi,\zeta,\vec\nu}(y+\vec\nu_j)=H_{\Psi,\zeta,\vec\nu}(y)\quad\forall
y\in\R^N\;\;\text{such that}\;\;|y\cdot\vec\nu|<1/2\,,\\
\frac{\partial}{\partial\vec\nu}H_{\Psi,\zeta,\vec\nu}(y)=0\quad\forall
y\in\R^N\;\;\text{such that}\;\;|y\cdot\vec\nu|=1/2\,,
\end{cases}
\end{equation}
and
\begin{multline}\label{L2009Ddef2hhhjjjj77788hhhkkkkllkjjjjkkkhhhhffggdddkkkgjhikhhhjjddddhdkgkkkhgghjhjhjkjtjytrjghjghjjjgjhjg}
\mathcal{S}\big(\varphi^+,\varphi^-,I_{\vec \nu}\big):=
\bigg\{\zeta\in C^n(\R^N,\R^m):\;\;\zeta(y)=\varphi^-\;\text{ if }\;y\cdot\vec \nu\leq-1/2,\\
\zeta(y)=\varphi^+\;\text{ if }\; y\cdot\vec \nu\geq 1/2\;\text{ and
}\;\zeta\big(y+\vec \nu_j\big)=\zeta(y)\;\;\forall j=2,\ldots,
N\bigg\}\,,
\end{multline}
Here
$I_{\vec \nu}:=\big\{y\in\R^N:\;|y\cdot \vec\nu_j|<1/2\;\;\forall
j=1,\ldots, N\big\}$,
where $\{\vec\nu_1,\ldots,\vec\nu_N\}\subset\R^N$ is an orthonormal
base in $\R^N$ such that $\vec\nu_1:=\vec \nu$.
\end{theorem}
\begin{theorem}\label{dehgfrygfrgygenbgggggggggggggkgkgthtjtfjkjkkhhkjfjgfgghnewgghvh}
Let $\O\subset\R^N$ and $G$, $W$, $\Psi$ and $\varphi$ be the same
as in Theorem
\ref{dehgfrygfrgygenbgggggggggggggkgkgthtjtfjkjkkhhkggkjj}.
Moreover, assume that there exists a constant $C>0$ and $p\geq 1$
such that $\big|\Psi(b)\big|\leq C\big(|b|^{p/2}+1\big)$ for every
$b\in\R^m$ and $|a_n|^p/C\leq G(a_n,\ldots, a_2,a_1,b)+W(b)\leq
C\big(\sum_{j=1}^n|a_j|^p+|b|^p+1\big)$ for every $a_j\in\R^{m\times
N^j}$ and $b\in\R^m$.
Then for every sequence $\{\varphi_\e\}_{\e>0}\subset
W^{n,p}_{loc}(\O,\R^m)$
such that $\varphi_\e\to \varphi$ in $L^p_{loc}(\O,\R^m)$ as $\e\to
0^+$, we have
\begin{multline}\label{a1a2a3a4a5a6a7s8hhjhjjhjjjjjjkkkkgenhjhhhhjtjurtgfhfhfjfjfjhjjjnjhjjjhdgfgffgnewjhjkgjk}
\liminf_{\e\to 0^+}I_\e(\varphi_\e):=\liminf_{\e\to
0^+}\Bigg\{\frac{1}{\e}\bigg(\int_\O G\Big(\e^n\nabla^n
\varphi_\e(x),\ldots,\e\nabla
\varphi_\e(x),\varphi_\e(x)\Big)+W\big(\varphi_\e(x)\big)\bigg)dx+\frac{1}{\e}\int_{\R^N}|\nabla
\bar H_{\Psi(\varphi_\e)}|^2dx\Bigg\}\\ \geq\int_{\O\cap
J_\varphi}\hat E_0\Big(\varphi^+(x),\varphi^-(x),\vec \nu(x)\Big)d
\mathcal H^{N-1}(x)\,,
\end{multline}
where $\bar H_{u}$ is defined by \er{defH12} and
\begin{multline}\label{L2009hhffff12kkkhjhjghghgvgvggcjhggghtgjuthjhhhjjkgkjggjknewjhgh}
\hat E_{0}\Big(\varphi^+,\varphi^-,\vec \nu\Big)\;:=\; \inf\Bigg\{
\liminf_{\e\to 0^+}\int_{I_{\vec \nu}}\frac{1}{\e}
\bigg(G\big(\e^n\nabla^n \zeta_\e(y),\ldots,\e\nabla
\zeta_\e(y),\zeta_\e(y)\big)+W\big(\zeta_\e(y)\big)+ \big|\nabla
H^0_{\Psi,\zeta_\e,\vec\nu}(y)\big|^2\bigg)\,dy:
\\
\zeta_\e\in\mathcal{S}\big(\varphi^+,\varphi^-,I_{\vec\nu}\big)
\;\;\text{such
that}\;\;\zeta_\e(y)\to\chi(y,\varphi^+,\varphi^-,\vec\nu)\;\;\text{in}\;\;L^p(I_{\vec
\nu},\R^m)\Bigg\}\,,
\end{multline}
where $H^0_{\Psi,\zeta,\vec\nu}\in W^{1,2}_{0}(I_{\vec \nu},\R^k)$
satisfies
\begin{equation}\label{defH12slbgghfjfghjhjhjhjjkjkjhhjjgjggnewghghf}
\Delta_y H^0_{\Psi,\zeta,\vec\nu}(y)=\Div_y
\Psi\big(\zeta(y)\big)\quad\text{in }I_{\vec\nu}\,,
\end{equation}
$\mathcal{S}\big(\varphi^+,\varphi^-,I_{\vec\nu}\big)$ is defined by
\er{L2009Ddef2hhhjjjj77788hhhkkkkllkjjjjkkkhhhhffggdddkkkgjhikhhhjjddddhdkgkkkhgghjhjhjkjtjytrjghjghjjjgjhjg}
and \begin{equation*}
\chi(y,\varphi^+,\varphi^-,\vec\nu):=\begin{cases}
\varphi^+\quad\text{if}\;\;y\cdot\vec\nu>0\,,\\
\varphi^-\quad\text{if}\;\;y\cdot\vec\nu<0\,.
\end{cases}
\end{equation*}
Here $I_{\vec \nu}:=\big\{y\in\R^N:\;|y\cdot
\vec\nu_j|<1/2\;\;\forall j=1,\ldots, N\big\}$ where
$\{\vec\nu_1,\ldots,\vec\nu_N\}\subset\R^N$ is an orthonormal base
in $\R^N$ such that $\vec\nu_1:=\vec \nu$.
\end{theorem}
As the boundary conditions for $H_{\Psi,\zeta,\vec\nu}$ in
\er{defH12slbgghfjfghjhjhjhjjkjkjhhjjjkj} are different from those
for $H^0_{\Psi,\zeta,\vec\nu}$, there is a natural question either
in general upper bound obtained in Theorem
\ref{dehgfrygfrgygenbgggggggggggggkgkgthtjtfjkjkkhhkggkjj} coincides
with the lower bound obtained in Theorem
\ref{dehgfrygfrgygenbgggggggggggggkgkgthtjtfjkjkkhhkjfjgfgghnewgghvh}.
The answer yes will mean that we will find the full $\Gamma$-limit
of $I_\e$ in the case of $BV\cap L^\infty$ limiting functions. The
equivalent question is either
\begin{equation}\label{hdiohdo3gjkjgg}
\hat E_{per}\big(\varphi^+,\varphi^-,\vec \nu\big)=\hat
E_{0}\big(\varphi^+,\varphi^-,\vec \nu\big)\,,
\end{equation}
where $E_{per}(\cdot)$ is defined in
\er{L2009hhffff12kkkhjhjghghgvgvggcjhggghtgjuthjhhhjjhlhhj} and
$E_{0}(\cdot)$ is defined by
\er{L2009hhffff12kkkhjhjghghgvgvggcjhggghtgjuthjhhhjjkgkjggjknewjhgh}.
As we showed in \cite{polmag} this is indeed the case when $m=1$.
Moreover, in the later case the optimal profiles are one
dimensional. It can be shown that the question in the general case
is equivalent to the question of equality of upper and lower bound
arisen in \cite{PI}, \cite{PII}.

Section \ref{hgvfghfghoj} is devoted to the variational formulation
of the Method of Vanishing Viscosity for systems of Conservation
Laws. Let $F(u)=\{F_{ij}(u)\}\in
C^1({\mathbb{R}}^k,{\mathbb{R}}^{k\times N})$ with $F(0)=0$.
Consider a Cauchy problem for a system of conservation laws
\begin{equation}\label{equ24gggg}
\begin{cases}\partial_t u+div_x F(u)=0\quad\forall
(x,t)\in{\mathbb{R}}^N\times (0,+\infty)
\\
u(x,0)=v_0(x)\quad\forall x\in\mathbb{R}^N,
\end{cases}
\end{equation}
We say that the function $\eta(u)\in C^2(\mathbb{R}^k,\mathbb{R})$
is an entropy for the system \eqref{equ24gggg} and
$\Psi(u):=\big(\Psi_1(u),\ldots,\Psi_N(u)\big)\in
C^1(\mathbb{R}^k,\mathbb{R}^N)$ is an entropy flux associated with
$\eta$ if we have
$\nabla_u\Psi_j(u)=\nabla_u\eta(u)\cdot\nabla_u \bar F_j(u)$ for all
$u\in \mathbb{R}^k$ and $j=1,\ldots,N$,
where $\bar
F_j(u):=\big(F_{1j}(u),\ldots,F_{kj}(u)\big):\mathbb{R}^k\to\mathbb{R}^k\;$
$\forall j=1,\ldots,N$.

 Since we have a lack of uniqueness of solution to \eqref{equ24gggg},
we need to choose an admissible solution. Due to the method of
vanishing viscosity, given a fixed entropy $\eta(u)$, a solution
$u(x,t)$ to \eqref{equ24gggg} is admissible, if
$u=\lim_{\varepsilon\to 0^+}u_\varepsilon$ in $L^q$, where
$u_\varepsilon(x,t)$ is a solution to the following system:
\begin{equation}\label{bdgvfghffhkkgkkgg}
\begin{cases}\varepsilon\Delta_x\big(
\nabla_u\eta\big(u(x,t)\big)\big)=\partial_t u(x,t)+div_x
F\big(u(x,t)\big)\quad \forall(x,t)\in{\mathbb{R}}^N\times(0,T),
\\ u(x,0)=v_0(x)\quad\quad\quad\quad\quad\forall x\in{\mathbb{R}}^N,
\end{cases}
\end{equation}
(usually, for symmetric conservation laws, one chooses
$\eta(u)=|u|^2/2$ and thus $\Delta_x\big(
\nabla_u\eta\big(u(x,t)\big)\big)=\Delta_x u(x,t)$).
Consider the following family of energy functionals
$\{I_{\e,F}(u)\}$, defined for $u(x,t):\R^N\times[0,T]\to\R^k$ by
\begin{equation}\label{conskgjjk}
I_{\e,F}(u):=\int_0^T\int_{\R^N}\bigg(\e\Big|\nabla_x
\big\{\nabla_u\eta\big(u(x,t)\big)\big\}\Big|^2+\frac{1}{\e}\Big|\nabla_x
H_{F,u}(x,t)\Big|^2\bigg)\,dxdt+\int_{\R^N}\eta\big(u(x,T)\big)\,dx\,,
\end{equation}
where $H_{F,u}(x,t)$
satisfies
\begin{equation}\label{conssstgkfgjhgj}
\begin{cases}\Delta_x H_{F,u}(x,t)=\partial_t u(x,t)+div_x
F\big(u(x,t)\big)\,,\\
\nabla_x H_{F,u}(x,t)\in L^2\big(0,T;L^2(\R^N,\R^{k\times
N})\big)\,,
\end{cases}
\end{equation}
and we assume that
\begin{equation}\label{jhniguiguyyjkjgkkhjh}
u(x,t)\in L^2\big(0,T;W^{1,2}_0(\R^N,\R^k)\big)\cap
C\big(0,T;L^2(\R^N,\R^k)\big)\cap L^\infty\;\;\text{and}\;\;
\partial_t u(x,t)\in L^2\big(0,T;W^{-1,2}(\R^N,\R^k)\big)\,,
\end{equation}
Since
\begin{multline*}
-\int_0^T\int_{\R^N}\nabla_x
\big\{\nabla_u\eta\big(u(x,t)\big)\big\}:\nabla_x
H_{F,u}(x,t)\,dxdt=\int_0^T\int_{\R^N}
\nabla_u\eta\big(u(x,t)\big)\cdot\Delta_x H_{F,u}(x,t)\,dxdt\\=
\int_0^T\int_{\R^N} \nabla_u\eta\big(u(x,t)\big)\cdot\Big(\partial_t
u(x,t)+div_x
F\big(u(x,t)\big)\Big)\,dxdt\\=\int_{\R^N}\bigg(\int_0^T
\partial_t\big\{\eta\big(u(x,t)\big)\big\}\,dt\bigg)dx+ \int_0^T\int_{\R^N}\sum_{j=1}^{N}\nabla_u\eta\big(u(x,t)\big)\cdot
\nabla_u\bar F_j\big(u(x,t)\big)\cdot\frac{\partial u}{\partial
x_j}(x,t)\,dxdt\\=\int_{\R^N}\Big(\eta\big(u(x,T)\big)-\eta\big(u(x,0)\big)\Big)\,dx+
\int_0^T\int_{\R^N}\Div_x\Psi\big(u(x,t)\big)\,dxdt=\int_{\R^N}\Big(\eta\big(u(x,T)\big)-\eta\big(u(x,0)\big)\Big)\,dx\,,
\end{multline*}
we can rewrite the expression of $I_{\e,F}(u)$ as
\begin{equation}\label{conslllgkggf}
I_{\e,F}(u)=\int_0^T\int_{\R^N}\frac{1}{\e}\bigg|\e\nabla_x
\big\{\nabla_u\eta\big(u(x,t)\big)\big\}-\nabla_x
H_{F,u}(x,t)\bigg|^2\,dxdt+\int_{\R^N}\eta\big(u(x,0)\big)\,dx\,,
\end{equation}
Thus if there exists a solution of \er{bdgvfghffhkkgkkgg}
for some $v_0(x)\in L^2(\R^N,\R^k)\cap L^\infty$ then, by
\er{conslllgkggf}, $u(x,t)$ is also a minimizer to
\begin{equation}\label{bdgvfghffhkkjkgjkgjkgkgkgkgjkgj}
\inf\big\{I_{\e,F}(u): u(x,0)=v_0(x)\big\}\,.
\end{equation}
Moreover, in this case,
\begin{equation}\label{bdgvfghffhkkjkgjkgjkgkgkgkjhnghjgffh}
\inf\big\{I_{\e,F}(u):
u(x,0)=v_0(x)\big\}=\int_{\R^N}\eta\big(v_0(x)\big)\,dx\,,
\end{equation}
and the function $u(x,t):\R^N\times[0,T]\to\R^k$ is a minimizer to
\er{bdgvfghffhkkjkgjkgjkgkgkgkgjkgj} if and only if $u(x,t)$ is a
solution to \er{bdgvfghffhkkgkkgg}. On the other hand, it is clear
that if minimizers $u_\e$ of
\er{bdgvfghffhkkjkgjkgjkgkgkgkjhnghjgffh} strongly converges in
$L^q$ to some $u$, then $u$ is a solution of
\begin{equation}\label{equ24ggggkhh}
\begin{cases}
\partial_t u+div_x F(u)=0\quad\forall (x,t)\in\R^N\times
(0,T)\,\\ u(x,0)=v_0(x)\,.
\end{cases}
\end{equation}
Thus, it is a natural question in the Method of Vanishing Viscosity
for Conservation Laws to know the $\Gamma$-limit of the functional
\begin{equation}\label{consdfbhgyjykukukluklgjghgkgkgk}
J_{\e,F,v_0}(u)=
\begin{cases}I_{\e,F}(u)\quad\quad
\text{if}\quad u(x,0)\equiv v_0(x),\\
+\infty\quad\quad\quad\;\text{otherwise}\,.\end{cases}
\end{equation}
Indeed, if for a given function $u$ we have $\Gamma-\liminf_{\e\to
0^+} I_{\e,F,v_0}(u)<+\infty$ then there exists a sequence $u_\e$
such that $u_\e\to u$ in $L^q$ and $\liminf_{\e\to
0^+}I_{\e,F,v_0}(u_\e)<+\infty$ and therefore, by \er{conskgjjk} we
obtain $\nabla_x H_{F,u}\equiv 0$ for the given $u$ that implies
\er{equ24ggggkhh}, i.e. the functional $\Gamma-\liminf_{\e\to 0^+}
I_{\e,F,v_0}(u)$ attains finite values only on the weak solutions of
\er{equ24ggggkhh}. On the other hand, from the theory of
$\Gamma$-limits it is clear that the admissible solution $u$ is a
minimizer of the $\Gamma-\liminf_{\e\to 0^+} I_{\e,F,v_0}(u)$. So,
as we have a lack of uniqueness of solution to \er{equ24ggggkhh}, we
get a necessary condition for a given weak solution of
\er{equ24ggggkhh} to be admissible: the admissible solution is
obligated to be a minimizer of the $\Gamma$-limit energy
$\Gamma-\liminf_{\e\to 0^+} I_{\e,F,v_0}$.
The question of $\Gamma$-limit for $I_{\e,F,v_0}$ was arisen in
\cite{P4}. In \cite{polmag} we found its upper bound, achieved by
one-dimensional profiles. Moreover, we showed that this bound
coincides with the $\Gamma$-limit in the case $k=1$ i.e. in the case
of scalar Conservation Law. In this paper we improve the upper bound
in the case of systems and we construct also the lower bound.

 As before, we can reduce the problem \er{conskgjjk} to local problems studied in \cite{PI} and \cite{PII}.
Indeed assume that $h_0:\R^N\to\R^{k\times N}$ satisfies $div_x
h_0(x)\equiv v_0(x)$.
Then set $L_u:\R^N\times(0,T)\to\R^{k\times N}$ by
\begin{equation}\label{defH12sldfgvrftghbjhghjkkkfff}
L_u(x,t):=h_0(x)+\int_0^t\Big\{\nabla_x
H_{F,u}(x,s)-F\big(u(x,s)\big)\Big\}ds\,.
\end{equation}
where $H_{F,u}(x,t)$ satisfies \er{conssstgkfgjhgj}. So
$L_u(x,0)=h_0(x)$ and $\partial_t L_u(x,t)=\nabla_x
H_{F,u}(x,t)-F\big(u(x,t)\big)$. Thus $div_x L_u(x,0)=v_0(x)$ and
$\partial_t div_x L_u(x,t)=\Delta_x H_{F,u}(x,t)-div_x
F\big(u(x,t)\big)=\partial_t u(x,t)$. Therefore we get
\begin{equation}\label{defH12sldfgvrftghbjhghjkkkjkhifff}
div_x L_u(x,t)=u(x,t),\quad div_x
L_u(x,0)=v_0(x)\quad\text{and}\quad\nabla_x H_{F,u}(x,t)=\partial_t
L_u(x,t)+F\big(div_x L_u(x,t)\big)\,.
\end{equation}
Then we can rewrite the energy in \er{conskgjjk} as
\begin{multline}\label{consyjfff}
I_{\e,F}(u)=\frac{1}{2}\int_0^T\int_{\R^N}\bigg(\e\Big|\nabla_x\big\{
\nabla\eta(u)\big\}\Big|^2+\frac{1}{\e}\big|\nabla_x
H_{F,u}\big|^2\bigg)\,dxdt+\frac{1}{2}\int_{\R^N}\eta\big(u(x,T)\big)\,dx=\\
\frac{1}{2}\int_0^T\int_{\R^N}\bigg(\e\Big|\nabla_x\big\{\nabla\eta(
div_x L_u)\big\}\Big|^2+\frac{1}{\e}\Big|\partial_t L_u+F\big(div_x
L_u\big)\Big|^2\bigg)\,dxdt+\frac{1}{2}\int_{\R^N}\eta\big(div_x
L_u(x,T)\big)\,dx\,.
\end{multline}
On the other hand, define
\begin{multline}\label{consyjnkgfngjgjfff}
\bar I_{\e,F}(L):=
\frac{1}{2}\int_0^T\int_{\R^N}\bigg(\e\Big|\nabla_x\big\{\nabla\eta(
div_x L)\big\}\Big|^2+\frac{1}{\e}\Big|\partial_t L+F\big(div_x
L\big)\Big|^2\bigg)\,dxdt+\frac{1}{2}\int_{\R^N}\eta\big(div_x
L(x,T)\big)\,dx\\ \text{if}\quad div_x L(x,0)=v_0(x)\,.
\end{multline}
Then if for every  $L\in L^2\big(\R^N\times(0,T),\R^{k\times
N}\big)$ such that $\partial_t L, div_x L,\nabla_x div_x L\in L^2$,
we set $u(x,t):=div_x L(x,t)$ then $u\in
L^2\big(\R^N\times(0,T),\R^k\big)\cap
L^2\big(0,T;W^{-1,2}(\R^N,\R^k)\big)$, $u(x,0)=v_0(x)$ and $\nabla_x
H_{F,u}=\big(R+\partial_t L+F(div_x L)\big)$, where $div_x R=0$.
Thus in particular
$$\int_0^T\int_{\R^N}|\nabla_x
H_{F,u}|^2dxdt=\int_0^T\int_{\R^N}\Big|\partial_t L_u+F\big(div_x
L_u\big)\Big|^2dxd\leq \int_0^T\int_{\R^N}\Big|\partial_t
L+F\big(div_x L\big)\Big|^2dxdt\,.$$ Thus, as before, we obtain that
the $\Gamma-\liminf$ and $\Gamma-\limsup$ of $I_{\e,F}$ when
$u_\e\to u$ are the same as the $\Gamma-\liminf$ and
$\Gamma-\limsup$ of $\bar I_{\e,F}$ as $(div_x L_\e,\partial_t
L_\e)\to (u,-F(u))$.

 Applying the results of \cite{PI} and \cite{PII} we obtain the
following theorems about the upper and lower bounds
\begin{theorem}\label{prcnlkkmainthggg}
Let $F(u)\in C^1(\R^k,\R^{k\times N})$ and $\eta(u)\in C^3(\R^k,\R)$
be an entropy for system \er{equ24gggg}, which satisfies
$\eta(u)\geq 0$ and $\eta(0)=0$. Furthermore, let $u(x,t)\in
BV\big(\R^N\times(0,T),\R^k\big)\cap
L^\infty\big(0,T;L^2(\R^N,\R^k)\big)\cap L^\infty$ be such that
$u(x,t)$ is continuous in $[0,T]$ as a function of $t$ with the
values in $L^1(\R^N,\R^k)$ and satisfies the following Conservation
Law on the strip:
\begin{equation}\label{conslawgggggh}
\partial_t v(x,t)+div_x F\big(v(x,t)\big)=0\quad\quad\forall (x,t)\in\R^N\times(0,T)\,.
\end{equation}
Then there exists a sequence of functions $\big\{\bar
v_\e(x,t)\big\}_{\e>0}\subset
L^2_{loc}\big(\R^N\times(0,T),\R^{k\times N}\big)$ such that $\bar
u_\e(x,t):=\Div_x \bar v_\e(x,t)\in
L^2\big(0,T;H^1_0(\R^N,\R^k)\big)\cap
C\big([0,T];L^2(\R^N,\R^k)\big)\cap L^\infty$
and $\bar L_\e(x,t):=-\partial_t \bar v_\e(x,t)\in
L^2\big(\R^N\times(0,T),\R^{k\times N}\big)$;
$\bar u_\e\to u$ in $\bigcap_{q\geq 1}L^q
\big(\R^N\times(0,T);\R^k\big)$;  $\bar L_\e\to F(u)$ in
$L^2\big(\R^N\times(0,T),\R^{k\times N}\big)$;
$\partial_t \bar u_\e+\Div_x \bar L_\e\equiv 0$, $\bar
u_\e(x,0)=u(x,0)$ and
\begin{multline}\label{glavnkjjghjbmainmainthggg}
\lim_{\e\to
0^+}\Bigg\{\int_0^T\int_{\R^N}\Bigg(\e\Big|\nabla_x\big\{\nabla_u
\eta\big(\bar u_\e(x,t)\big)\big\}\Big|^2+\frac{1}{\e}\Big|\nabla_x
H_{F,\bar u_\e}(x,t)\Big|^2\Bigg)\,dxdt
+\int_{\R^N}\eta\big(\bar u_\e(x,T)\big)\,dx\Bigg\}\leq\\
\lim_{\e\to
0^+}\Bigg\{\int_0^T\int_{\R^N}\Bigg(\e\Big|\nabla_x\big\{\nabla_u
\eta\big(\bar u_\e(x,t)\big)\big\}\Big|^2+\frac{1}{\e}\Big| \bar
L_\e(x,t)-F\big(\bar
u_\e(x,t)\big)\Big|^2\Bigg)\,dxdt+\int_{\R^N}\eta\big(\bar
u_\e(x,T)\big)\,dx\Bigg\}\\= \lim_{\e\to
0^+}\Bigg\{\int_0^T\int_{\R^N}\Bigg(\e\bigg|\nabla_x
\Big\{\nabla_u\eta\big(\Div_x \bar
v_\e(x,t)\big)\Big\}\bigg|^2+\frac{1}{\e}\Big|
\partial_t \bar v_\e(x,t)+F\big(\Div_x
\bar v_\e(x,t)\big)\Big|^2\Bigg)\,dxdt
\\+\int_{\R^N}\eta\big(\Div_x\bar v_\e(x,T)\big)\,dx\Bigg\}
\leq\int_{J_{u}}\hat E_0\Big(u^+(x,t),
u^-(x,t),\vec\nu(x,t)\Big)\,\partial\mathcal{H}^N(x,t)+\int_{\R^N}\eta\big(u(x,T)\big)\,dx\,,
\end{multline}
where $H_{F,u}(x,t)$ satisfies \er{conssstgkfgjhgj},
\begin{multline}
\label{L2009limew03zeta71288888Contggiuuggyyyynew88789999vprop78899shtrihkkkllyhjyukjkkmmmklklklhhhhkkffgghhjjjkkkllkkhhhjjuhhiuijkffffudrrt7rddrggg}
\hat E_0\big(u^+, u^-,\vec\nu\big):=\inf\Bigg\{\int_{\bar I_{\vec
\nu}}\bigg(L\Big|\nabla_y
\Big\{\nabla_u\eta\big(\zeta(y,s)\big)\Big\}\Big|^2+\frac{1}{L}\Big|
\gamma(y,s)-F\big(\zeta(y,s)\big)\Big|^2\bigg)\,dyds:\\
L\in(0,+\infty),\,\zeta\in\mathcal{Z}^{(2)}\big(u^+,u^-,\vec\nu\big),\;
\gamma\in\mathcal{Z}^{(3)}\big(F(u^+),F(u^-),\vec\nu\big),\;\partial_s\zeta(y,s)+\Div_y\gamma(y,s)\equiv
0\Bigg\}=\\ \hat E_1\big(u^+,
u^-,\vec\nu\big):=\inf\Bigg\{\int_{\bar I_{\vec
\nu}}\bigg(L\Big|\nabla_y \Big\{\nabla_u\eta\big(\Div_y
\xi(y,s)\big)\Big\}\Big|^2+\frac{1}{L}\Big|
\partial_s\xi(y,s)+F\big(\Div_y
\xi(y,s)\big)\Big|^2\bigg)\,dyds:\\
L\in(0,+\infty),\,\xi\in\mathcal{Z}^{(1)}\big(u^+,u^-,\vec\nu\big)\Bigg\}
\,,
\end{multline}
with
\begin{multline}\label{L2009Ddef2hhhjjjj77788hhhkkkkllkjjjjkkkhhhhffggdddkkkgjhikhhhjjhhhhgbjhjngyuittrhvghggg}
\mathcal{Z}^{(1)}\big(u^+,u^-,\vec\nu\big):=\\ \Bigg\{\xi(y,s)\in
\mathcal{D}'(\R^N\times\R,\R^{k\times N}):\;\Div_y\xi(y,s)\in
C^1(\R^N\times\R,\R^{k}),\;\partial_s\xi(y,s)\in
C^1(\R^N\times\R,\R^{k\times N}),\\
\big(\Div_y\xi,-\partial_s\xi\big)(y,s)=\big(u^-,F(u^-)\big)\;\text{
if }\;y\cdot\vec\nu\leq-1/2,\;
\big(\Div_y\xi,-\partial_s\xi\big)(y,s)=\big(u^+,F(u^+)\big)\;\text{
if }\; y\cdot\vec\nu\geq 1/2\\ \text{ and
}\;\big(\Div_y\xi,-\partial_s\xi\big)\big((y,s)+\vec\nu_j\big)=\big(\Div_y\xi,-\partial_s\xi\big)(y,s)\;\;\forall
j=2,\ldots, (N+1)\Bigg\}\,,
\end{multline}
\begin{multline}\label{L2009Ddef2hhhjjjj77788hhhkkkkllkjjjjkkkhhhhffggdddkkkgjhikhhhjjhhhhgbjhjnjhhuuuuikfuytrtrghgyggg}
\mathcal{Z}^{(2)}\big(u^+,u^-,\vec\nu\big):= \Bigg\{\zeta(y,s)\in
C^1(\R^N\times\R,\R^{k}):\; \zeta(y,s)=u^-\;\text{ if
}\;y\cdot\vec\nu\leq-1/2,\\
\zeta(y,s)=u^+\;\text{ if }\; y\cdot\vec\nu\geq 1/2\; \text{ and
}\;\zeta\big((y,s)+\vec\nu_j\big)=\zeta(y,s)\;\;\forall j=2,\ldots,
(N+1)\Bigg\}\,,
\end{multline}
\begin{multline}\label{L2009Ddef2hhhjjjj77788hhhkkkkllkjjjjkkkhhhhffggdddkkkgjhikhhhjjhhhhgbjhjnijiuhujijikkkkdfurrrertyrtggg}
\mathcal{Z}^{(3)}\big(A,B,\vec\nu\big):= \Bigg\{\gamma(y,s)\in
C^1(\R^N\times\R,\R^{k\times N}):\;
\gamma(y,s)=B\;\text{ if }\;y\cdot\vec\nu\leq-1/2,\\
\gamma(y,s)=A\;\text{ if }\; y\cdot\vec\nu\geq 1/2\; \text{ and
}\;\gamma\big((y,s)+\vec\nu_j\big)=\gamma(y,s)\;\;\forall
j=2,\ldots, (N+1)\Bigg\}\,.
\end{multline}
Here $\bar I_{\vec \nu}:=\big\{y\in\R^{N+1}:\;|y\cdot
\vec\nu_j|<1/2\;\;\;\forall j=1,\ldots, (N+1)\big\}$ where
$\{\vec\nu_1,
\ldots,\vec\nu_N,
\vec\nu_{N+1}\}\subset\R^{N+1}$ is an orthonormal base in $\R^{N+1}$
such that $\vec\nu_1:=\vec \nu$.
\end{theorem}

\begin{theorem}\label{prcnlkkmainthhuiuiggg}
Let $F(u)\in C^1(\R^k,\R^{k\times N})$ and $\eta(u)\in C^3(\R^k,\R)$
be a convex entropy for the corresponding system
\er{conslawgggggh}, which satisfies $\eta(u)\geq 0$, $\eta(0)=0$ and
$\big|F(u)\big|\leq C\big(|u|+1\big)$ $\forall u\in\R^k$, for some
constant $C>0$. Furthermore, let $u(x,t)\in
BV\big(\R^N\times(0,T),\R^k\big)\cap
L^\infty\big(0,T;L^2(\R^N,\R^k)\big)\cap L^\infty$ be such that
$u(x,t)$ is continuous in $[0,T]$ as a function of $t$ with the
values in $L^1(\R^N,\R^k)$ and satisfies \er{conslawgggggh}.
Then for every sequence of functions
$u_\e(x,t)
\in L^2\big(0,T;H^1_0(\R^N,\R^k)\big)\cap
C\big([0,T];L^2(\R^N,\R^k)\big)\cap L^\infty$
and $L_\e(x,t)
\in L^2\big(\R^N\times(0,T),\R^{k\times N}\big)$ such that
$u_\e\to u$ in $L^2 \big(\R^N\times(0,T);\R^k\big)$, $L_\e\to F(u)$
in $L^2 \big(\R^N\times(0,T);\R^{k\times N}\big)$
and $\partial_t u_\e+\Div_x L_\e\equiv 0$, we have
\begin{multline}\label{glavnkjjghjbmainmainthgjhguyhgggg}
\lim_{\e\to
0^+}\Bigg\{\int_0^T\int_{\R^N}\Bigg(\e\Big|\nabla_x\big\{\nabla_u
\eta\big(u_\e(x,t)\big)\big\}\Big|^2+\frac{1}{\e}\Big|
L_\e(x,t)-F\big(u_\e(x,t)\big)\Big|^2\Bigg)\,dxdt+\int_{\R^N}\eta\big(u_\e(x,T)\big)\,dx\Bigg\}=\\
\lim_{\e\to 0^+}\Bigg\{\int_0^T\int_{\R^N}\Bigg(\e\bigg|\nabla_x
\Big\{\nabla_u\eta\big(\Div_x
v_\e(x,t)\big)\Big\}\bigg|^2+\frac{1}{\e}\Big|
\partial_t v_\e(x,t)+F\big(\Div_x
v_\e(x,t)\big)\Big|^2\Bigg)\,dxdt
\\+\int_{\R^N}\eta\big(\Div_x v_\e(x,T)\big)\,dx\Bigg\}\\ \geq\lim_{\e\to
0^+}\Bigg\{\int_0^T\int_{\R^N}\Bigg(\e\Big|\nabla_x\big\{\nabla_u
\eta\big(u_\e(x,t)\big)\big\}\Big|^2+\frac{1}{\e}\Big|\nabla_x
H_{F,u_\e}(x,t)\Big|^2\Bigg)\,dxdt+\int_{\R^N}\eta\big(u_\e(x,T)\big)\,dx\Bigg\}\\
\geq\int_{J_{u}}\hat I_0\Big(u^+(x,t),
u^-(x,t),\vec\nu(x,t)\Big)\,\partial\mathcal{H}^N(x,t)+\int_{\R^N}\eta\big(u(x,T)\big)\,dx\,,
\end{multline}
where $H_{F,u}(x,t)$ satisfies \er{conssstgkfgjhgj}, $v_\e(x,t)\in
L^2_{loc}\big(\R^N\times(0,T),\R^{k\times N}\big)$ is such that
$u_\e(x,t)=\Div_x v_\e(x,t)$ and $L_\e(x,t)=-\partial_t v_\e(x,t)$,
\begin{multline}
\label{L2009limew03zeta71288888Contggiuuggyyyynew88789999vprop78899shtrihkkkllyhjyukjkkmmmklklklhhhhkkffgghhjjjkkkllkkhhhjjuhhiuijkffffudrrt7rddrjkgggyhggg}
\hat I_0\big(u^+, u^-,\vec\nu\big):=\inf\Bigg\{\liminf_{\e\to
0^+}\int_{\bar I_{\vec \nu}}\bigg(\e\Big|\nabla_y
\Big\{\nabla_u\eta\big(\zeta_\e(y,s)\big)\Big\}\Big|^2+\frac{1}{\e}\Big|
\gamma_\e(y,s)-F\big(\zeta_\e(y,s)\big)\Big|^2\bigg)\,dyds:\\
\zeta_\e\to\chi(y,u^+,u^-,\vec\nu)\;\text{in}\;L^2\big(I_{\vec\nu},\R^{k}\big),\;
\gamma_\e\to\hat\chi(y,u^+,u^-,\vec\nu,F)\;\text{in}\;L^2\big(I_{\vec\nu},\R^{k\times
N}\big)\\
\text{and}\;\partial_s\zeta_\e(y,s)+\Div_y\gamma_\e(y,s)\equiv
0\Bigg\}\\
=\inf\Bigg\{\liminf_{\e\to 0^+}\int_{\bar I_{\vec
\nu}}\bigg(\e\Big|\nabla_y \Big\{\nabla_u\eta\big(\Div_y
\xi_\e(y,s)\big)\Big\}\Big|^2+\frac{1}{\e}\Big|
\partial_s\xi_\e(y,s)+F\big(\Div_y
\xi_\e(y,s)\big)\Big|^2\bigg)\,dyds:\\
\Div_x\xi_\e\to\chi(y,u^+,u^-,\vec\nu)\;\text{in}\;L^2\big(I_{\vec\nu},\R^{k}\big)\;-\partial_s\xi_\e\to\hat\chi(y,u^+,u^-,\vec\nu,F)
\;\text{in}\;L^2\big(I_{\vec\nu},\R^{k\times N}\big)\Bigg\} \,,
\end{multline}
\begin{equation*} \chi(y,u^+,u^-,\vec\nu):=\begin{cases}
u^+\quad\text{if}\;\;y\cdot\vec\nu>0\,,\\
u^-\quad\text{if}\;\;y\cdot\vec\nu<0\,,
\end{cases}
\end{equation*}
and
\begin{equation*} \hat\chi(y,u^+,u^-,\vec\nu,F):=\begin{cases}
F(u^+)\quad\text{if}\;\;y\cdot\vec\nu>0\,,\\
F(u^-)\quad\text{if}\;\;y\cdot\vec\nu<0\,.
\end{cases}
\end{equation*}
\end{theorem}
\begin{remark}\label{vyuguigiugbuikkk2888}
In what follows we use some special notations and apply some basic
theorems about $BV$ functions. For the convenience of the reader we
put these notations and theorems in Appendix.
\end{remark}

\section{The non-local problems related to Micromagnetics}
\begin{lemma}\label{gbhfhgdfdghfddhg}
Let $\O\subset\R^N$ be an open set (possibly unbounded) with locally
Lipschitz's boundary and let
\begin{equation}\label{vhfhgjkhjjkhjhgkgkjggjkgjk}
L_{k,N}(\O):=\Big\{L(x)\in L^2(\O,\R^{k\times N}),\;\;div
\,L(x)\equiv 0\Big\}\,.
\end{equation}
Given $M(x)\in L^2(\O,\R^{k\times N})$ consider the functional
$J_{M}(\big(L(x)\big):L_{k,N}(\O)\to\R$ defined by
\begin{equation}\label{defH12sldfgvrftghbjhghvnvjkjklkhhhgffgfg}
J_{M}\big(L(x)\big):=\int_{\O}\Big|L(x)+M(x)\Big|^2dx\,,
\end{equation}
and consider the variational problem
\begin{equation}\label{defH12sldfgvrftghbjhghvnv}
J_0:=\inf\bigg\{J_{M}\big(L(x)\big)
:\;\;L\in L_{k,N}(\O)\bigg\}\,.
\end{equation}
Then there exists a unique minimizer $L_0(x)\in L_{k,N}(\O)$ to
\er{defH12sldfgvrftghbjhghvnv}, i.e. $J_{M}\big(L_0(x)\big)=J_0$.
Moreover we have $L_0(x)+M(x)\equiv\nabla H_M(x)$ where $H_M(x)\in
W^{1,2}_{loc}(\R^N,\R^k)$ is a function which satisfies
\begin{equation}\label{defH12slvvvbhjhjhjhj}
\begin{cases}\Delta_x H_M(x)=\Div_x M(x)\quad\text{in }\;\R^N\,,\\
H_M(x)=0\quad\text{on}\;\;\partial\O\,,
\\
\nabla_x H_M(x)\in L^2(\O,\R^{k\times N})\,.
\end{cases}
\end{equation}
\end{lemma}
\begin{proof}
Clearly $L_{k,N}(\O)$ is a closed subspace of the Hilbert space
$L^2(\O,\R^{k\times N})$. Therefore, clearly there exists a closed
subspace $V\subset L^2(\O,\R^{k\times N})$, such that $V$ is an
orthogonal complement of $L_{k,N}(\O)$ in $L^2(\O,\R^{k\times N})$,
i.e.
\begin{equation}\label{yfjhgkhklhllhk}
\int_\O u(x): v(x)\,dx=0\,,\quad\forall u(x)\in
L_{k,N}(\O),\;\,\forall v(x)\in V\,,
\end{equation}
and for every $h(x)\in L^2(\O,\R^{k\times N})$ there exist uniquely
defined $u_h(x)\in L_{k,N}(\O)$ and $v_h(x)\in V$ such that
$h(x)\equiv u_h(x)+v_h(x)$. Thus, in particular, there exist
uniquely defined $P(x)\in L_{k,N}(\O)$ and $Q(x)\in V$ such that
$-M(x)\equiv P(x)+Q(x)$. On the other hand for arbitrary $L(x)\in
L_{k,N}(\O)$, using \er{yfjhgkhklhllhk} we have
\begin{multline}\label{vhjffjuguggug}
J_{M}\big(L(x)\big):=\int_{\O}\Big|L(x)+M(x)\Big|^2dx=\int_{\O}\Big|\big(P(x)-L(x)\big)+Q(x)\Big|^2dx
=\\
\int_{\O}\Big|P(x)-L(x)\Big|^2dx+\int_{\O}\Big|Q(x)\Big|^2dx+2\int_\O
\big(P(x)-L(x)\big):Q(x)dx=\\
\int_{\O}\Big|P(x)-L(x)\Big|^2dx+\int_{\O}\Big|P(x)+M(x)\Big|^2dx=
J_{M}\big(P(x)\big)+\int_{\O}\Big|P(x)-L(x)\Big|^2dx\,.
\end{multline}
Thus $L_0(x):=P(x)$ is unique minimizer to
\er{defH12sldfgvrftghbjhghvnv}. Moreover, since $P(x)\in
L_{k,N}(\O)$ and $Q(x)\in V$ we have
\begin{equation}\label{yfjhgkhklhllhkhjhjhj}
\int_\O Q(x): L(x)\,dx=0\,,\quad\forall L(x)\in L_{k,N}(\O)\,.
\end{equation}
Thus since $L_0(x)+M(x)\equiv-Q(x)$ we obtain
\begin{equation}\label{yfjhgkhklhllhkhjhjhjjhjhjjjjk}
\int_\O \big(L_0(x)+M(x)\big): L(x)\,dx=0\,,\quad\forall L(x)\in
L_{k,N}(\O)\,.
\end{equation}
In particular,
\begin{equation}\label{yfjhgkhklhllhkhjhjhjjhjhjjjjkhlhhj}
\int_\O \big(L_0(x)+M(x)\big): \delta(x)\,dx=0\,,\quad\forall
\delta(x)\in C^1_c(\ov\O,\R^{k\times N})\;\;\text{such
that}\;\Div\delta(x)\equiv 0\,.
\end{equation}
Thus clearly there exists a function $H_M(x)\in
W^{1,2}_{loc}(\O,\R^k)$ such that $H_M(x)\in W^{1,2}(G,\R^k)$ for
every bounded open subset $G\subset\O$ and $\nabla_x H_M(x)\equiv
L_0(x)+M(x)$ on $\O$. Thus in particular
$$\Delta_x H_M(x)=\Div_x\big(\nabla_x H_M\big)\equiv
\Div L_0(x)+\Div M(x)=\Div_x M(x)\,.$$ Moreover obviously $\nabla_x
H_M(x)\in L^2(\O,\R^{k\times N})$. Finally by
\er{yfjhgkhklhllhkhjhjhjjhjhjjjjkhlhhj} we have
\begin{multline}\label{yfjhgkhklhllhkhjhjhjjhjhjjjjkhlhhjhhjhghgggh}
\int_{\partial\O}H_M(x)\cdot\Big\{\delta(x)\cdot \vec
n(x)\Big\}\,d\mathcal{H}^{N-1}(x)=\int_\O \Big\{\nabla_x H_M(x):
\delta(x)+H_M(x)\cdot \Div\delta_x(x)\Big\}\,dx=0\\ \forall
\delta(x)\in C^1_c(\ov\O,\R^{k\times N})\;\;\text{such
that}\;\Div\delta_x(x)\equiv 0\,.
\end{multline}
Thus $H_M(x)=0$ on $\partial\O$ in the sense of trace. This
completes the proof.
\end{proof}
\begin{lemma}\label{vfgcfgjfvhgkjkj}
Let $\O\subset\R^N$ be an open set with locally Lipschitz's boundary
and $G(\psi,x):\R^m\times\O\to\R^{k\times N}$ be a measurable
function, continuous by the first argument $\psi$, which satisfies
\begin{equation}\label{ghhghjghjhbfggffg}
\big|G(\psi,x)\big|\leq
C_0|\psi|^{p/2}+h_0(x)\quad\forall\psi\in\R^m,\;x\in\O\,,
\end{equation}
for some constant $C_0>0$, $p\geq 1$ and $h_0\in L^2(\O,\R)$.
Furthermore, for every $\e>0$ consider the functional
$E_\e\big(\psi(x)\big):L^p(\O,\R^m)\to[0,+\infty)\cup\{+\infty\}$
which (possibly) can attain the infinite values. Next for every
$\e>0$ consider the functional
$P_\e\big(\psi(x)\big):L^p(\O,\R^m)\to[0,+\infty)\cup\{+\infty\}$,
defined by
\begin{equation}\label{fhfhjffhjhf}
P_\e\big(\psi(x)\big):=E_\e\big(\psi(x)\big)+\frac{1}{\delta_\e}\int_{\R^N}\Big|\nabla_x
V_{G,\psi}(x)\Big|^2 dx\,,
\end{equation}
where $\delta_\e>0$ satisfies $\lim_{\e\to 0^+}\delta_\e=0$ and
given $\psi(x)\in L^p(\O,\R^m)$, $V_{G,\psi}(x):\R^N\to\R^k$ is
defined by
\begin{equation}\label{defH12slvvvb}
\begin{cases}\Delta_x V_{G,\psi}(x)=\Div_x\Big\{\chi_\O(x)\cdot G\big(\psi(x),x\big)\Big\}\quad\text{in }\R^N\,,\\
\nabla_x V_{G,\psi}(x)\in L^2(\R^N,\R^{k\times N})\,,
\end{cases}
\end{equation}
with $\chi_\O(x):=1$ if $x\in\O$ and $\chi_\O(x):=0$ if
$x\in\R^N\setminus\O$. Furthermore, for every $\e>0$ consider the
functional $Q_\e\big(\psi(x),L(x)\big):L^p(\O,\R^m)\times
L^2(\R^N,\R^{k\times N}) \to[0,+\infty)\cup\{+\infty\}$ defined by
\begin{equation}\label{fhfhjffhjhfgyukuoijkkj}
Q_\e\big(\psi(x),L(x)\big):=\begin{cases}E_\e\big(\psi(x)\big)+\frac{1}{\delta_\e}\int_{\R^N}\Big|L(x)+\chi_\O(x)\cdot
G\big(\psi(x),x\big)\Big|^2 dx\quad\quad\text{if}\quad \Div_x
L(x)=0\,,\\+\infty\quad\quad\quad\text{otherwise}\,.
\end{cases}
\end{equation}
Next for every $\varphi(x)\in L^p(\O,\R^m)$, such that
$\Div_x\big\{\chi_\O(x)\cdot G\big(\varphi(x),x\big)\big\}=0$ in
$\R^N$ set
\begin{equation}\label{hfghjfhfhfh}
\begin{split}
\underline P(\varphi):=\inf\Bigg\{\liminf_{\e\to
0^+}P_\e\big(\psi_\e(x)\big):\;\;\psi_\e(x)\to\varphi(x)\;\;\text{in}\;\;L^p(\O,\R^m)\Bigg\}\,,\\
\ov P(\varphi):=\inf\Bigg\{\limsup_{\e\to
0^+}P_\e\big(\psi_\e(x)\big):\;\;\psi_\e(x)\to\varphi(x)\;\;\text{in}\;\;L^p(\O,\R^m)\Bigg\}\,,
\end{split}
\end{equation}
and
\begin{equation}\label{hfghjfhfhfhj;k;jkhkhhk}
\begin{split}
\underline Q(\varphi):=\inf\Bigg\{\liminf_{\e\to
0^+}Q_\e\big(\psi_\e(x),L_\e(x)\big):\;\;\psi_\e(x)\to\varphi(x)\;\;\text{in}\;\;L^p(\O,\R^m)\\
\text{and}\;\;L_\e(x)\to\Big(-\chi_\O(x)\cdot
G\big(\varphi(x),x\big)\Big)\;\;\text{in}\;\;L^2(\R^N,\R^{k\times
N})\Bigg\}\,,\\ \ov Q(\varphi):=\inf\Bigg\{\limsup_{\e\to
0^+}Q_\e\big(\psi_\e(x),L_\e(x)\big):
\psi_\e(x)\to\varphi(x)\;\;\text{in}\;\;L^p(\O,\R^m)\\
\text{and}\;\;L_\e(x)\to\Big(-\chi_\O(x)\cdot
G\big(\varphi(x),x\big)\Big)\;\;\text{in}\;\;L^2(\R^N,\R^{k\times
N})\Bigg\}\,,
\end{split}
\end{equation}
Then we have the following equalities
\begin{equation}\label{vhfhvffgjfghkjjjjjkjkhk}
\underline P(\varphi)=\underline Q(\varphi)\quad\text{and}\quad \ov
P(\varphi)=\ov Q(\varphi)
\,.
\end{equation}
\end{lemma}
\begin{proof}
Fix some $\varphi(x)\in L^p(\O,\R^m)$ such that
$\Div_x\big\{\chi_\O(x)\cdot G\big(\varphi(x),x\big)\big\}=0$ in
$\R^N$. Then by \er{ghhghjghjhbfggffg} we have
$G\big(\varphi(x),x\big)\in L^2(\O,\R^{k\times N})$ and thus
\begin{equation}\label{jdfhhghfhg}
-\chi_\O(x)\cdot G\big(\varphi(x),x\big)\in L^2(\R^N,\R^{k\times
N})\,.
\end{equation}
Next fix some sequence $\{\psi_\e(x)\}\subset L^p(\O,\R^m)$ such
that $\psi_\e(x)\to\varphi(x)$ in $L^p(\O,\R^m)$ as $\e\to 0^+$.
Then by \er{ghhghjghjhbfggffg} we have
\begin{equation}\label{jdfhhghfhgfhkgfhkdhbnbjlbhbjjnm}
-\chi_\O(x)\cdot G\big(\psi_\e(x),x\big)\to-\chi_\O(x)\cdot
G\big(\varphi(x),x\big)\quad\text{in}\;\; L^2(\R^N,\R^{k\times
N})\,.
\end{equation}
On the other hand by Lemma \ref{gbhfhgdfdghfddhg} together with
\er{fhfhjffhjhf}, \er{defH12slvvvb} and \er{fhfhjffhjhfgyukuoijkkj}
clearly we have
\begin{equation}\label{fighfghhffjhhhj}
Q_\e\big(\psi_\e(x),L(x)\big)\geq
P_\e\big(\psi_\e(x)\big)\quad\forall L(x)\in L^2(\R^N,\R^{k\times
N})\,.
\end{equation}
Moreover, if we set
\begin{equation}\label{fighfghhffjhhhjhhjhjjjkkj}
\hat L_{\e}(x):=\nabla_x V_{G,\psi}(x)-\chi_\O(x)\cdot
G\big(\psi_\e(x),x\big)\quad\forall x\in\R^N
\end{equation}
then $\hat L_{\e}(x)\in L^2(\R^N,\R^{k\times N})$, $\Div \hat
L_\e(x)\equiv 0$ and
\begin{equation}\label{fighfghhffjhhhjhfjdgdgdgdh}
Q_\e\big(\psi_\e(x),\hat L_\e(x)\big)= P_\e\big(\psi_\e(x)\big)\,.
\end{equation}
In particular since $\delta_\e\to 0$ by
\er{fighfghhffjhhhjhfjdgdgdgdh}, \er{fhfhjffhjhfgyukuoijkkj} and
\er{jdfhhghfhgfhkgfhkdhbnbjlbhbjjnm}, for arbitrary sequence
$\e_n\to 0^+$ as $n\to+\infty$, we must have
\begin{equation}\label{fighfghhffjhhhjhfjdgdgdgdhhjkhjlkjlnghjg}
\hat L_{\e_n}(x)\to -\chi_\O(x)\cdot
G\big(\varphi(x),x\big)\;\;\text{in}\;\;L^2(\R^N,\R^{k\times
N})\quad\quad\text{if}\quad\limsup_{n\to+\infty}P_{\e_n}\big(\psi_{\e_n}(x)\big)<+\infty\,,
\end{equation}
Moreover, by \er{fighfghhffjhhhjhfjdgdgdgdh}, in this case we have
\begin{equation}\label{fighfghhffjhhhjhfjdgdgdgdhhjkhjlkjlnghjgjyjyjujjmhkjhjkh}
\limsup_{n\to+\infty}Q_{\e_n}\big(\psi_{\e_n}(x),\hat
L_{\e_n}(x)\big)=\limsup_{n\to+\infty}P_{\e_n}\big(\psi_{\e_n}(x)\big)<+\infty\,.
\end{equation}
On the other hand by \er{fighfghhffjhhhj} for every sequence $\hat
L_{\e_n}(x)\to -\chi_\O(x)\cdot G\big(\varphi(x),x\big)$ in
$L^2(\R^N,\R^{k\times N})$ we obviously have
\begin{equation}\label{fighfghhffjhhhjhfjdgdgdgdhhjkhjlkjlnghjghjhjhjnjhhj}
\limsup_{n\to+\infty}Q_{\e_n}\big(\psi_{\e_n}(x),L_{\e_n}(x)\big)=+\infty\quad\text{if}\quad\limsup_{n\to+\infty}P_{\e_n}\big(\psi_{\e_n}(x)\big)=+\infty\,.
\end{equation}
Therefore, by \er{fighfghhffjhhhjhfjdgdgdgdhhjkhjlkjlnghjg},
\er{fighfghhffjhhhjhfjdgdgdgdhhjkhjlkjlnghjgjyjyjujjmhkjhjkh} and
\er{fighfghhffjhhhjhfjdgdgdgdhhjkhjlkjlnghjghjhjhjnjhhj} in any case
there exists a sequence $\tilde L_{\e_n}(x)\to -\chi_\O(x)\cdot
G\big(\varphi(x),x\big)$ in $L^2(\R^N,\R^{k\times N})$ so that
\begin{equation}\label{fighfghhffjhhhjhfjdgdgdgdhhjkhjlkjlnghjgjyjyjujjmhkjhjkhgjuhik}
\limsup_{n\to+\infty}Q_{\e_n}\big(\psi_{\e_n}(x),\tilde
L_{\e_n}(x)\big)=\limsup_{n\to+\infty}P_{\e_n}\big(\psi_{\e_n}(x)\big)\,.
\end{equation}
Thus since $\e_n\to 0^+$ $\psi_\e\to\varphi$ were chosen arbitrary,
by
\er{fighfghhffjhhhjhfjdgdgdgdhhjkhjlkjlnghjgjyjyjujjmhkjhjkhgjuhik}
we deduce
\begin{equation}\label{vhfhvffgjfghkjjjjjkjkhkhhgjghgh}
\underline P(\varphi)\geq\underline Q(\varphi)\quad\text{and}\quad
\ov P(\varphi)\geq \ov Q(\varphi)\quad\quad\forall \varphi\in
L^p(\O,\R^m)\,.
\end{equation}
On the other hand, by \er{fighfghhffjhhhj}, clearly
\begin{equation}\label{vhfhvffgjfghkjjjjjkjkhkhhgjghghmkjjk}
\underline P(\varphi)\leq\underline Q(\varphi)\quad\text{and}\quad
\ov P(\varphi)\leq \ov Q(\varphi)\quad\quad\forall \varphi\in
L^p(\O,\R^m)\,.
\end{equation}
This completes the proof.
\end{proof}
Next plugging Lemma \ref{gbhfhgdfdghfddhg} into Theorem
4.2 in \cite{PI} we deduce the following upper bound result for
problem with a non-local term.
\begin{theorem}\label{dehgfrygfrgygenbgggggggggggggkgkgthtjtfjkjkkhhk}
Let $\O\subset\R^N$ be an open set with locally Lipschitz's
boundary. Furthermore, let $G\in C^1(\R^m\times\R^q,\R^{k\times N})$
and $F\in C^1\big(\R^{m\times N^n}
\times\ldots\times\R^{m\times N}\times \R^m\times\R^q,\R\big)$, be
such that $F\geq 0$. Next let $\varphi\in BV(\R^N,\R^{m})\cap
L^\infty$ and $f\in BV_{loc}(\R^N,\R^q)\cap L^\infty$ be such that
$\|D \varphi\|(\partial\Omega)=0$, $F\big(0,\ldots,0,
\varphi(x),f(x)\big)=0$ for a.e. $x\in\O$, $\Div_x
G\big(\varphi(x),f(x)\big)=0$ in $\O$ and
$G\big(\varphi(x),f(x)\big)\cdot\vec n(x)=0$ on $\partial\O$.
Then for every $\delta>0$ there exists a sequence
$\{\psi_\e\}_{\e>0}\subset C^\infty(\R^N,\R^m)$ such that
$\int_\O\psi_\e(x)dx=\int_\O\varphi(x)dx$, $\lim_{\e\to
0^+}\psi_\e=\varphi$ in $L^p$ and $\lim_{\e\to
0^+}\e^j\nabla^j\psi_\e=0$ in $L^p$ for every $p\geq 1$ and every
$j\in\{1,\ldots,n\}$, $\{\e^n\nabla^n\psi_\e\}_{\e>0},\ldots,
\{\e\nabla\psi_\e\}_{\e>0}$ and $\{\psi_\e\}_{\e>0}$ are a bounded
in $L^\infty$ sequences, and we have
\begin{multline}\label{a1a2a3a4a5a6a7s8hhjhjjhjjjjjjkkkkgenhjhhhhjtjurtgfhfhfjfjfjhjjjnjhjjj}
\limsup_{\e\to 0^+}\Bigg\{\frac{1}{\e}\int_\O F\Big(\e^n\nabla^n
\psi_\e(x),\ldots,\e\nabla
\psi_\e(x),\psi_\e(x),f(x)\Big)dx+\frac{1}{\e}\int_{\R^N}|\nabla
V_{G,\psi_\e}|^2dx\Bigg\} \\ \leq\int_{\O\cap J_\varphi}\hat
E_{per}\Big(\varphi^+(x),\varphi^-(x),\vec \nu(x),x\Big)d \mathcal
H^{N-1}(x)+\delta\,,
\end{multline}
where $V_{G,\psi}:\R^N\to\R^k$ is defined by
\begin{equation}\label{defH12slbgghfjfg}
\begin{cases}\Delta_x V_{G,\psi}(x)=\Div_x\Big\{\chi_\O(x) G\big(\psi(x),f(x)\big)\Big\}\quad\text{in }\R^N\,,\\
\nabla V_{G,\psi}\in L^2(\R^N,\R^{k\times N})\,,
\end{cases}
\end{equation}
\begin{multline}\label{L2009hhffff12kkkhjhjghghgvgvggcjhggghtgjuthjhhhjj}
\hat E_{per}\big(\varphi^+,\varphi^-,\vec \nu,x\big)\;:=\;
\inf\Bigg\{\int_{I_{\vec \nu}}\frac{1}{L} F\Big(L^n\nabla^n
\zeta(y),\ldots,L\nabla \zeta(y),\zeta(y), \sigma_{f,x}(y)
\Big)dy+
\int_{I_{\vec \nu}}\frac{1}{L} \big|\nabla
H_{G,\zeta,x,\vec\nu}(y)\big|^2 dy:
\\
L\in(0,+\infty)\,,\;\zeta\in
\mathcal{S}(\varphi^+,\varphi^-,I_{\vec\nu})\Bigg\}\,,
\end{multline}
where $H_{G,\zeta,x,\vec\nu}\in W^{2,2}_{loc}(\R^N,\R^k)$ satisfies
\begin{equation}\label{defH12slbgghfjfghjhjhjhjjkjkjhhj}
\begin{cases}\Delta_y H_{G,\zeta,x,\vec\nu}(y)=\Div_y G\big(\zeta(y),\sigma_{f,x}(y)\big)\quad\text{in }I_{\vec\nu}\,,\\
H_{G,\zeta,x,\vec\nu}(y+\vec\nu_j)=H_{G,\zeta,x,\vec\nu}(y)\quad\forall
y\in\R^N\;\;\text{such that}\;\;|y\cdot\vec\nu|<1/2\,,\\
\frac{\partial}{\partial\vec\nu}H_{G,\zeta,x,\vec\nu}(y)=0\quad\forall
y\in\R^N\;\;\text{such that}\;\;|y\cdot\vec\nu|=1/2\,,
\end{cases}
\end{equation}
with
\begin{equation}\label{bhggkjghkhjkghgjggkgj}
\sigma_{f,x}(y):=
\begin{cases}
f^+(x)\quad\text{if}\;\;y\cdot\vec\nu>0\,,\\
f^-(x)\quad\text{if}\;\;y\cdot\vec\nu<0\,,
\end{cases}
\end{equation}
and
\begin{multline}\label{L2009Ddef2hhhjjjj77788hhhkkkkllkjjjjkkkhhhhffggdddkkkgjhikhhhjjddddhdkgkkkhgghjhjhjkjtjytrjghjghjjjg}
\mathcal{S}\big(\varphi^+,\varphi^-,I_{\vec \nu}\big):=
\bigg\{\zeta\in C^n(\R^N,\R^m):\;\;\zeta(y)=\varphi^-\;\text{ if }\;y\cdot\vec \nu\leq-1/2,\\
\zeta(y)=\varphi^+\;\text{ if }\; y\cdot\vec \nu\geq 1/2\;\text{ and
}\;\zeta\big(y+\vec \nu_j\big)=\zeta(y)\;\;\forall j=2,\ldots,
N\bigg\}.
\end{multline}
Here $I_{\vec \nu}:=\Big\{y\in\R^N:\;|y\cdot
\vec\nu_j|<1/2\;\;\forall j=1,\ldots, N\Big\}$,
where $\{\vec\nu_1,
\ldots,\vec\nu_N\}\subset\R^N$ is an
orthonormal base in $\R^N$ such that $\vec\nu_1:=\vec \nu$.
\end{theorem}
\begin{proof}
Since $\Div_x G\big(\varphi(x),f(x)\big)=0$ in $\O$ and
$G\big(\varphi(x),f(x)\big)\cdot\vec n(x)=0$ on $\partial\O$ we
easily deduce that $$\Div_x \{\chi_\O(x)
G\big(\varphi(x),f(x)\big)\}=0\quad\text{in}\;\;\R^N$$ in the sense
of distribution. Then by Theorem 4.2
in \cite{PI} we deduce that for every $\delta>0$ there exist
sequences $\{\psi_\e\}_{\e>0}\subset C^\infty(\R^N,\R^m)$ and
$\{L_\e\}_{\e>0}\subset C^\infty(\R^N,\R^{k\times N})$ such that
$\int_\O\psi_\e(x)dx=\int_\O\varphi(x)dx$, $\lim_{\e\to
0^+}\psi_\e=\varphi$ in $L^p$ and $\lim_{\e\to
0^+}\e^j\nabla^j\psi_\e=0$ in $L^p$ for every $p\geq 1$ and every
$j\in\{1,\ldots,n\}$, $\{\e^n\nabla^n\psi_\e\}_{\e>0},\ldots,
\{\e\nabla\psi_\e\}_{\e>0}$, $\{\psi_\e\}_{\e>0}$ are bounded in
$L^\infty$ sequences, $\Div L_\e\equiv 0$ in $\R^N$, $L_\e\to
\chi_\O G\big(\varphi,f\big)$ in $L^2$ and we have
\begin{multline}\label{a1a2a3a4a5a6a7s8hhjhjjhjjjjjjkkkkgenhjhhhhjtjurtgfhfhfjfjfjhjjjnjhjjjkkkklkkk}
\limsup_{\e\to 0^+}\Bigg\{\frac{1}{\e}\int_{\R^N} F\Big(\e^n\nabla^n
\psi_\e(x),\ldots,\e\nabla
\psi_\e(x),\psi_\e(x),f(x)\Big)\chi_\O(x)dx+\frac{1}{\e}\int_{\R^N}\Big|L_\e-\chi_\O
G\big(\psi_\e,f\big)\Big|^2dx\Bigg\}\\ \leq \int_{\O\cap J_\varphi}
\tilde E_{per}\Big(\varphi^+(x),\varphi^-(x),\vec \nu(x),x\Big)d
\mathcal H^{N-1}(x)+\delta\,,
\end{multline}
where
\begin{multline}\label{L2009hhffff12kkkhjhjghghgvgvggcjhggghtgjuthjhhhjjgjjhgj}
\tilde E_{per}\big(\varphi^+,\varphi^-,\vec \nu,x\big):=
\inf\Bigg\{\int_{I_{\vec \nu}}\frac{1}{L} F\Big(L^n\nabla^n
\zeta(y),\ldots,L\nabla \zeta(y),\zeta(y),\sigma_{f,x}(y)
\Big)dy+\\
\int_{I_{\vec \nu}}\frac{1}{L}
\Big|\xi(y)-G\big(\zeta(y),\sigma_{f,x}(y)\big)\Big|^2dy: \;\;\;
L\in(0,+\infty)\,,\;\zeta\in
\mathcal{S}(\varphi^+,\varphi^-,I_{\vec\nu})\,,\;\xi\in\mathcal{S}_0(\varphi^+,f^+,\varphi^-,f^-,I_{\vec\nu})\Bigg\}\,,
\end{multline}
with
\begin{multline}\label{L2009Ddef2hhhjjjj77788hhhkkkkllkjjjjkkkhhhhffggdddkkkgjhikhhhjjddddhdkgkkkhgghjhjhjkjtjytrjghjghjjjgvvhjvhjjvh}
\mathcal{S}\big(\varphi^+,\varphi^-,I_{\vec \nu}\big):=
\bigg\{\zeta\in C^n(\R^N,\R^m):\;\;\zeta(y)=\varphi^-\;\text{ if }\;y\cdot\vec \nu\leq-1/2,\\
\zeta(y)=\varphi^+\;\text{ if }\; y\cdot\vec \nu\geq 1/2\;\text{ and
}\;\zeta\big(y+\vec \nu_j\big)=\zeta(y)\;\;\forall j=2,\ldots,
N\bigg\}\quad\text{and}\\
\mathcal{S}_0\big(\varphi^+,f^+,\varphi^-,f^-,I_{\vec \nu}\big):=
\bigg\{\xi\in C^n(\R^N,\R^{k\times N}):\;\;\Div_y\xi(y)=0,\;\;\xi(y)=G(\varphi^-,f^-)\;\text{ if }\;y\cdot\vec \nu\leq-1/2,\\
\xi(y)=G(\varphi^+,f^+)\;\text{ if }\; y\cdot\vec \nu\geq
1/2\;\text{ and }\;\xi\big(y+\vec \nu_j\big)=\xi(y)\;\;\forall
j=2,\ldots, N\bigg\}
\end{multline}
Thus using Lemma \ref{gbhfhgdfdghfddhg}, by
\er{a1a2a3a4a5a6a7s8hhjhjjhjjjjjjkkkkgenhjhhhhjtjurtgfhfhfjfjfjhjjjnjhjjjkkkklkkk}
we deduce
\begin{multline}\label{a1a2a3a4a5a6a7s8hhjhjjhjjjjjjkkkkgenhjhhhhjtjurtgfhfhfjfjfjhjjjnjhjjjjhgjgjg}
\limsup_{\e\to 0^+}\Bigg\{\frac{1}{\e}\int_\O F\Big(\e^n\nabla^n
\psi_\e(x),\ldots,\e\nabla
\psi_\e(x),\psi_\e(x),f(x)\Big)dx+\frac{1}{\e}\int_{\R^N}|\nabla
V_{G,\psi_\e}|^2dx\Bigg\} \\ \leq\int_{\O\cap J_\varphi}\tilde
E_{per}\Big(\varphi^+(x),\varphi^-(x),\vec \nu(x),x\Big)d \mathcal
H^{N-1}(x)+\delta\,,
\end{multline}
where $V_{G,\psi}:\R^N\to\R^k$ is defined by \er{defH12slbgghfjfg}.
Therefore, in order to complete the proof of the Theorem it is
sufficient to prove that we always have
\begin{equation}\label{L2009hhffff12kkkhjhjghghgvgvggcjhggghtgjuthjhhhjjgjjhgjgjkgkgk}
\tilde E_{per}\Big(\varphi^+,\varphi^-,\vec \nu,x\Big)\;=\;\hat
E_{per}\Big(\varphi^+,\varphi^-,\vec \nu,x\Big)\,,
\end{equation}
(see the definitions of the corresponding quantities in
\er{L2009hhffff12kkkhjhjghghgvgvggcjhggghtgjuthjhhhjjgjjhgj} and
\er{L2009hhffff12kkkhjhjghghgvgvggcjhggghtgjuthjhhhjj}). So fix some
$\zeta\in \mathcal{S}(\varphi^+,\varphi^-,I_{\vec\nu})$ and $L>0$.
Then it is sufficient to prove that
\begin{multline}\label{L2009hhffff12kkkhjhjghghgvgvggcjhggghtgjuthjhhhjjgjjhgjgnnmgmn}
\inf\Bigg\{\int_{I_{\vec \nu}}
\Big|\xi(y)-G\big(\zeta(y),\sigma_{f,x}(y)\big)\Big|^2\,dy:
\;\xi\in\mathcal{S}_0(\varphi^+,f^+,\varphi^-,f^-,I_{\vec\nu})\Bigg\}
=\int_{I_{\vec \nu}}
\big|\nabla
H_{G,\zeta,x,\vec\nu}(y)\big|^2\,dy\,,
\end{multline}
where $H_{G,\zeta,x,\vec\nu}\in W^{2,2}_{loc}(\R^N,\R^k)$ satisfies
\er{defH12slbgghfjfghjhjhjhjjkjkjhhj}. Indeed set
\begin{multline}\label{L2009Ddef2hhhjjjj77788hhhkkkkllkjjjjkkkhhhhffggdddkkkgjhikhhhjjddddhdkgkkkhgghjhjhjkjtjytrjghjghjjjgvvhjvhjjvhhnjhdjhghggh}
\mathcal{T}\big(\varphi^+,f^+,\varphi^-,f^-,I_{\vec \nu}\big):=
\Bigg\{\xi\in L^2_{loc}(\R^N,\R^{k\times N}):\;\Div_y\xi(y)=0,\;\xi(y)=G(\varphi^-,f^-)\;\text{ if }\;y\cdot\vec \nu<-1/2,\\
\xi(y)=G(\varphi^+,f^+)\;\text{ if }\; y\cdot\vec \nu> 1/2\;\text{
and }\;\xi\big(y+\vec \nu_j\big)=\xi(y)\;\;\forall j=2,\ldots,
N\Bigg\}\supset\mathcal{S}_0\big(\varphi^+,f^+,\varphi^-,f^-,I_{\vec
\nu}\big)\,.
\end{multline}
Then clearly by the density arguments we have
\begin{multline}\label{L2009hhffff12kkkhjhjghghgvgvggcjhggghtgjuthjhhhjjgjjhgjgnnmgmnghghgdgh}
\inf\Bigg\{\int_{I_{\vec \nu}}
\Big|\xi(y)-G\big(\zeta(y),\sigma_{f,x}(y)\big)\Big|^2\,dy:
\;\xi\in\mathcal{S}_0(\varphi^+,f^+,\varphi^-,f^-,I_{\vec\nu})\Bigg\}
\\=\inf\Bigg\{\int_{I_{\vec \nu}}
\Big|\xi(y)-G\big(\zeta(y),\sigma_{f,x}(y)\big)\Big|^2\,dy:
\;\xi\in\mathcal{T}(\varphi^+,f^+,\varphi^-,f^-,I_{\vec\nu})\Bigg\}\,.
\end{multline}
On the other hand, clearly there exists a minimizer to the r.h.s. of
\er{L2009hhffff12kkkhjhjghghgvgvggcjhggghtgjuthjhhhjjgjjhgjgnnmgmnghghgdgh},
i.e.
$\exists\xi_0\in\mathcal{T}(\varphi^+,f^+,\varphi^-,f^-,I_{\vec\nu})$,
such that
\begin{multline}\label{L2009hhffff12kkkhjhjghghgvgvggcjhggghtgjuthjhhhjjgjjhgjgnnmgmnghghgdghbbbvbvmvv}
\int_{I_{\vec \nu}}
\Big|\xi_0(y)-G\big(\zeta(y),\sigma_{f,x}(y)\big)\Big|^2dy
=\inf\Bigg\{\int_{I_{\vec \nu}}
\Big|\xi(y)-G\big(\zeta(y),\sigma_{f,x}(y)\big)\Big|^2dy:
\;\xi\in\mathcal{T}(\varphi^+,f^+,\varphi^-,f^-,I_{\vec\nu})\Bigg\}.
\end{multline}
Moreover, $\xi_0$ clearly satisfies
\begin{multline}\label{L2009hhffff12kkkhjhjghghgvgvggcjhggghtgjuthjhhhjjgjjhgjgnnmgmnghghgdghbbbvbvmvvbbnbnbnbkh}
\int_{I_{\vec \nu}}
\Big(\xi_0(y)-G\big(\zeta(y),\sigma_{f,x}(y)\big)\Big):\theta(y)\,dy\quad\text{for
every}\;\theta\in L^2_{loc}(\R^N,\R^{k\times N})\\ \text{such
that}\;\Div_y\theta(y)=0,\;\theta(y)=0\;\text{ if }\;|y\cdot\vec
\nu|>1/2,\;\text{ and }\;\theta\big(y+\vec
\nu_j\big)=\theta(y)\;\,\forall j=2,\ldots, N\,.
\end{multline}
In particular there exists $H\in W^{1,2}_{loc}(\R^N,\R^{k})$, such
that $H\big(y+\vec \nu_j\big)=H(y)\;\,\forall y\in\R^N,\;\forall
j=2,\ldots, N$ and $\nabla_y
H(y)=G\big(\zeta(y),\sigma_{f,x}(y)\big)-\xi_0(y)$ on $N_{\vec \nu}$
where $N_{\vec \nu}:=\big\{y\in\R^N:\;|y\cdot\vec \nu|<1/2\big\}$.
On the other hand since
$\xi_0\in\mathcal{T}(\varphi^+,f^+,\varphi^-,f^-,I_{\vec\nu})$ we
clearly have $\nabla_y H(y)\chi_{N_{\vec
\nu}}(y)=G\big(\zeta(y),\sigma_{f,x}(y)\big)-\xi_0(y)$ for every
$y\in\R^N$. Thus since $\Div_y \xi_0(y)\equiv 0$ we obtain
$\Div_y\big\{\nabla_y H(y)\chi_{N_{\vec \nu}}(y)\big\}=\Div_y
G\big(\zeta(y),\sigma_{f,x}(y)\big)$ on $\R^N$. Therefore,
$H(y)\equiv H_{G,\zeta,x,\vec\nu}$ where $H_{G,\zeta,x,\vec\nu}$
satisfies \er{defH12slbgghfjfghjhjhjhjjkjkjhhj}. Plugging it into
\er{L2009hhffff12kkkhjhjghghgvgvggcjhggghtgjuthjhhhjjgjjhgjgnnmgmnghghgdghbbbvbvmvv}
we deduce
\begin{equation}\label{L2009hhffff12kkkhjhjghghgvgvggcjhggghtgjuthjhhhjjgjjhgjgnnmgmnghghgdghbbbvbvmvvjhghgjk}
\int_{I_{\vec \nu}}
\Big|\nabla H_{G,\zeta,x,\vec\nu}(y)\Big|^2\,dy
=\inf\Bigg\{\int_{I_{\vec \nu}}
\Big|\xi(y)-G\big(\zeta(y),\sigma_{f,x}(y)\big)\Big|^2\,dy:
\;\xi\in\mathcal{T}(\varphi^+,f^+,\varphi^-,f^-,I_{\vec\nu})\Bigg\}\,.
\end{equation}
Therefore, using
\er{L2009hhffff12kkkhjhjghghgvgvggcjhggghtgjuthjhhhjjgjjhgjgnnmgmnghghgdgh}
and
\er{L2009hhffff12kkkhjhjghghgvgvggcjhggghtgjuthjhhhjjgjjhgjgnnmgmnghghgdghbbbvbvmvvjhghgjk}
we infer
\er{L2009hhffff12kkkhjhjghghgvgvggcjhggghtgjuthjhhhjjgjjhgjgnnmgmn}.
This completes the proof.
\end{proof}
Similarly plugging Lemma \ref{gbhfhgdfdghfddhg} into Theorem 2.3 in
\cite{PII}
we deduce the
following abstract lower bound result for problem with a non-local
term.
\begin{theorem}\label{dehgfrygfrgygenbgggggggggggggkgkgthtjtfjkjkkhhkjfjgfgghnew}
Let $\O\subset\R^N$ be an open set with locally Lipschitz's
boundary. Furthermore, let $p\geq 1$ and $F\in C^0\big(\R^{m\times
N^n}\times\ldots\times\R^{m\times N}\times \R^m,\R\big)$ and $G\in
C^1(\R^m,\R^{k\times N})$, be such that $F\geq 0$ and there exists a
constant $C>0$ such that $\big|G(b)\big|\leq C\big(|b|^{p/2}+1\big)$
for every $b\in\R^m$ and $|a_n|^p/C\leq F(a_n,\ldots, a_2,a_1,b)\leq
C\big(\sum_{j=1}^n|a_j|^p+|b|^p+1\big)$ for every $a_j\in\R^{m\times
N^j}$ and $b\in\R^m$. Next let $\varphi\in BV(\O,\R^{m})\cap
L^p(\O,\R^m)$ be such that $F\big(0,\ldots,0, \varphi(x)\big)=0$ for
a.e. $x\in\O$, $\Div_x G\big(\varphi(x)\big)=0$ in $\O$ and
$G\big(\varphi(x)\big)\cdot\vec n(x)=0$ on $\partial\O$.
Then for every sequence $\{\varphi_\e\}_{\e>0}\subset
W^{n,p}_{loc}(\O,\R^m)$
such that $\varphi_\e\to \varphi$ in $L^p_{loc}(\O,\R^m)$ as $\e\to
0^+$, we have
\begin{multline}\label{a1a2a3a4a5a6a7s8hhjhjjhjjjjjjkkkkgenhjhhhhjtjurtgfhfhfjfjfjhjjjnjhjjjhdgfgffgnew}
\liminf_{\e\to 0^+}\Bigg\{\frac{1}{\e}\int_\O F\Big(\e^n\nabla^n
\varphi_\e(x),\ldots,\e\nabla
\varphi_\e(x),\varphi_\e(x)\Big)dx+\frac{1}{\e}\int_{\R^N}|\nabla
V_{G,\varphi_\e}|^2dx\Bigg\}\geq\\ \int_{\O\cap J_\varphi}\hat
E_0\Big(\varphi^+(x),\varphi^-(x),\vec \nu(x)\Big)d \mathcal
H^{N-1}(x)\,,
\end{multline}
where $V_{G,\psi}:\R^N\to\R^k$ is defined by
\begin{equation}\label{defH12slbgghfjfghhghfnew}
\begin{cases}\Delta_x V_{G,\psi}(x)=\Div_x\Big\{\chi_\O G\big(\psi(x)\big)\Big\}\quad\text{in }\R^N\,,\\
\nabla V_{G,\psi}\in L^2(\R^N,\R^{k\times N})\,,
\end{cases}
\end{equation}
\begin{multline}\label{L2009hhffff12kkkhjhjghghgvgvggcjhggghtgjuthjhhhjjkgkjggjknew}
\hat E_{0}\Big(\varphi^+,\varphi^-,\vec \nu,x\Big)\;:=\; \inf\Bigg\{
\liminf_{\e\to 0^+}\int_{I_{\vec \nu}}\frac{1}{\e}
\Big(F\big(\e^n\nabla^n \zeta_\e(y),\ldots,\e\nabla
\zeta_\e(y),\zeta_\e(y)\big)+ \big|\nabla
H^0_{G,\zeta_\e,\vec\nu}(y)\big|^2\Big)\,dy:
\\
\zeta_\e\in\mathcal{S}^{(n)}_2\big(\varphi^+,\varphi^-,I_{\vec\nu}\big)
\;\;\text{such
that}\;\;\zeta_\e(y)\to\chi(y,\varphi^+,\varphi^-,\vec\nu)\;\;\text{in}\;\;L^p(I_{\vec
\nu},\R^m)\Bigg\}\,,
\end{multline}
where $H^0_{G,\zeta,\vec\nu}\in W^{1,2}_{0}(I_{\vec \nu},\R^k)$
satisfies
\begin{equation}\label{defH12slbgghfjfghjhjhjhjjkjkjhhjjgjggnew}
\Delta_y H^0_{G,\zeta,\vec\nu}(y)=\Div_y
G\big(\zeta(y)\big)\quad\text{in }I_{\vec\nu}\,,
\end{equation}
\begin{multline}\label{L2009Ddef2hhhjjjj77788hhhkkkkllkjjjjkkkhhhhffggdddkkkgjhikhhhjjddddhdkgkkknewjhgjgj}
\mathcal{S}^{(n)}_2\big(\varphi^+,\varphi^-,I_{\vec\nu}\big):=
\bigg\{\zeta\in C^n(\R^N,\R^m):\;\;\zeta(y)=\varphi^-\;\text{ if }\;y\cdot\vec\nu\leq-1/2,\\
\zeta(y)=\varphi^+\;\text{ if }\; y\cdot\vec\nu\geq 1/2\;\text{ and
}\;\zeta\big(y+\vec \nu_j\big)=\zeta(y)\;\;\forall j=2,\ldots,
N\bigg\}\,,
\end{multline}
and \begin{equation*}
\chi(y,\varphi^+,\varphi^-,\vec\nu):=\begin{cases}
\varphi^+\quad\text{if}\;\;y\cdot\vec\nu>0\,,\\
\varphi^-\quad\text{if}\;\;y\cdot\vec\nu<0\,.
\end{cases}
\end{equation*}
Here $I_{\vec \nu}:=\big\{y\in\R^N:\;|y\cdot
\vec\nu_j|<1/2\;\;\;\forall j=1,\ldots, N\big\}$ where
$\{\vec\nu_1,
\ldots,\vec\nu_N\}\subset\R^N$ is an
orthonormal base in $\R^N$ such that $\vec\nu_1:=\vec \nu$.
\end{theorem}
\begin{proof}
Let $\{\varphi_\e\}_{\e>0}\subset W^{n,p}_{loc}(\O,\R^m)$
be such that $\varphi_\e\to \varphi$ in $L^p_{loc}(\O,\R^m)$ as
$\e\to 0^+$. Without loss of generality we may assume
\begin{equation}\label{a1a2a3a4a5a6a7s8hhjhjjhjjjjjjkkkkgenhjhhhhjtjurtgfhfhfjfjfjhjjjnjhjjjhdgfgffgjhfjhjffjfjfjfffnew}
D:=\liminf_{\e\to 0^+}\Bigg\{\frac{1}{\e}\int_\O F\Big(\e^n\nabla^n
\varphi_\e(x),\ldots,\e\nabla
\varphi_\e(x),\varphi_\e(x)\Big)dx+\frac{1}{\e}\int_{\R^N}|\nabla
V_{G,\varphi_\e}|^2dx\Bigg\}<+\infty\,.
\end{equation}
Next set $L_\e:=\chi_\O G\big(\varphi_\e(x)\big)-\nabla
V_{G,\varphi_\e}$. Then we have $\Div L_\e\equiv 0$ on $\O$ and we
have
\begin{multline}\label{a1a2a3a4a5a6a7s8hhjhjjhjjjjjjkkkkgenhjhhhhjtjurtgfhfhfjfjfjhjjjnjhjjjhdgfgffgjhfjhjffjfjfjfffjkjkjnew}
D=\liminf_{\e\to 0^+}\Bigg\{\frac{1}{\e}\int_\O F\Big(\e^n\nabla^n
\varphi_\e(x),\ldots,\e\nabla
\varphi_\e(x),\varphi_\e(x)\Big)dx+\frac{1}{\e}\int_{\R^N}\Big|L_\e-
\chi_\O(x)G\big(\varphi_\e(x)\big)\Big|^2dx\Bigg\}\\ \geq
\liminf_{\e\to 0^+}\Bigg\{\frac{1}{\e}\int_\O F\Big(\e^n\nabla^n
\varphi_\e(x),\ldots,\e\nabla
\varphi_\e(x),\varphi_\e(x)\Big)dx+\frac{1}{\e}\int_{\O}\Big|L_\e-
G\big(\varphi_\e(x)\big)\Big|^2dx\Bigg\}\,.
\end{multline}
Thus applying Theorem 2.3
in \cite{PII} we
deduce
\begin{multline}\label{a1a2a3a4a5a6a7s8hhjhjjhjjjjjjkkkkgenhjhhhhjtjurtgfhfhfjfjfjhjjjnjhjjjhdgfgffgvvbjvghcvghnew}
D\geq\liminf_{\e\to 0^+}\Bigg\{\frac{1}{\e}\int_\O
F\Big(\e^n\nabla^n \varphi_\e(x),\ldots,\e\nabla
\varphi_\e(x),\varphi_\e(x)\Big)dx+\frac{1}{\e}\int_{\O}\Big|L_\e-
G\big(\varphi_\e(x)\big)\Big|^2dx\Bigg\}\geq\\ \int_{\O\cap
J_\varphi}\tilde E_0\Big(\varphi^+(x),\varphi^-(x),\vec \nu(x)\Big)d
\mathcal H^{N-1}(x)\,,
\end{multline}
where
\begin{multline}\label{L2009hhffff12kkkhjhjghghgvgvggcjhggghtgjuthjhhhjjkgkjggjkhgcfghvvgghkghghjknew}
\tilde E_{0}\Big(\varphi^+,\varphi^-,\vec \nu,x\Big)\;:=\;
\inf\Bigg\{ \liminf_{\e\to 0^+}\int_{I_{\vec \nu}}\frac{1}{\e}
\Big(F\big(\e^n\nabla^n \zeta_\e(y),\ldots,\e\nabla
\zeta_\e(y),\zeta_\e(y)\big)+
\Big|\xi_\e(y)-G\big(\zeta_\e(x)\big)\Big|^2\Big)\,dy:
\\
\zeta_\e\in W^{n,p}(I_{\vec \nu},\R^m)
,\,\xi_\e(y)\in W^{n,2}(I_{\vec \nu},\R^{k\times N})
\;\;\text{such that}\;
\Div_y\xi_\e(y)\equiv 0,\\
\zeta_\e(y)\to\chi(y,\varphi^+,\varphi^-,\vec\nu)\;\;\text{in}\;\;L^p(I_{\vec
\nu},\R^m)\;
\text{and}\;\xi_\e(y)\to\chi\big(y,G(\varphi^+),G(\varphi^-),\vec\nu\big)\;\;\text{in}\;\;L^2(I_{\vec
\nu},\R^{k\times N})\Bigg\}\,.
\end{multline}
On the other hand, by
Theorem 3.1 in \cite{PII}
we obtain
\begin{multline}\label{L2009hhffff12kkkhjhjghghgvgvggcjhggghtgjuthjhhhjjkgkjggjkhgcfghvvgghkghghjknewfhjhffjffj}
\tilde E_{0}\Big(\varphi^+,\varphi^-,\vec \nu,x\Big)\;\geq\;
\inf\Bigg\{ \liminf_{\e\to 0^+}\int_{I_{\vec \nu}}\frac{1}{\e}
\Big(F\big(\e^n\nabla^n \zeta_\e(y),\ldots,\e\nabla
\zeta_\e(y),\zeta_\e(y)\big)+
\Big|\xi_\e(y)-G\big(\zeta_\e(x)\big)\Big|^2\Big)\,dy:
\\
\zeta_\e\in
\mathcal{S}^{(n)}_2\big(\varphi^+,\varphi^-,I_{\vec\nu}\big)
,\,\xi_\e(y)\in W^{n,2}(I_{\vec \nu},\R^{k\times N})
\;\;\text{such that}\;
\Div_y\xi_\e(y)\equiv 0,\\
\zeta_\e(y)\to\chi(y,\varphi^+,\varphi^-,\vec\nu)\;\;\text{in}\;\;L^p(I_{\vec
\nu},\R^m)\;
\text{and}\;\xi_\e(y)\to\chi\big(y,G(\varphi^+),G(\varphi^-),\vec\nu\big)\;\;\text{in}\;\;L^2(I_{\vec
\nu},\R^{k\times N})\Bigg\}\,.
\end{multline}
Therefore, using Lemma \ref{gbhfhgdfdghfddhg} we obtain
\begin{equation}\label{bxhvgjvdfjghfgfghghfnew}
\tilde E_{0}\big(\varphi^+,\varphi^-,\vec \nu,x\big)\geq\hat
E_{0}\big(\varphi^+,\varphi^-,\vec \nu,x\big)\,,
\end{equation}
where $\hat E_{0}\big(\varphi^+,\varphi^-,\vec \nu,x\big)$ is
defined by
\er{L2009hhffff12kkkhjhjghghgvgvggcjhggghtgjuthjhhhjjkgkjggjknew}.
Thus, plugging \er{bxhvgjvdfjghfgfghghfnew} into
\er{a1a2a3a4a5a6a7s8hhjhjjhjjjjjjkkkkgenhjhhhhjtjurtgfhfhfjfjfjhjjjnjhjjjhdgfgffgvvbjvghcvghnew}
we deduce
\begin{multline}\label{a1a2a3a4a5a6a7s8hhjhjjhjjjjjjkkkkgenhjhhhhjtjurtgfhfhfjfjfjhjjjnjhjjjhdgfgffgvvbjvghcvghghnew}
D=\liminf_{\e\to 0^+}\Bigg\{\frac{1}{\e}\int_\O F\Big(\e^n\nabla^n
\varphi_\e(x),\ldots,\e\nabla
\varphi_\e(x),\varphi_\e(x)\Big)dx+\frac{1}{\e}\int_{\R^N}|\nabla
V_{G,\varphi_\e}|^2dx\Bigg\}\\ \geq \int_{\O\cap J_\varphi}\hat
E_0\Big(\varphi^+(x),\varphi^-(x),\vec \nu(x)\Big)d \mathcal
H^{N-1}(x)\,,
\end{multline}
and the result follows.
\end{proof}

\section{The problem related to the theory of Conservation
Laws}\label{hgvfghfghoj}

\subsection{Some definitions and preliminaries}
\begin{definition}\label{3bdf}
For a given Banach space $X$ with the associated norm $\|\cdot\|_X$
and a real interval $(a,b)$ we denote by $L^q(a,b;X)$ the linear
space of (equivalence classes of) strongly measurable (i.e
equivalent to some strongly Borel mapping)
functions $f:(a,b)\to X$ such that the functional
\begin{equation*}
\|f\|_{L^q(a,b;X)}:=
\begin{cases}
\Big(\int_a^b\|f(t)\|^q_X dt\Big)^{1/q}\quad\text{if }\;1\leq
q<\infty\\
{\text{es$\,$sup}}_{t\in (a,b)}\|f(t)\|_X\quad\text{if }\; q=\infty
\end{cases}
\end{equation*}
is finite. It is known that this functional defines a norm with
respect to which $L^q(a,b;X)$ becomes a Banach space.
\end{definition}
\begin{definition}\label{dXY}
Let $\O\subset\R^N$ be an open set. We denote by $\tilde
H^1_0(\O,\R^k)$ the closure of $C_c^\infty(\O,\R^k)$ with respect to
the norm $|||\f|||:=\big(\int_{\O}|\nabla\f|^2dx\big)^{1/2}$ (this
space differ from $W^{1,2}_0(\O,\R^k)$ only in the case of unbounded
domain $\O$) and denote by $\tilde H^{-1}(\O,\R^k)$ the space dual
to $\tilde H^1_0(\O,\R^k)$.
\end{definition}
\begin{remark}\label{rem}
It is obvious that $u\in\mathcal{D}'(\O,\R^k)$ belongs to $\tilde
H^{-1}(\O,\R^k)$ if and only if there exists $w\in \tilde
H^1_0(\O,\R^k)$ such that
\begin{equation}\label{ghhhgghhg}
\int_{\O}\nabla w:\nabla \delta\, dx=
-<u,\delta>\quad\quad\forall\delta\in C^\infty_c(\O,\R^k)\,,
\end{equation}
Note that \er{ghhhgghhg} is equivalent to that $\Delta w=u$ as
distributions. Moreover,
$$|||w|||=\sup\limits_{\delta\in \tilde H^1_0(\O,\R^k),\;|||\delta|||\leq
1}<u,\delta>=|||u|||_{-1}\,.$$ Finally observe that $w$ is uniquely
defined by $u$.
\end{remark}
\begin{remark}\label{remkllk}
It is obvious that $u\in \mathcal{D}'\big(\O\times(0,T),\R^k\big)$
belongs to $L^2\big(0,T;\tilde H^{-1}(\O,\R^k)\big)$ if and only if
there exists $w\in L^2\big(0,T;\tilde H^1_0(\O,\R^k)\big)$ such that
\begin{equation}\label{ghhhgghhgjjh}
\int_0^T\int_{\O}\nabla_x w(x,t):\nabla_x \delta(x,t)\, dxdt=
-\big<u,\delta\big>\quad\quad\forall\delta\in
C_c^\infty\big(\O\times(0,T),\R^k\big)\,.
\end{equation}
Note that \er{ghhhgghhgjjh} is equivalent to that $\Delta_x w=u$ as
distributions. Moreover,
$$\big\|w\big\|_{L^2(0,T;\tilde H^1_0(\O,\R^k))}=\big\|u\big\|_{L^2(0,T;\tilde H^{-1}(\O,\R^k))}\,.$$
Finally observe that $w$ is uniquely defined by $u$.
\end{remark}
%
%
%

\begin{lemma}\label{vfgcfgjfvhgkjkjkkk}
Let $\O\subset\R^N$ be an open set (possibly unbounded) with locally
Lipschitz's boundary and let $T>0$. Furthermore, let
$G(\psi,x,t):\R^k\times\O\times(0,T)\to\R^{k\times N}$ be a
measurable function, continuous by the first argument $\psi$, which
satisfies
\begin{equation}\label{ghhghjghjhbfggffgkkk}
\big|G(\psi,x,t)\big|\leq
C_0|\psi|^{p/2}+h_0(x,t)\quad\forall\psi\in\R^k,\;x\in\O\,,
\end{equation}
for some constant $C_0>0$, $p\geq 1$ and $h_0\in
L^2\big(\O\times(0,T),\R\big)$. Furthermore, for every $\e>0$
consider the functional
$E_\e\big(\psi(x,t)\big):L^p\big(\O\times(0,T),\R^k\big)\to[0,+\infty)\cup\{+\infty\}$
which (possibly) can attain the infinite values. Next for every
$\e>0$ consider the functional
$P_\e\big(\psi(x,t)\big):L^p\big(\O\times(0,T),\R^k\big)\to[0,+\infty)\cup\{+\infty\}$,
defined by
\begin{equation}\label{fhfhjffhjhfkkk}
P_\e\big(\psi(x,t)\big):=\begin{cases}E_\e\big(\psi(x,t)\big)+\frac{1}{\delta_\e}\int_0^T\int_{\O}\Big|\nabla_x
V_{G,\psi}(x,t)\Big|^2 dxdt\quad\text{if}\quad
\partial_t\psi(x,t)\in L^2\big(0,T;\tilde
H^{-1}(\O,\R^k)\big)\,,\\+\infty\quad\text{otherwise}\,,
\end{cases}
\end{equation}
where $\delta_\e>0$ satisfies $\lim_{\e\to 0^+}\delta_\e=0$ and
given $\psi(x,t)\in L^p\big(\O\times(0,T),\R^k\big)$,
$V_{G,\psi}(x,t):\R^N\to\R^k$ is defined by
\begin{equation}\label{defH12slvvvbkkk}
\begin{cases}\Delta_x V_{G,\psi}(x,t)=\partial_t\psi(x,t)+\Div_x G\big(\psi(x,t),x,t\big)\quad\text{in }\O\times(0,T)\,,\\
V_{G,\psi}(x,t)\in L^2\big(0,T;H^1_0(\O,\R^k)\big)\,,
\end{cases}
\end{equation}
Furthermore, for every $\e>0$ consider the functional
$Q_\e\big(\psi(x,t),L(x,t)\big):L^p\big(\O\times(0,T),\R^k\big)\times
L^2\big(\O\times(0,T),\R^{k\times N}\big)
\to[0,+\infty)\cup\{+\infty\}$ defined by
\begin{multline}\label{fhfhjffhjhfgyukuoijkkjkkk}
Q_\e\big(\psi(x,t),L(x,t)\big):=\\=
\begin{cases}E_\e\big(\psi(x,t)\big)+\frac{1}{\delta_\e}\int_0^T\int_{\O}\big|L(x,t)-
G\big(\psi(x,t),x,t\big)\big|^2 dxdt\quad\quad\text{if}\quad
\partial_t\psi(x,t)+\Div_x
L(x,t)=0\,,\\+\infty\quad\quad\quad\text{otherwise}\,.
\end{cases}
\end{multline}
Finally for every $\e>0$ consider the functional
$R_\e(u):\mathcal{D}'\big(\O\times(0,T),\R^{k\times N}\big)
\to[0,+\infty)\cup\{+\infty\}$ defined by
\begin{multline}\label{fhfhjffhjhfgyukuoijkkjkkkjhh}
R_\e\big(u\big):=
\begin{cases}E_\e\big(\Div_x u(x,t)\big)+\frac{1}{\delta_\e}\int_0^T\int_{\O}\big|\partial_t
u(x,t)+G\big(\Div_x u(x,t),x,t\big)\big|^2 dxdt\\
\quad\quad\quad\quad\quad\quad\quad\quad\quad
\text{if}\;\;\partial_t u\in L^2(\O\times(0,T),\R^{k\times
N})\;\;\text{and}\;\;\Div_x u\in
L^p(\O\times(0,T),\R^{k})\,,\\+\infty\quad\quad\quad\quad\quad\quad\quad\quad\quad\quad\quad\quad\text{otherwise}\,.
\end{cases}
\end{multline}
Next for every $\varphi(x,t)\in L^p\big(\O\times(0,T),\R^k\big)$,
such that $\partial_t \varphi(x,t)+\Div_x
G\big(\varphi(x,t),x,t\big)=0$ on $\O$, set
\begin{equation}\label{hfghjfhfhfhkkk}
\begin{split}
\underline P(\varphi):=\inf\Bigg\{\liminf_{\e\to
0^+}P_\e\big(\psi_\e(x,t)\big):\;\;\psi_\e(x,t)\to\varphi(x,t)\;\;\text{in}\;\;L^p\big(\O\times(0,T),\R^k\big)\Bigg\}\,,\\
\ov P(\varphi):=\inf\Bigg\{\limsup_{\e\to
0^+}P_\e\big(\psi_\e(x,t)\big):\;\;\psi_\e(x,t)\to\varphi(x,t)\;\;\text{in}\;\;L^p\big(\O\times(0,T),\R^k\big)\Bigg\}\,,
\end{split}
\end{equation}
\begin{equation}\label{hfghjfhfhfhj;k;jkhkhhkkkkbhhhhjhjh}
\begin{split}
\underline Q(\varphi):=\inf\Bigg\{\liminf_{\e\to
0^+}Q_\e\big(\psi_\e(x,t),L_\e(x,t)\big):\;\;\psi_\e(x,t)\to\varphi(x,t)\;\;\text{in}\;L^p\big(\O\times(0,T),\R^k\big)\\
\text{and}\;\;L_\e(x,t)\to
G\big(\varphi(x,t),x,t\big)\;\;\text{in}\;\;L^2\big(\O\times(0,T),\R^{k\times
N}\big)\Bigg\}\,,\\ \ov Q(\varphi):=\inf\Bigg\{\limsup_{\e\to
0^+}Q_\e\big(\psi_\e(x,t),L_\e(x,t)\big):\;\;\psi_\e(x,t)\to\varphi(x,t)\;\;\text{in}\;L^p\big(\O\times(0,T),\R^k\big)\\
\text{and}\;\;L_\e(x,t)\to
G\big(\varphi(x,t),x,t\big)\;\;\text{in}\;\;L^2\big(\O\times(0,T),\R^{k\times
N}\big)\Bigg\}\,,
\end{split}
\end{equation}
and
\begin{equation}\label{hfghjfhfhfhj;k;jkhkhhkkkk}
\begin{split}
\underline R(\varphi):=\inf\Bigg\{\liminf_{\e\to
0^+}R_\e\big(u_\e(x,t)\big):\;\;\Div_x u_\e(x)\to\varphi(x)\;\;\text{in}\;\;L^p\big(\O\times(0,T),\R^k\big)\\
\text{and}\;\;\partial_t u_\e(x)\to\Big(-
G\big(\varphi(x,t),x,t\big)\Big)\;\;\text{in}\;\;L^2\big(\O\times(0,T),\R^{k\times
N}\big)\Bigg\}\,,\\ \ov R(\varphi):=\inf\Bigg\{\limsup_{\e\to
0^+}R_\e\big(u_\e(x,t)\big):\;\;\Div_x u_\e(x)\to\varphi(x)\;\;\text{in}\;\;L^p\big(\O\times(0,T),\R^k\big)\\
\text{and}\;\;\partial_t u_\e(x)\to\Big(-
G\big(\varphi(x,t),x,t\big)\Big)\;\;\text{in}\;\;L^2\big(\O\times(0,T),\R^{k\times
N}\big)\Bigg\}\,.
\end{split}
\end{equation}
Then we have the following equalities
\begin{equation}\label{vhfhvffgjfghkjjjjjkjkhkkkkjhggjgj}
\underline P(\varphi)=\underline Q(\varphi)=\underline
R(\varphi)\quad\text{and}\quad \ov P(\varphi)=\ov Q(\varphi)=\ov
R(\varphi)\,.
\end{equation}
\end{lemma}
\begin{proof}
Fix some $\psi(x,t)\in L^p\big(\O\times(0,T),\R^k\big)$, such that
$\partial_t\psi(x,t)\in L^2\big(0,T;\tilde H^{-1}(\O,\R^k)\big)$.
Then by \er{ghhghjghjhbfggffgkkk}, we clearly have
$G\big(\psi(x,t),x,t\big)\in L^2\big(\O\times(0,T),\R^{k\times
N}\big)$ and then by Remark \ref{remkllk} we have $\Div_x
G\big(\psi(x,t),x,t\big)\in L^2\big(0,T;\tilde
H^{-1}(\O,\R^k)\big)$. Thus
by Remark \ref{remkllk} there exist uniquely defined
$H_{\psi}(x,t)\in L^2\big(0,T;\tilde H^1_0(\O,\R^k)\big)$,
$D_{G,\psi}(x,t)\in L^2\big(0,T;\tilde H^1_0(\O,\R^k)\big)$ and
$V_{G,\psi}(x,t)\in L^2\big(0,T;\tilde H^1_0(\O,\R^k)\big)$ such
that
\begin{equation}\label{defH12slvvvbkkkvhgjhvg}
\begin{cases}
\Delta_x H_{\psi}(x,t)=\partial_t\psi(x,t)\quad\text{in
}\O\times(0,T)\,,\\
\Delta_x D_{G,\psi}(x,t)=\Div_x
G\big(\psi(x,t),x,t\big)\quad\text{in }\O\times(0,T)\,,\\ \Delta_x
V_{G,\psi}(x,t)=\partial_t\psi(x,t)+\Div_x
G\big(\psi(x,t),x,t\big)\quad\text{in }\O\times(0,T)\,,
\end{cases}
\end{equation}
and clearly
\begin{equation}\label{defH12slvvvbkkkvhgjhvgnbjjvvvhcghcghgh}
V_{G,\psi}(x,t)\equiv H_{\psi}(x,t)+ D_{G,\psi}(x,t)\quad\text{in
}\O\times(0,T)\,.
\end{equation}
Moreover, by Lemma \ref{gbhfhgdfdghfddhg} for every $U(x,t)\in
L^2\big(\O\times(0,T),\R^{k\times N}\big)$, such that $\Div_x
U(x,t)\equiv 0$ we have
\begin{equation}\label{fhfhjffhjhfgyukuoijkkjkkkhhhjjjhhgjhhyjtjkhjkjkjk}
\int_0^T\int_{\O}\big|U(x,t)-\nabla_x H_{\psi}(x,t)-
G\big(\psi(x,t),x,t\big)\big|^2 dxdt\geq
\int_0^T\int_{\O}\big|\nabla_x V_{\psi}(x,t)\big|^2 dxdt\,,
\end{equation}
and if we denote
\begin{equation}\label{xfgvfhgjgjykuikuk}
U_{G,\psi}(x,t):=G\big(\psi(x,t),x,t\big)-\nabla_x
D_{G,\psi}(x,t)\quad\forall(x,t)\in\O\times (0,T)\,,
\end{equation}
then $U_{G,\psi}(x,t)\in L^2\big(\O\times(0,T),\R^{k\times N}\big)$,
$\Div_x U_{G,\psi}(x,t)\equiv 0$ and
\begin{equation}\label{fhfhjffhjhfgyukuoijkkjkkkhhhjjjhhgjhhyjtjkhjkjkjkhjhjh}
\int_0^T\int_{\O}\big|U_{G,\psi}(x,t)-\nabla_x H_{\psi}(x,t)-
G\big(\psi(x,t),x,t\big)\big|^2 dxdt=\int_0^T\int_{\O}\big|\nabla_x
V_{\psi}(x,t)\big|^2 dxdt\,.
\end{equation}
In particular, by
\er{fhfhjffhjhfgyukuoijkkjkkkhhhjjjhhgjhhyjtjkhjkjkjk} for every
$L(x,t)\in L^2\big(\O\times(0,T),\R^{k\times N}\big)$, such that
$\Div_x L(x,t)+\partial_t\psi(x,t)\equiv 0$ we have
\begin{equation}\label{fhfhjffhjhfgyukuoijkkjkkkhhhjjjhhgjhhyjtjkhjkjkjkvncccxxcxcdf}
\int_0^T\int_{\O}\big|L(x,t)- G\big(\psi(x,t),x,t\big)\big|^2
dxdt\geq \int_0^T\int_{\O}\big|\nabla_x V_{\psi}(x,t)\big|^2 dxdt\,,
\end{equation}
and if we denote
\begin{equation}\label{xfgvfhgjgjykuikukyuyuuyuyiigf}
L_{G,\psi}(x,t):=G\big(\psi(x,t),x,t\big)-\nabla_x
V_{G,\psi}(x,t)\quad\forall(x,t)\in\O\times (0,T)\,,
\end{equation}
then $L_{G,\psi}(x,t)\in L^2\big(\O\times(0,T),\R^{k\times N}\big)$,
$\Div_x L_{G,\psi}(x,t)+\partial_t\psi(x,t)\equiv 0$ and
\begin{equation}\label{fhfhjffhjhfgyukuoijkkjkkkhhhjjjhhgjhhyjtjkhjkjkjkhjhjhjhjhghgg}
\int_0^T\int_{\O}\big|L_{G,\psi}(x,t)-
G\big(\psi(x,t),x,t\big)\big|^2 dxdt=\int_0^T\int_{\O}\big|\nabla_x
V_{\psi}(x,t)\big|^2 dxdt\,.
\end{equation}
Finally define
\begin{equation}\label{xfgvfhgjgjykuikukyuyuuyuyiichgvdfhjh}
u_\psi(x,t):=K_\psi(x)-\int_0^t
L_{G,\psi}(x,s)\,ds\quad\quad\quad\quad\forall(x,t)\in\O\times
(0,T)\,,
\end{equation}
where $K_\psi(x):\O\to\R^{k\times N}$ satisfies $\Div_x
K_\psi(x)\equiv\psi(x,0)$. Then since $\Div_x
L_{G,\psi}(x,t)+\partial_t\psi(x,t)\equiv 0$ we deduce that
\begin{equation}\label{xfgvfhgjgjykuikukyuyuuyuyiichgvdfhjhhtjtyvvjhjjhj}
\partial_t u_\psi(x,t):=-
L_{G,\psi}(x,t)\quad\text{and}\quad \Div_x
u_\psi(x,t)=\psi(x,t)\quad\quad\forall(x,t)\in\O\times (0,T)\,.
\end{equation}
Therefore, by \er{fhfhjffhjhfkkk}, \er{fhfhjffhjhfgyukuoijkkjkkk}
and \er{fhfhjffhjhfgyukuoijkkjkkkjhh}, using
\er{fhfhjffhjhfgyukuoijkkjkkkhhhjjjhhgjhhyjtjkhjkjkjkhjhjhjhjhghgg}
and \er{xfgvfhgjgjykuikukyuyuuyuyiichgvdfhjhhtjtyvvjhjjhj} we deduce
\begin{equation}\label{xfgvfhgjgjykuikukyuyuuyuyiichgvdfhjhhtjtyvvjhjjhjhnbnunjkjklk}
P_\e\big(\psi(x,t)\big)=Q_\e\big(\psi(x,t),L_{G,\psi}(x,t)\big)=R_\e\big(u_\psi\big)\,.
\end{equation}
Moreover, by \er{fhfhjffhjhfgyukuoijkkjkkk} and
\er{fhfhjffhjhfgyukuoijkkjkkkhhhjjjhhgjhhyjtjkhjkjkjkvncccxxcxcdf}
we have
\begin{equation}\label{xfgvfhgjgjykuikukyuyuuyuyiichgvdfhjhhtjtyvvjhjjhjhnbnunjkjklkjuuhyuuyhjgj}
P_\e\big(\psi(x,t)\big)\leq
Q_\e\big(\psi(x,t),L(x,t)\big)\quad\quad\forall L(x,t)\in
L^2\big(\O\times(0,T),\R^{k\times N}\big)\,.
\end{equation}
Next fix some $\varphi(x,t)\in L^p\big(\O\times(0,T),\R^k\big)$,
such that
\begin{equation}\label{fgfjhgjjghgvjhh}
\partial_t \varphi(x,t)+\Div_x
G\big(\varphi(x,t),x,t\big)=0\quad\text{in}\;\;\O\times(0,T)\,.
\end{equation}
Then, using the fact that given $\psi(x,t)\in
L^p\big(\O\times(0,T),\R^k\big)$ and $L(x,t)\in
L^2\big(\O\times(0,T),\R^{k\times N}\big)$, such that $\Div_x
L(x,t)+\partial_t\psi(x,t)\equiv 0$ we always have
$\partial_t\psi(x,t)\in L^2\big(0,T;\tilde H^{-1}(\O,\R^k)\big)$, by
\er{xfgvfhgjgjykuikukyuyuuyuyiichgvdfhjhhtjtyvvjhjjhjhnbnunjkjklkjuuhyuuyhjgj}
we obtain
\begin{equation}\label{vhfhvffgjfghkjjjjjkjkhkkkk}
\underline P(\varphi)\leq\underline Q(\varphi)\quad\text{and}\quad
\ov P(\varphi)\leq\ov Q(\varphi)\,.
\end{equation}
Furthermore, fix some sequence $\{\psi_\e(x,t)\}\subset
L^p\big(\O\times(0,T),\R^k\big)$ such that
$\psi_\e(x,t)\to\varphi(x,t)$ in $L^p\big(\O\times(0,T),\R^k\big)$
as $\e\to 0^+$ and for a subsequence $\e_n\downarrow 0$ we have
$\lim_{n\to+\infty}P_{\e_n}\big(\psi_{\e_n}(x,t)\big)<+\infty$. Then
since $\delta_{\e_n}\to 0^+$, by \er{fhfhjffhjhfkkk} we deduce
\begin{equation}\label{jdfhhghfhgfhkgfhkdhbnbjlbhbjjnmkkk}
\lim_{n\to +\infty}V_{G,\psi_{\e_n}}(x,t)=0\quad\text{in}\quad
L^2\big(\O\times(0,T),\R^{k\times N}\big)\,.
\end{equation}
On the other hand, since $\psi_\e(x,t)\to\varphi(x,t)$ in
$L^p\big(\O\times(0,T),\R^k\big)$, by \er{ghhghjghjhbfggffgkkk} we
have
\begin{equation}\label{jdfhhghfhgfhkgfhkdhbnbjlbhbjjnmkkkhhh}
\lim_{\e\to 0^+}G\big(\psi_\e(x,t),x,t\big)=
G\big(\varphi(x,t),x,t\big)\quad\text{in}\;\;
L^2\big(\O\times(0,T),\R^{k\times N}\big)\,.
\end{equation}
Thus by \er{xfgvfhgjgjykuikukyuyuuyuyiigf},
\er{jdfhhghfhgfhkgfhkdhbnbjlbhbjjnmkkk} and
\er{jdfhhghfhgfhkgfhkdhbnbjlbhbjjnmkkkhhh} we deduce that
\begin{equation}\label{jdfhhghfhgfhkgfhkdhbnbjlbhbjjnmkkkhhhjhhigyygygh}
\lim_{n\to+\infty}L_{G,\psi_{\e_n}}(x,t)=
G\big(\varphi(x,t),x,t\big)\quad\text{in}\;\;
L^2\big(\O\times(0,T),\R^{k\times N}\big)\,,
\end{equation}
and therefore, by
\er{xfgvfhgjgjykuikukyuyuuyuyiichgvdfhjhhtjtyvvjhjjhj} we have
\begin{multline}\label{xfgvfhgjgjykuikukyuyuuyuyiichgvdfhjhhtjtyvvjhjjhjfggfgfjhhk}
\lim_{n\to+\infty}\partial_t u_{\psi_{\e_n}}(x,t)=-
G\big(\varphi(x,t),x,t\big)\;\;\text{in}\;\;L^2\big(\O\times(0,T),\R^{k\times
N}\big)\quad\text{and}\\ \lim_{n\to+\infty}\Div_x
u_{\psi_{\e_n}}(x,t)=\varphi(x,t)\;\;\text{in}\;\;L^p\big(\O\times(0,T),\R^{k}\big)\,.
\end{multline}
Moreover, by
\er{xfgvfhgjgjykuikukyuyuuyuyiichgvdfhjhhtjtyvvjhjjhjhnbnunjkjklk}
we have
\begin{equation}\label{xfgvfhgjgjykuikukyuyuuyuyiichgvdfhjhhtjtyvvjhjjhjhnbnunjkjklkjfgj}
\lim_{n\to+\infty}P_{\e_n}\big(\psi_{\e_n}(x,t)\big)=
\lim_{n\to+\infty}Q_{\e_n}\big(\psi_{\e_n}(x,t),L_{G,\psi_{\e_n}}(x,t)\big)
=\lim_{n\to+\infty}R_{\e_n}\big(u_{\psi_{\e_n}}\big)\,.
\end{equation}
Thus since the sequence $\{\psi_\e\}$ was arbitrary, we get
\begin{equation}\label{vhfhvffgjfghkjjjjjkjkhkkkkfffj}
\underline P(\varphi)\geq\underline Q(\varphi)\quad\text{and}\quad
\ov P(\varphi)\geq\ov Q(\varphi)\,.
\end{equation}
and
\begin{equation}\label{vhfhvffgjfghkjjjjjkjkhkkkkfffjgjggggg}
\underline P(\varphi)\geq\underline R(\varphi)\quad\text{and}\quad
\ov P(\varphi)\geq\ov R(\varphi)\,.
\end{equation}
Thus plugging \er{vhfhvffgjfghkjjjjjkjkhkkkkfffj} into
\er{vhfhvffgjfghkjjjjjkjkhkkkk} we obtain
\begin{equation}\label{vhfhvffgjfghkjjjjjkjkhkkkkjhkjgbjk}
\underline P(\varphi)=\underline Q(\varphi)\quad\text{and}\quad \ov
P(\varphi)=\ov Q(\varphi)\,.
\end{equation}
Finally fix arbitrary sequence
$\{u_\e\}_{\e>0}\subset\mathcal{D}'\big(\O\times(0,T),\R^{k\times
N}\big)$ such that
$$\Div_x
u_\e(x)\to\varphi(x)\;\text{in}\;L^p\big(\O\times(0,T),\R^k\big)\;\;
\text{and}\;\;\partial_t u_\e(x)\to\Big(-
G\big(\varphi(x,t),x,t\big)\Big)\;\text{in}\;L^2\big(\O\times(0,T),\R^{k\times
N}\big)\,.$$ Then if we set
\begin{equation}\label{xfgvfhgjgjykuikukyuyuuyuyiichgvdfhjhhtjtyvvjhjjhjgjkgjkgjk}
L_\e(x,t):=-\partial_t u_\e(x,t) \quad\text{and}\quad
\psi_\e(x,t):=\Div_x u_\e(x,t)\quad\quad\forall(x,t)\in\O\times
(0,T)\,,
\end{equation}
we obtain $\Div L_\e+\partial_t\psi_\e\equiv 0$ and
$$\psi_\e(x,t)\to\varphi(x)\;\text{in}\;L^p\big(\O\times(0,T),\R^k\big)\;\;
\text{and}\;\;L_\e(x,t)\to
G\big(\varphi(x,t),x,t\big)\;\text{in}\;L^2\big(\O\times(0,T),\R^{k\times
N}\big)\,.$$ Moreover, by
\er{xfgvfhgjgjykuikukyuyuuyuyiichgvdfhjhhtjtyvvjhjjhjgjkgjkgjk} we
have
$$Q_\e\big(\psi_\e(x,t),L_\e(x,t)\big)=R_\e\big(u_\e\big)\,.$$
Therefore, since the sequence $\{u_\e\}_{\e>0}$ was arbitrary, we
deduce
\begin{equation}\label{vhfhvffgjfghkjjjjjkjkhkkkkfffjgjgggggjkjgkjkklhj}
\underline Q(\varphi)\leq\underline R(\varphi)\quad\text{and}\quad
\ov Q(\varphi)\leq\ov R(\varphi)\,.
\end{equation}
Thus by plugging \er{vhfhvffgjfghkjjjjjkjkhkkkkjhkjgbjk},
\er{vhfhvffgjfghkjjjjjkjkhkkkkfffjgjggggg} and
\er{vhfhvffgjfghkjjjjjkjkhkkkkfffjgjgggggjkjgkjkklhj} we finally
deduce \er{vhfhvffgjfghkjjjjjkjkhkkkkjhggjgj}.
\end{proof}

\begin{definition}\label{bvclllffdftyyuuyuyy}
Let $F(u)=\{F_{ij}(u)\}\in C^1(\R^k,\R^{k\times N})$. Set $\bar
F_j(u):=\big(F_{1j}(u),\ldots,F_{kj}(u)\big):\R^k\to\R^k$ $\forall
j\in\{1,\ldots,N\}$. Consider the system of Conservation Laws
\begin{equation}\label{equ24}
\partial_t u+div_x
F(u)=0\quad\forall (x,t)\in\R^N\times (0,T)\,.
\end{equation}
We say that the function $\eta(u)\in C^1(\R^k,\R)$ is an entropy for
the system \er{equ24} and
$\Psi(u):=\big(\Psi_1(u),\ldots,\Psi_N(u)\big)\in C^1(\R^k,\R^N)$ is
an entropy flux associated with $\eta$ if we have
\begin{equation}\label{equ25}
\nabla_u\Psi_j(u)=\nabla_u\eta(u)\cdot\nabla_u \bar F_j(u)\quad\quad
\forall\, u\in\R^k,j\in\{1,\ldots,N\}\,.
\end{equation}
\end{definition}
Let $F(u)=\{F_{ij}(u)\}\in C^1(\R^k,\R^{k\times N})$ and $\eta(u)\in
C^2(\R^k,\R)$ be an entropy for the system \er{equ24}, which
satisfies $\eta(u)\geq 0$ and $\eta(0)=0$, and
$\Psi(u):=\big(\Psi_1(u),\ldots,\Psi_N(u)\big)\in C^1(\R^k,\R^N)$ be
the corresponding entropy flux associated with $\eta$. Considered
the following family of energy functionals $\{I_{\e,F}(u)\}$,
defined for $u(x,t):\R^N\times[0,T]\to\R^k$ by
\begin{equation}\label{cons}
I_{\e,F}(u):=\int_0^T\int_{\R^N}\bigg(\e\Big|\nabla_x
\big\{\nabla_u\eta\big(u(x,t)\big)\big\}\Big|^2+\frac{1}{\e}\Big|\nabla_x
H_{F,u}(x,t)\Big|^2\bigg)\,dxdt+\int_{\R^N}\eta\big(u(x,T)\big)\,dx\,,
\end{equation}
where $H_{F,u}(x,t)\in L^2\big(0,T;\tilde H^1_0(\R^N,\R^k)\big)$
satisfies
\begin{equation}\label{conssst}
\Delta_x H_{F,u}(x,t)=\partial_t u(x,t)+div_x F\big(u(x,t)\big)\,,
\end{equation}
and we assume that
\begin{equation}\label{jhniguiguyyjkj}
u(x,t)\in L^2\big(0,T;\tilde H^1_0(\R^N,\R^k)\big)\cap
C\big(0,T;L^2(\R^N,\R^k)\big)\cap L^\infty\;\;\text{and}\;\;
\partial_t u(x,t)\in L^2\big(0,T;\tilde H^{-1}(\R^N,\R^k)\big)\,,
\end{equation}
Since
\begin{multline*}
-\int_0^T\int_{\R^N}\nabla_x
\big\{\nabla_u\eta\big(u(x,t)\big)\big\}:\nabla_x
H_{F,u}(x,t)\,dxdt=\int_0^T\int_{\R^N}
\nabla_u\eta\big(u(x,t)\big)\cdot\Delta_x H_{F,u}(x,t)\,dxdt\\=
\int_0^T\int_{\R^N} \nabla_u\eta\big(u(x,t)\big)\cdot\Big(\partial_t
u(x,t)+div_x
F\big(u(x,t)\big)\Big)\,dxdt\\=\int_{\R^N}\bigg(\int_0^T
\partial_t\big\{\eta\big(u(x,t)\big)\big\}\,dt\bigg)dx+ \int_0^T\int_{\R^N}\sum_{j=1}^{N}\nabla_u\eta\big(u(x,t)\big)\cdot
\nabla_u\bar F_j\big(u(x,t)\big)\cdot\frac{\partial u}{\partial
x_j}(x,t)\,dxdt\\=\int_{\R^N}\Big(\eta\big(u(x,T)\big)-\eta\big(u(x,0)\big)\Big)\,dx+
\int_0^T\int_{\R^N}\Div_x\Psi\big(u(x,t)\big)\,dxdt=\int_{\R^N}\Big(\eta\big(u(x,T)\big)-\eta\big(u(x,0)\big)\Big)\,dx\,,
\end{multline*}
we can rewrite the expression of $I_{\e,F}(u)$ as
\begin{equation}\label{conslll}
I_{\e,F}(u)=\int_0^T\int_{\R^N}\frac{1}{\e}\bigg|\e\nabla_x
\big\{\nabla_u\eta\big(u(x,t)\big)\big\}-\nabla_x
H_{F,u}(x,t)\bigg|^2\,dxdt+\int_{\R^N}\eta\big(u(x,0)\big)\,dx\,,
\end{equation}
Thus if there exists a solution to
\begin{equation}\label{bdgvfghffhkk}
\begin{cases}\e\Delta_x
\big\{\nabla_u \eta\big(u(x,t)\big)\big\}=\partial_t u(x,t)+div_x
F\big(u(x,t)\big)\quad \forall(x,t)\in\R^N\times(0,T),
\\ u(x,0)=v_0(x)\quad\quad\quad\quad\quad\forall x\in\R^N\,.
\end{cases}
\end{equation}
for some $v_0(x)\in L^2(\R^N,\R^k)\cap L^\infty$ then, by
\er{conslll}, $u(x,t)$ is also a minimizer to
\begin{equation}\label{bdgvfghffhkkjkgjkgjkgkgkgk}
\inf\big\{I_{\e,F}(u): u(x,0)=v_0(x)\big\}\,.
\end{equation}
Moreover, in this case,
\begin{equation}\label{bdgvfghffhkkjkgjkgjkgkgkgkjhnghj}
\inf\big\{I_{\e,F}(u):
u(x,0)=v_0(x)\big\}=\int_{\R^N}\eta\big(v_0(x)\big)\,dx\,,
\end{equation}
and the function $u(x,t):\R^N\times[0,T]\to\R^k$ is a minimizer to
\er{bdgvfghffhkkjkgjkgjkgkgkgk} if and only if $u(x,t)$ is a
solution to \er{bdgvfghffhkk}. Thus it is a natural question in the
Method of Vanishing Viscosity for Conservation Laws to know the
$\Gamma$-limit of the functional
\begin{equation}\label{consdfbhgyjykukukluklgjgh}
J_{\e,F,v_0}(u)=
\begin{cases}I_{\e,F}(u)\quad\quad
\text{if}\quad u(x,0)\equiv v_0(x),\\
+\infty\quad\quad\quad\;\text{otherwise}\,.\end{cases}
\end{equation}
\begin{lemma}\label{difdigdfighfghdfhgdhgh}
Consider $F\in C^1(\R^k,\R^{k\times N})$ satisfying $F(0)=0$. Next
let $u(x,t)\in BV\big(\R^N\times(0,T),\R^k\big)\cap
L^\infty\big(0,T;L^2(\R^N,\R^k)\big)\cap L^\infty$ be such that
$u(x,t)$ is continuous in $[0,T]$ as a function of $t$ with the
values in $L^\infty(\R^N,\R^k)$ with respect to $L^\infty$-weak$^*$
topology and satisfies the following Conservation Law on the strip:
\begin{equation}\label{conslaw}
\partial_t u(x,t)+div_x F\big(u(x,t)\big)=0\quad\quad\forall (x,t)\in\R^N\times(0,T)\,.
\end{equation}
Then we can construct $\bar v(x,t)\in
L^2_{loc}\big(\R^N\times(-2T,2T),\R^{k\times N}\big)$, such that if
we set $\bar u(x,t):=\Div_x\bar v(x,t)$ and $\bar
L(x,t):=-\partial_t\bar v(x,t)$, then $\bar u(x,t)\in
BV\big(\R^N\times(-2T,2T),\R^k\big)\cap
L^\infty\big(-2T,2T;L^2(\R^N,\R^k)\big)\cap L^\infty$, $\bar u(x,t)$
is continuous in $[-2T,2T]$ as a function of $t$ with the values in
$L^1(\R^N,\R^k)$ with respect to $L^1$-strong topology, $\bar
L(x,t)\in BV\big(\R^N\times(-2T,2T),\R^{k\times N}\big)\cap
L^\infty\big(-2T,2T;L^2(\R^N,\R^{k\times N})\big)\cap L^\infty$,
$\bar L(x,t)$ is continuous in $[-2T,2T]$ as a function of $t$ with
the values in $L^1(\R^N,\R^{k\times N})$ with respect to
$L^1$-strong topology, $\partial_t\bar u(x,t)+\Div_x \bar
L(x,t)\equiv 0$ in $\R^N\times(-2T,2T)$, and we have
\begin{equation}\label{fjhfhfutufghjfj}
\begin{cases}
\bar u(x,t)=u(x,t)\\
\bar L(x,t)=F\big(u(x,t)\big)\\
\bar u(x,-t)=2u(x,0)-u(x,t)
\\
\bar L(x,-t)=F\big(u(x,t)\big)
\end{cases}
\quad\forall (x,t)\in\R^N\times(0,T).
\end{equation}
\end{lemma}
\begin{proof}
Define $v(x,t):\R^N\times[0,T]\to\R^{k\times N}$ by
\begin{equation}\label{xfgvfhgjgjykuikukyuyuuyuyiichgvdfhjhjhnjhhjkkk}
v(x,t):=K(x)-\int_0^t
F\big(u(x,s)\big)ds\quad\quad\quad\quad\forall(x,t)\in\R^N\times\R\,,
\end{equation}
where $K(x)\in L^2_{loc}(\R^N,\R^{k\times N})$ satisfies $\Div_x
K(x)\equiv u(x,0)$.
Then using \er{conslaw} we deduce that
\begin{equation}\label{xfgvfhgjgjykuikukyuyuuyuyiichgvdfhjhhtjtyvvjhjjhjjbjhjhjjgkkk}
\partial_t v(x,t):=-
F\big(u(x,t)\big)\quad\text{and}\quad \Div_x
v(x,t)=u(x,t)\quad\quad\forall(x,t)\in\R^N\times(0,T)\,.
\end{equation}
Next define  $\hat v(x,t):\R^N\times[0,2T]\to\R^{k\times N}$ by
\begin{equation}\label{fjhfhfutufghjfjfngnhgkgj}
\hat v(x,t):=\begin{cases}
v(x,t)\quad\forall(x,t)\in\R^N\times[0,T]\\
2v(x,T)-v(x,2T-t)\quad\forall(x,t)\in\R^N\times[T,2T],
\end{cases}
\end{equation}
and set $\hat u(x,t):=\Div_x\hat v(x,t)$ and $\hat
L(x,t):=-\partial_t \hat v(x,t)$. Then, clearly $\hat u(x,t)\in
BV\big(\R^N\times(0,2T),\R^k\big)\cap
L^\infty\big(0,2T;L^2(\R^N,\R^k)\big)\cap L^\infty$, $\hat u(x,t)$
is continuous in $[0,2T]$ as a function of $t$ with the values in
$L^1(\R^N,\R^k)$ with respect to $L^1$-strong topology, $\hat
L(x,t)\in BV\big(\R^N\times(0,2T),\R^{k\times N}\big)\cap
L^\infty\big(0,2T;L^2(\R^N,\R^{k\times N})\big)\cap L^\infty$, $\hat
L(x,t)$ is continuous in $[0,2T]$ as a function of $t$ with the
values in $L^1(\R^N,\R^{k\times N})$ with respect to $L^1$-strong
topology, $\partial_t\hat u(x,t)+\Div_x \hat L(x,t)\equiv 0$ in
$\R^N\times(0,2T)$, and by
\er{xfgvfhgjgjykuikukyuyuuyuyiichgvdfhjhhtjtyvvjhjjhjjbjhjhjjgkkk}
we have
\begin{equation}\label{fjhfhfutufghjfjjgkg}
\begin{cases}
\hat u(x,t)=u(x,t)\\
\hat L(x,t)=F\big(u(x,t)\big)\\
\hat u(x,2T-t)=2u(x,T)-u(x,t)
\\
\hat L(x,2T-t)=F\big(u(x,t)\big)
\end{cases}
\quad\forall (x,t)\in\R^N\times(0,T).
\end{equation}
Next define  $\bar v(x,t):\R^N\times[-2T,2T]\to\R^{k\times N}$ by
\begin{equation}\label{fjhfhfutufghjfjfngnhgkgjjkggjk}
\bar v(x,t):=\begin{cases}
\hat v(x,t)\quad\forall(x,t)\in\R^N\times[0,2T]\\
2v(x,0)-\hat v(x,-t)\quad\forall(x,t)\in\R^N\times[-2T,0],
\end{cases}
\end{equation}
and set $\bar u(x,t):=\Div_x\bar v(x,t)$ and $\bar
L(x,t):=-\partial_t \bar v(x,t)$. Then, clearly $\bar u(x,t)\in
BV\big(\R^N\times(-2T,2T),\R^k\big)\cap
L^\infty\big(-2T,2T;L^2(\R^N,\R^k)\big)\cap L^\infty$, $\bar u(x,t)$
is continuous in $[-2T,2T]$ as a function of $t$ with the values in
$L^1(\R^N,\R^k)$ with respect to $L^1$-strong topology, $\bar
L(x,t)\in BV\big(\R^N\times(-2T,2T),\R^{k\times N}\big)\cap
L^\infty\big(-2T,2T;L^2(\R^N,\R^{k\times N})\big)\cap L^\infty$,
$\bar L(x,t)$ is continuous in $[-2T,2T]$ as a function of $t$ with
the values in $L^1(\R^N,\R^{k\times N})$ with respect to
$L^1$-strong topology, $\partial_t\bar u(x,t)+\Div_x \bar
L(x,t)\equiv 0$ in $\R^N\times(-2T,2T)$, and by
\er{fjhfhfutufghjfjjgkg} we deduce \er{fjhfhfutufghjfj}.
\end{proof}
\begin{lemma}\label{vjkjkjhnjhjhnjnjkbhjj}
Let $F(u)\in C^1(\R^k,\R^{k\times N})$ satisfying $F(0)=0$ and
$\eta(u)\in C^3(\R^k,\R)$ be an entropy for the corresponding system
\er{equ24}, which satisfies $\eta(u)\geq 0$ and $\eta(0)=0$.
Furthermore, let $u(x,t)$ be same as in Lemma
\ref{difdigdfighfghdfhgdhgh}.
and $\kappa(h)\in C_c^\infty(\R^{N+1})$ be a radial function, such
that $\int_{\R^{N+1}}\kappa(h)dh=1$. Then for every $\delta>0$ there
exists a sequence of functions $\big\{v_\e(x,t)\big\}_{\e>0}\in
C^\infty(\R^N\times\R,\R^{k\times N})$ such that $u_\e(x,t):=\Div_x
v_\e(x,t)\in W^{1,2}\big(\R^N\times(0,T),\R^k\big)\cap L^\infty$ and
$L_\e(x,t):=-\partial_t v_\e(x,t)\in
W^{1,2}\big(\R^N\times(0,T),\R^{k\times N}\big)\cap L^\infty$;
$\{u_\e\}$, $\{L_\e\}$ and $\{\e\nabla_x u_\e\}$ are bounded in
$L^\infty$ sequences; $u_\e\to u$, $L_\e\to F(u)$ and $\e \nabla_x
u_\e\to 0$ as $\e\to 0^+$ in $L^q \big(\R^N\times(0,T)\big)$;
$\partial_t u_\e+\Div_x L_\e\equiv 0$ and
\begin{multline}\label{glavnkjjghjbmainmainthfuyytufrt}
\lim_{\e\to
0^+}\int_0^T\int_{\R^N}\Bigg(\e\Big|\nabla_x\big\{\nabla_u
\eta\big(u_\e(x,t)\big)\big\}\Big|^2+\frac{1}{\e}\Big|\nabla_x H_{F,u_\e}(x,t)\Big|^2\Bigg)\,dxdt\leq\\
\lim_{\e\to
0^+}\int_0^T\int_{\R^N}\Bigg(\e\Big|\nabla_x\big\{\nabla_u
\eta\big(u_\e(x,t)\big)\big\}\Big|^2+\frac{1}{\e}\Big|
L_\e(x,t)-F\big(u_\e(x,t)\big)\Big|^2\Bigg)\,dxdt=\\
\lim_{\e\to 0^+}\int_0^T\int_{\R^N}\Bigg(\e\bigg|\nabla_x
\Big\{\nabla_u\eta\big(\Div_x
v_\e(x,t)\big)\Big\}\bigg|^2+\frac{1}{\e}\Big|
\partial_t v_\e(x,t)+F\big(\Div_x
v_\e(x,t)\big)\Big|^2\Bigg)\,dxdt\\
\leq\int_{J_{u}}\hat E_0\Big(u^+(x,t),
u^-(x,t),\vec\nu(x,t)\Big)\,\partial\mathcal{H}^N(x,t)+\delta\,,
\end{multline}
where $H_{F,u_\e}(x,t)\in L^2\big(0,T;\tilde H^1_0(\R^N,\R^k)\big)$
satisfies
\begin{equation}\label{conssstcghfghkjgbj}
\Delta_x H_{F,u_\e}(x,t)=\partial_t u_\e(x,t)+\Div_x
F\big(u_\e(x,t)\big)\,,
\end{equation}
\begin{multline}
\label{L2009limew03zeta71288888Contggiuuggyyyynew88789999vprop78899shtrihkkkllyhjyukjkkmmmklklklhhhhkkffgghhjjjkkkllkkhhhjjuhhiuijk}
\hat E_0\big(u^+, u^-,\vec\nu\big):=\inf\Bigg\{\int_{\bar I_{\vec
\nu}}\bigg(L\Big|\nabla_y
\Big\{\nabla_u\eta\big(\zeta(y,s)\big)\Big\}\Big|^2+\frac{1}{L}\Big|
\gamma(y,s)-F\big(\zeta(y,s)\big)\Big|^2\bigg)\,dyds:\\
L\in(0,+\infty),\,\zeta\in\mathcal{Z}^{(2)}\big(u^+,u^-,\vec\nu\big),\;
\gamma\in\mathcal{Z}^{(3)}\big(F(u^+),F(u^-),\vec\nu\big),\;\partial_s\zeta(y,s)+\Div_y\gamma(y,s)\equiv
0\Bigg\}=\\ \hat E_1\big(u^+,
u^-,\vec\nu\big):=\inf\Bigg\{\int_{\bar I_{\vec
\nu}}\bigg(L\Big|\nabla_y \Big\{\nabla_u\eta\big(\Div_y
\xi(y,s)\big)\Big\}\Big|^2+\frac{1}{L}\Big|
\partial_s\xi(y,s)+F\big(\Div_y
\xi(y,s)\big)\Big|^2\bigg)\,dyds:\\
L\in(0,+\infty),\,\xi\in\mathcal{Z}^{(1)}\big(u^+,u^-,\vec\nu\big)\Bigg\}
\,,
\end{multline}
with
\begin{multline}\label{L2009Ddef2hhhjjjj77788hhhkkkkllkjjjjkkkhhhhffggdddkkkgjhikhhhjjhhhhgbjhjn}
\mathcal{Z}^{(1)}\big(u^+,u^-,\vec\nu\big):=\\ \Bigg\{\xi(y,s)\in
\mathcal{D}'(\R^N\times\R,\R^{k\times N}):\;\Div_y\xi(y,s)\in
C^1(\R^N\times\R,\R^{k}),\;\partial_s\xi(y,s)\in
C^1(\R^N\times\R,\R^{k\times N}),\\
\big(\Div_y\xi,-\partial_s\xi\big)(y,s)=\big(u^-,F(u^-)\big)\;\text{
if }\;y\cdot\vec\nu\leq-1/2,\;
\big(\Div_y\xi,-\partial_s\xi\big)(y,s)=\big(u^+,F(u^+)\big)\;\text{
if }\; y\cdot\vec\nu\geq 1/2\\ \text{ and
}\;\big(\Div_y\xi,-\partial_s\xi\big)\big((y,s)+\vec\nu_j\big)=\big(\Div_y\xi,-\partial_s\xi\big)(y,s)\;\;\forall
j=2,\ldots, (N+1)\Bigg\}\,,
\end{multline}
\begin{multline}\label{L2009Ddef2hhhjjjj77788hhhkkkkllkjjjjkkkhhhhffggdddkkkgjhikhhhjjhhhhgbjhjnjhhuuuuik}
\mathcal{Z}^{(2)}\big(u^+,u^-,\vec\nu\big):= \Bigg\{\zeta(y,s)\in
C^1(\R^N\times\R,\R^{k}):\; \zeta(y,s)=u^-\;\text{ if
}\;y\cdot\vec\nu\leq-1/2,\\
\zeta(y,s)=u^+\;\text{ if }\; y\cdot\vec\nu\geq 1/2\; \text{ and
}\;\zeta\big((y,s)+\vec\nu_j\big)=\zeta(y,s)\;\;\forall j=2,\ldots,
(N+1)\Bigg\}\,,
\end{multline}
\begin{multline}\label{L2009Ddef2hhhjjjj77788hhhkkkkllkjjjjkkkhhhhffggdddkkkgjhikhhhjjhhhhgbjhjnijiuhujijikkkk}
\mathcal{Z}^{(3)}\big(A,B,\vec\nu\big):= \Bigg\{\gamma(y,s)\in
C^1(\R^N\times\R,\R^{k\times N}):\;
\gamma(y,s)=B\;\text{ if }\;y\cdot\vec\nu\leq-1/2,\\
\gamma(y,s)=A\;\text{ if }\; y\cdot\vec\nu\geq 1/2\; \text{ and
}\;\gamma\big((y,s)+\vec\nu_j\big)=\gamma(y,s)\;\;\forall
j=2,\ldots, (N+1)\Bigg\}\,.
\end{multline}
Here $\bar I_{\vec \nu}:=\big\{y\in\R^{N+1}:\;|y\cdot
\vec\nu_j|<1/2\;\;\;\forall j=1,\ldots, (N+1)\big\}$ where
$\{\vec\nu_1,
\ldots,\vec\nu_N,
\vec\nu_{N+1}\}\subset\R^{N+1}$ is an orthonormal base in $\R^{N+1}$
such that $\vec\nu_1:=\vec \nu$. Moreover, there exist $\sigma>0$
and $R>0$ (depending on $\delta$), such that for every $0<\e<1$ and
for every $(x,t)\in\big(\R^N\times\R\big)\setminus
\big(\{x\in\R^N:|x|<R\}\times(\sigma, T-\sigma)\big)$ we have
$u_\e(x,t)=u^{(0)}_\e(x,t)$ and $L_\e(x,t)=L^{(0)}_\e(x,t)$ where
\begin{equation*}
u^{(0)}_\e(x,t)=\frac{1}{\e^{N+1}}\int_\R\int_{\R^N}\kappa\bigg(\frac{y-x}{\e},\frac{s-t}{\e}\bigg)\bar
u(y,s)dyds,\;\;
L^{(0)}_\e(x,t)=\frac{1}{\e^{N+1}}\int_\R\int_{\R^N}\kappa\bigg(\frac{y-x}{\e},\frac{s-t}{\e}\bigg)\bar
L(y,s)dyds,
\end{equation*}
where $\bar u$ and $\bar L$ be the same as in Lemma
\ref{difdigdfighfghdfhgdhgh}.
\end{lemma}
\begin{proof}
Let $\bar v$, $\bar u$ and $\bar L$ be as in Lemma
\ref{difdigdfighfghdfhgdhgh}. In particular,
$$\big\|D_{(x,t)}\bar
u\big\|\Big(\big\{(x,t)\in\R^N\times\R:\,t(T-t)=0\big\}\Big)=\big\|D_{(x,t)}\bar
L\big\|\Big(\big\{(x,t)\in\R^N\times\R:\,t(T-t)=0\big\}\Big)=0.$$
Thus, using Theorem 4.1
in \cite{PI} we deduce that for every
$\delta>0$ there exists a sequences of functions
$\{u_\e(x,t)\}_{e>0}\subset C^\infty(\R^N\times\R,\R^{k})\cap
W^{1,2}\big(\R^N\times(0,T),\R^k\big)\cap L^\infty$ and
$\{L_\e(x,t)\}_{e>0}\subset C^\infty(\R^N\times\R,\R^{k\times
N})\cap W^{1,2}\big(\R^N\times(0,T),\R^{k\times N}\big)\cap
L^\infty$; such that $\{u_\e\}$, $\{L_\e\}$ and $\{\e\nabla_x
u_\e\}$ are bounded in $L^\infty$; $u_\e\to u$, $L_\e\to F(u)$ and
$\e \nabla_x u_\e\to 0$ as $\e\to 0^+$ in $L^q
\big(\R^N\times(0,T)\big)$;
$\partial_t u_\e+\Div_x L_\e\equiv 0$ and
\begin{multline}\label{glavnkjjghjbmainmainthfuyytufrtghkgjghgffghgg}
\lim_{\e\to
0^+}\int_0^T\int_{\R^N}\Bigg(\e\Big|\nabla_x\big\{\nabla_u
\eta\big(u_\e(x,t)\big)\big\}\Big|^2+\frac{1}{\e}\Big|
L_\e(x,t)-F\big(u_\e(x,t)\big)\Big|^2\Bigg)\,dxdt\\=\lim_{\e\to
0^+}\int_0^T\int_{\R^N}\frac{1}{\e}\Bigg(\Big|\nabla^2_u
\eta\big(u_\e(x,t)\big)\cdot\big(\e\nabla_x
u_\e(x,t)\big)\Big|^2+\Big|
L_\e(x,t)-F\big(u_\e(x,t)\big)\Big|^2\Bigg)\,dxdt\\
\leq\int_{J_{u}}\hat E_0\Big(u^+(x,t),
u^-(x,t),\vec\nu(x,t)\Big)\,\partial\mathcal{H}^N(x,t)+\delta\,,
\end{multline}
where $\hat E_0\big(u^+, u^-,\vec\nu\big)$ is defined by
\er{L2009limew03zeta71288888Contggiuuggyyyynew88789999vprop78899shtrihkkkllyhjyukjkkmmmklklklhhhhkkffgghhjjjkkkllkkhhhjjuhhiuijk}.
%
%
%
%
Moreover, there exist $\sigma>0$ and $R>0$, such that for every
$0<\e<1$ and every $(x,t)\in\big(\R^N\times\R\big)\setminus
\big(\{x\in\R^N:|x|<R\}\times(\sigma, T-\sigma)\big)$ we have
$u_\e(x,t)=u^{(0)}_\e(x,t)$ and $L_\e(x,t)=L^{(0)}_\e(x,t)$.
Furthermore, by Lemma \ref{gbhfhgdfdghfddhg}, or by
\er{fhfhjffhjhfgyukuoijkkjkkkhhhjjjhhgjhhyjtjkhjkjkjkvncccxxcxcdf},
we obtain
\begin{multline}\label{glavnkjjghjbmainmainthfuyytufrtghkgjghgffghghjfgfg}
\lim_{\e\to
0^+}\int_0^T\int_{\R^N}\Bigg(\e\Big|\nabla_x\big\{\nabla_u
\eta\big(u_\e(x,t)\big)\big\}\Big|^2+\frac{1}{\e}\Big|\nabla_x H_{F,u_\e}(x,t)\Big|^2\Bigg)\,dxdt\leq\\
\lim_{\e\to
0^+}\int_0^T\int_{\R^N}\Bigg(\e\Big|\nabla_x\big\{\nabla_u
\eta\big(u_\e(x,t)\big)\big\}\Big|^2+\frac{1}{\e}\Big|
L_\e(x,t)-F\big(u_\e(x,t)\big)\Big|^2\Bigg)\,dxdt=\\
\leq\int_{J_{u}}\hat E_0\Big(u^+(x,t),
u^-(x,t),\vec\nu(x,t)\Big)\,\partial\mathcal{H}^N(x,t)+\delta\,.
\end{multline}
Next, define
\begin{equation}\label{xfgvfhgjgjykuikukyuyuuyuyiichgvdfhjhjhnjhhj}
v_\e(x,t):=\bar K_\e(x)-\int_0^t
L_\e(x,s)\,ds\quad\quad\quad\quad\forall(x,t)\in\R^N\times\R\,,
\end{equation}
where $\bar K_\e(x)\in C^\infty(\R^N,\R^{k\times N})$ satisfies
$\Div_x \bar K_\e(x)\equiv u_\e(x,0)$. Then clearly $v_\e(x,t)\in
C^\infty(\R^N\times\R,\R^{k\times N})$. Moreover, since $\Div_x
L_{\e}(x,t)+\partial_t u_\e(x,t)\equiv 0$ we deduce that
\begin{equation}\label{xfgvfhgjgjykuikukyuyuuyuyiichgvdfhjhhtjtyvvjhjjhjjbjhjhjjg}
\partial_t v_\e(x,t):=-
L_{\e}(x,t)\quad\text{and}\quad \Div_x
v_\e(x,t)=u_\e(x,t)\quad\quad\forall(x,t)\in\R^N\times\R\,.
\end{equation}
Next we will prove that
\begin{equation}
\label{L2009limew03zeta71288888Contggiuuggyyyynew88789999vprop78899shtrihkkkllyhjyukjkkmmmklklklhhhhkkffgghhjjjkkkllkkhhhjjuhhiuijkghghgjkjkjkjkpppgyygyttf}
\hat E_1\big(u^+, u^-,\vec\nu\big)=\hat E_0\big(u^+,
u^-,\vec\nu\big)\,,
\end{equation}
where $\hat E_1\big(u^+, u^-,\vec\nu\big)$ and $\hat E_0\big(u^+,
u^-,\vec\nu\big)$ are defined by
\er{L2009limew03zeta71288888Contggiuuggyyyynew88789999vprop78899shtrihkkkllyhjyukjkkmmmklklklhhhhkkffgghhjjjkkkllkkhhhjjuhhiuijk}.
%
%
%
%
Indeed, since for every
$\xi\in\mathcal{Z}^{(1)}\big(u^+,u^-,\vec\nu\big)$ we clearly have
$\Div_y \xi\in\mathcal{Z}^{(2)}\big(u^+,u^-,\vec\nu\big)$ and
$$-\partial_s\xi\in\mathcal{Z}^{(3)}\big(F(u^+),F(u^-),\vec\nu\big)$$
and since $\partial_s(\Div_y \xi)+\Div_y(-\partial_s\xi)$, we
clearly have
\begin{equation}
\label{L2009limew03zeta71288888Contggiuuggyyyynew88789999vprop78899shtrihkkkllyhjyukjkkmmmklklklhhhhkkffgghhjjjkkkllkkhhhjjuhhiuijkghghgjkjkjkjk}
\hat E_1\big(u^+, u^-,\vec\nu\big)
\geq\hat E_0\big(u^+, u^-,\vec\nu\big)
\,.
\end{equation}
On the other hand fix
$\zeta\in\mathcal{Z}^{(2)}\big(u^+,u^-,\vec\nu\big)$ and
$\gamma\in\mathcal{Z}^{(3)}\big(F(u^+),F(u^-),\vec\nu\big)$ such
that $\partial_s\zeta(y,s)+\Div_y\gamma(y,s)\equiv 0$. Then define
\begin{equation}\label{xfgvfhgjgjykuikukyuyuuyuyiichgvdfhjhjhnjhhjhhjhbkl}
\xi(y,s):=
Q(y)-\int_{0}^{s}\gamma(y,\tau)\,d\tau\quad\quad\quad\quad\forall(y,s)\in\R^N\times\R\,,
\end{equation}
where $Q(y)\in C^1(\R^N,\R^k)$ is an arbitrary function which
satisfies $\Div_y Q(y)=\zeta(y,0)$. Then clearly $\xi\in
C^1(\R^N\times\R,\R^{k\times N})$ and moreover, since
$\partial_s\zeta(y,s)+\Div_y\gamma(y,s)\equiv 0$, we easily deduce
that
\begin{equation}\label{xfgvfhgjgjykuikukyuyuuyuyiichgvdfhjhhtjtyvvjhjjhjjbjhjhjjgjklhhjjhbbnbnhh}
\partial_s \xi(y,s):=-
\gamma(y,s)\quad\text{and}\quad \Div_y
\xi(y,s)=\zeta(y,s)\quad\quad\forall(y,s)\in\R^N\times\R\,.
\end{equation}
Thus clearly $\hat E_1\big(u^+, u^-,\vec\nu\big) \leq\hat
E_0\big(u^+, u^-,\vec\nu\big)$
and plugging it into
\er{L2009limew03zeta71288888Contggiuuggyyyynew88789999vprop78899shtrihkkkllyhjyukjkkmmmklklklhhhhkkffgghhjjjkkkllkkhhhjjuhhiuijkghghgjkjkjkjk}
we deduce
\er{L2009limew03zeta71288888Contggiuuggyyyynew88789999vprop78899shtrihkkkllyhjyukjkkmmmklklklhhhhkkffgghhjjjkkkllkkhhhjjuhhiuijkghghgjkjkjkjkpppgyygyttf}.
\end{proof}

\begin{proof}[Proof of Theorem \ref{prcnlkkmainthggg}]
Let $\delta>0$ and $\kappa(r)\in C_c^\infty(\R^{N+1},\R)$ be a
radial function, such that $\int_{\R^{N+1}}\kappa(r)dr=1$ and
$\kappa\geq 0$.
Therefore, by Lemma \ref{vjkjkjhnjhjhnjnjkbhjj}, there exists a
sequence of functions $\big\{v_{\e}(x,t)\big\}_{\e>0}\in
C^\infty(\R^N\times\R,\R^{k\times N})$ such that
$u_{\e}(x,t):=\Div_x v_{\e}(x,t)\in
W^{1,2}\big(\R^N\times(0,T),\R^k\big)\cap L^\infty$ and
$L_{\e}(x,t):=-\partial_t v_{\e}(x,t)\in
W^{1,2}\big(\R^N\times(0,T),\R^{k\times N}\big)\cap L^\infty$;
$\{u_{\e}\}$, $\{L_{\e}\}$ and $\{\e\nabla_x u_{\e}\}$ are bounded
in $L^\infty$ sequences; $u_{\e}\to u$, $L_{\e}\to F(u)$ and $\e
\nabla_x u_{\e}\to 0$ as $\e\to 0^+$ in $L^q
\big(\R^N\times(0,T)\big)$;
$\partial_t u_{\e}+\Div_x L_{\e}\equiv 0$ and
\begin{multline}\label{glavnkjjghjbmainmainthfuyytufrthuihuiuuhykkk}
\lim_{\e\to
0^+}\int_0^T\int_{\R^N}\Bigg(\e\Big|\nabla_x\big\{\nabla_u
\eta\big(u_{\e}(x,t)\big)\big\}\Big|^2+\frac{1}{\e}\Big|\nabla_x H_{F,u_{\e}}(x,t)\Big|^2\Bigg)\,dxdt\leq\\
\lim_{\e\to
0^+}\int_0^T\int_{\R^N}\Bigg(\e\Big|\nabla_x\big\{\nabla_u
\eta\big(u_{\e}(x,t)\big)\big\}\Big|^2+\frac{1}{\e}\Big|
L_{\e}(x,t)-F\big(u_{\e}(x,t)\big)\Big|^2\Bigg)\,dxdt=\\
\lim_{\e\to 0^+}\int_0^T\int_{\R^N}\Bigg(\e\bigg|\nabla_x
\Big\{\nabla_u\eta\big(\Div_x
v_{\e}(x,t)\big)\Big\}\bigg|^2+\frac{1}{\e}\Big|
\partial_t v_{\e}(x,t)+F\big(\Div_x
v_{\e}(x,t)\big)\Big|^2\Bigg)\,dxdt\\
\leq\int_{J_{u}}\hat E_0\Big(u^+(x,t),
u^-(x,t),\vec\nu(x,t)\Big)\,\partial\mathcal{H}^N(x,t)+\delta\,,
\end{multline}
where $H_{F,u_{\e}}(x,t)\in L^2\big(0,T;\tilde
H^1_0(\R^N,\R^k)\big)$ satisfies
\begin{equation}\label{conssstcghfghkjgbjkjjijiojoiio}
\Delta_x H_{F,u_{\e}}(x,t)=\partial_t u_{\e}(x,t)+\Div_x
F\big(u_{h,\e}(x,t)\big)\,,
\end{equation}
Moreover, there exist $\sigma>0$ such that for every $0<\e<1$ and
every $(x,t)\in\big(\R^N\times\R\big)\setminus
\big(\{\R^N\}\times(\sigma, T-\sigma)\big)$ we have
$u_{\e}(x,t)=u^{(0)}_{\e}(x,t)$
where
\begin{equation}\label{gjkjgjghiuilkk}
u^{(0)}_{\e}(x,t)=\frac{1}{\e^{N+1}}\int_\R\int_{\R^N}\kappa\bigg(\frac{y-x}{\e},\frac{s-t}{\e}\bigg)\bar
u(y,s)dyds
\end{equation}
where $\bar u,\bar L$ are the same as in Lemma
\ref{difdigdfighfghdfhgdhgh}. We need just to modify $u_{\e}$
slightly in such a way that it will satisfy the condition
$u_{\e}(x,0)\equiv u(x,0)$.
Let $\chi_1(x,t)\in L^\infty\big(0,+\infty;L^2(\R^N,\R^k)\big)\cap
L^2(0,+\infty;\tilde H^1_0\big(\R^N,\R^k)\big)\cap L^\infty$ be the
solution of the
heat equation:
\begin{equation}\label{hitinuty87}
\begin{cases}
\Delta_x \chi_1=\partial_t \chi_1\,,\\ \chi_1(x,0)=u(x,0)\,.
\end{cases}
\end{equation}
and set $\chi_\e(x,t):=\chi_1(x,\e t)$. Then $\chi_\e(x,t)\in
L^\infty\big(0,T;L^2(\R^N,\R^k)\big)\cap L^2(0,T;\tilde
H^1_0\big(\R^N,\R^k)\big)\cap L^\infty$ and satisfies:
\begin{equation}\label{hitin}
\begin{cases}
\e\Delta_x \chi_\e=\partial_t \chi_\e\,,\\ \chi_\e(x,0)=u(x,0)\,.
\end{cases}
\end{equation}
It is clear that we may assume that $\chi_\e$ is $L^2$-strongly
continuous in $[0,T]$ as a function of $t$ and
$\chi_\e(x,0)=u(x,0)$. Moreover for every $0\leq\bar t\leq T$ we
have
\begin{equation}\label{energy}
2\int_0^{\bar
t}\int_{\R^N}\e\big|\nabla\chi_\e(x,s)\big|^2dxds=\int_{\R^N}u^2(x,0)dx-\int_{\R^N}\chi^2_\e(x,\bar
t)dx\,.
\end{equation}
Finally, by the well known maximum principle for the Heat Equation
we clearly have $\|\chi_\e(\cdot,t)\|_{L^\infty}\leq
\|u(\cdot,0)\|_{L^\infty}$. Next let $\theta(t)\in
C^\infty(\R,[0,1])$ be a cut-off function satisfying $\theta(t)=0$
for every $t\geq 1$ and $\theta(t)=1$ for every $t\leq 1/2$. For
every small  $\e>0$ define $\bar u_{\e}(x,t)\in
L^\infty\big(0,T;L^2(\R^N,\R^k)\big)\cap L^2\big(0,T;\tilde
H^1_0(\R^N,\R^k)\big)\cap L^\infty$ by
\begin{equation}
\label{oprdigjojo} \bar u_{\e}(x,t):=u_{\e}(x,t)+\theta(t/\e
)\big(\chi_\e(x,t)-u_{\e}(x,0)\big) =u_{\e}(x,t)+\theta(t/\e)
\big(\chi_\e(x,t)-u^{(0)}_{\e}(x,0)\big) \,.
\end{equation}
Then $\bar u_{\e}$ is $L^2$-strongly continuous in $[0,T]$ as a
function of $t$ and $\bar u_{\e}(x,0)=u(x,0)$. Moreover, $\bar
u_{\e}(x,t)=u_{\e}(x,t)$ whenever $t\geq \e$. Finally
\begin{equation}\label{dop3kkk}
\big\|\bar
u_{\e}-u_{\e}\big\|_{L^\infty}+\big\|u^{(0)}_{\e}\big\|_{L^\infty}\leq
C_0,
\end{equation}
where $C_0>0$ is a constant which does not depend on $\e$.

 Now we want to prove that
\begin{equation}
\label{limsupkkkjghg} \limsup\limits_{\e\to
0^+}\int_0^{\e}\int_{\R^N}\bigg\{\e\big|\nabla_x(\bar u_{\e}
-u_{\e})\big|^2+\frac{1}{\e}\Big|\nabla_x
D_{\e}(x,t)\Big|^2\bigg\}\,dxdt=0\,,
\end{equation}
where $D_{\e}(x,t)\in L^2\big(0,T;\tilde H^1_0(\R^N,\R^k)\big)$
satisfies
$\Delta_x D_{\e}(x,t)=\Big\{
\partial_t \big(\bar u_{\e}
-u_{\e}\big)+\Div_x \big(F(\bar u_{\e})-F(u_{\e})\big)\Big\}$,

 First of all by \er{oprdigjojo} and \er{energy} we observe that
\begin{multline}\label{limsupkkkjghggrad}
\limsup\limits_{\e\to 0^+}\int_0^{\e}\int_{\R^N}\e\big|\nabla_x(\bar
u_{\e} -u_{\e})\big|^2\,dxdt\leq \limsup\limits_{\e\to 0^+}2\e^2
\int_{\R^N}\big|\nabla_x u^{(0)}_{\e}(x,0)\big|^2\,dx+\\
\limsup\limits_{\e\to 0^+}\int_0^{\e}\int_{\R^N}2\e\big|\nabla_x
\chi_\e(x,t)\big|^2\,dxdt= \limsup\limits_{\e\to 0^+}2
\int_{\R^N}\bigg|\int_\R\int_{\R^N}\nabla_z\kappa(z,s)\otimes \bar
u(x+\e z,\e s)dzds\bigg|^2\,dx\\+\limsup\limits_{\e\to
0^+}\bigg(\int_{\R^N}\chi^2_\e(x,0)dx-\int_{\R^N}\chi^2_\e(x,\e)dx\bigg)=
\limsup\limits_{\e\to 0^+}\bigg(\int_{\R^N}
\chi^2_1(x,0)dx-\int_{\R^N}\chi^2_1(x,\e^2)dx\bigg) =0\,.
\end{multline}
%
%
%
%
%
%
%
%
On the other hand by \er{dop3kkk} we deduce
\begin{multline}\label{limsupkkkjghgvtch}
\limsup\limits_{\e\to
0^+}\int_0^{\e}\int_{\R^N}\frac{1}{\e}\Big|\nabla_x
P_{\e}(x,t)\Big|^2\,dxdt\leq\limsup\limits_{\e\to
0^+}\int_0^{\e}\int_{\R^N}\frac{1}{\e}\big|F(\bar u_{\e})
-F(u_{\e})\big|^2\,dxdt\leq C_1\limsup\limits_{\e\to
0^+}\int_0^{\e}\int_{\R^N}\frac{1}{\e}|\bar u_{\e} -u_{\e}|^2\,dxdt
\\ \leq C_1\limsup\limits_{\e\to
0^+}\int_0^{\e}\int_{\R^N}\frac{2}{\e}\Big(\big|\chi_\e(x,t)-\chi_\e(x,0)\big|^2+\big|u^{(0)}_{\e}(x,0)-u(x,0)\big|^2\Big)\,dxdt=\\
C_1\limsup\limits_{\e\to
0^+}\bigg(2\int_0^{1}\int_{\R^N}\Big(\big|\chi_1(x,\e^2
\tau)-\chi_1(x,0)\big|^2\,dxd\tau+\int_{\R^N}\big|u^{(0)}_{\e}(x,0)-u(x,0)\big|^2\Big)dx\bigg)
=0\,.
\end{multline}
where $P_{\e}(x,t)\in L^2\big(0,T;\tilde H^1_0(\R^N,\R^k)\big)$
satisfies
$\Delta_x P_{\e}(x,t)=\Div_x \big(F(\bar u_{\e})-F(u_{\e})\big)$.
Next, using \er{hitin} we infer
\begin{multline}\label{estdt}
\partial_t (\bar u_{\e}
-u_{\e})=\partial_t\Big\{\theta\big(t/\e\big)
\big(\chi_\e(x,t)-u^{(0)}_{\e}(x,0)\big)\Big\}=
\theta\big(t/\e\big)\partial_t\chi_\e(x,t)+
\e^{-1}\theta'\big(t/\e\big)
\big(\chi_\e(x,t)-u^{(0)}_{\e}(x,0)\big)\\=\theta\big(t/\e\big)\e\Delta_x\chi_\e(x,t)+
\e^{-1}\theta'\big(t/\e\big)
\big(\chi_\e(x,t)-u^{(0)}_{\e}(x,0)\big)\,.
\end{multline}
Then consider $Q_{\e}(x,t)\in L^2\big(0,T;\tilde
H^1_0(\R^N,\R^k)\big)$ such that
$\Delta_x Q_{\e}(x,t)=
\partial_t \big(\bar u_{\e}
-u_{\e}\big)$.
By \er{estdt} we obtain
\begin{equation}\label{conssstcghfghkjgbjggytffttfiuhuuhhijhihuiuihihufghjhkjhjjhjh}
\Delta_x Q_{\e}(x,t)=
\theta\big(t/\e\big)\e\Delta_x\chi_\e(x,t)+
\e^{-1}\theta'\big(t/\e\big)
\Big(\big(\chi_\e(x,t)-\chi_\e(x,0)\big)-\big(u^{(0)}_{\e}(x,0)-u(x,0)\big)\Big)\,.
\end{equation}
On the other hand using the fact that $\bar u(x,t)\equiv
2u(x,0)-\bar u(x,-t)$ (see Lemma \ref{difdigdfighfghdfhgdhgh}) and
the fact that $\kappa(\cdot,-s)=\kappa(\cdot,s)$, by
\er{gjkjgjghiuilkk} we deduce:
\begin{multline}\label{gkguyhiuhkgjjj}
u^{(0)}_{\e}(x,0)=
\frac{1}{\e^{N+1}}\int_{\R^N}\int_{-\infty}^{0}\kappa\bigg(\frac{y-x}{\e},\frac{s}{\e}\bigg)
\Big(2u(y,0)-\bar
u(y,-s)\Big)dsdy\\+\frac{1}{\e^{N+1}}\int_{\R^N}\int_{0}^{+\infty}\kappa\bigg(\frac{y-x}{\e},\frac{s}{\e}\bigg)\bar
u(y,s)dsdy=\frac{2}{\e^{N+1}}\int_{\R^N}\int_{0}^{+\infty}\kappa\bigg(\frac{y-x}{\e},\frac{s}{\e}\bigg)
u(y,0)dsdy=\\
\frac{1}{\e^{N}}\int_{\R^N}\bigg(2\int_{0}^{+\infty}\kappa\Big(\frac{y-x}{\e},\tau\Big)d\tau\bigg)
u(y,0)dy=\frac{1}{\e^{N}}\int_{\R^N}\omega\Big(\frac{y-x}{\e}\Big)
u(y,0)dy,
\end{multline}
where $\omega(z)\in C_c^\infty(\R^{N})$ satisfies
$\int_{\R^{N}}\omega(z)dz=1$, $\omega\geq 0$ and
$\omega(-z)\equiv\omega(z)$. Thus, for every $\rho\in(0,\e)$ we have
\begin{multline}
\label{aas5789aaa1gbghfhgnhjhbjhgh}
u^{(0)}_{\e}(x,0)-u^{(0)}_{\rho}(x,0)= \int_\rho^\e\frac{\partial
u^{(0)}_{\tau}(x,0)}{\partial\tau}\,d\tau =\int_\rho^\e
\frac{\partial}{\partial\tau}\bigg(\frac{1}{\tau^{N}}\int_{\R^N}\omega\Big(\frac{y-x}{\tau}\Big)u(y,0)dy\bigg)d\tau=\\
-\int_\rho^\e
\bigg(\frac{1}{\tau^{N+1}}\int_{\R^N}\bigg\{\frac{y-x}{\tau}\cdot\nabla\omega\Big(\frac{y-x}{\tau}\Big)
+N\omega\Big(\frac{y-x}{\tau}\Big)\bigg\}u(y,0)dy\bigg)d\tau=\\
\int_\rho^\e
\bigg(\frac{1}{\tau^{N}}\int_{\R^N}\bigg(\Div_x\bigg\{\frac{y-x}{\tau}\omega\Big(\frac{y-x}{\tau}\Big)\bigg\}
u(y,0)dy
\bigg)d\tau\\= \Div_x\int_\rho^\e
\bigg(\frac{1}{\tau^{N}}\int_{\R^N}
u(y,0)\otimes\bigg\{\frac{y-x}{\tau}\omega\Big(\frac{y-x}{\tau}\Big)\bigg\}dy\bigg)d\tau.
\end{multline}
Then, by
\er{aas5789aaa1gbghfhgnhjhbjhgh}, for every $\varphi(x)\in
C^\infty_c(\R^N,\R^k)$ we have
\begin{multline}
\label{aas5789aaa1gbghfhgnhjhbjhghbhjhhiuujhnj}
\int_{\R^N}\Big(u^{(0)}_{\e}(x,0)-u^{(0)}_{\rho}(x,0)\Big)\cdot\varphi(x)dx=\\
-\int_{\R^N}\int_\rho^\e \bigg(\frac{1}{\tau^{N}}\int_{\R^N}
u(y,0)\otimes\bigg\{\frac{y-x}{\tau}\omega\Big(\frac{y-x}{\tau}\Big)\bigg\}dy\bigg):\nabla\varphi(x)d\tau
dx.
\end{multline}
Thus letting $\rho\to 0^+$ in
\er{aas5789aaa1gbghfhgnhjhbjhghbhjhhiuujhnj}, we obtain
\begin{multline}
\label{aas5789aaa1gbghfhgnhjhbjhghbhjhhiuujhnjhhgguyjouh}
\int_{\R^N}\Big(u^{(0)}_{\e}(x,0)-u(x,0)\Big)\cdot\varphi(x)dx=
-\int_{\R^N}\int_0^\e \bigg(\int_{\R^N}\omega(z) u(x+\tau
z,0)\otimes zdz\bigg):\nabla\varphi(x)d\tau dx
=\\
-\e\int_{\R^N}\int_0^1 \bigg(\int_{\R^N}\omega(z) u(x+\e\tau
z,0)\otimes zdz\bigg):\nabla\varphi(x)d\tau dx.
\end{multline}
Thus
\begin{equation}
\label{aas5789aaa1gbghfhgnhjhbjhghbhjhhiuujhnjhhgguyjouhfhthjhh}
u^{(0)}_{\e}(x,0)-u(x,0)= \e\Div\bigg\{\int_0^1\int_{\R^N}\omega(z)
u(x+\e\tau z,0)\otimes zdzd\tau
\bigg\}\,,
\end{equation}
On the other hand by \er{hitin} we deduce
\begin{equation}\label{limsupkkkjghgdtheet}
\chi_\e(x,t)-\chi_\e(x,0)=\e\Delta_x\bigg(\int_0^t\chi_\e(x,s)ds\bigg)\,.
\end{equation}
Therefore, by
\er{aas5789aaa1gbghfhgnhjhbjhghbhjhhiuujhnjhhgguyjouhfhthjhh},
\er{limsupkkkjghgdtheet} and
\er{conssstcghfghkjgbjggytffttfiuhuuhhijhihuiuihihufghjhkjhjjhjh} we
obtain
\begin{multline}\label{conssstcghfghkjgbjggytffttfiuhuuhhijhihuiuihihufghjhkjhjjhjhnkjuhkjgjkjk}
\Delta_x Q_{\e}(x,t)=
\text{div}_x\Bigg\{\theta\big(t/\e\big)\e\nabla_x\chi_\e(x,t)+
\theta'\big(t/\e\big)
\nabla_x\bigg(\int_0^t\chi_\e(x,s)ds\bigg)\\
-\theta'\big(t/\e\big)\bigg(\int_0^1\int_{\R^N}\omega(z) u(x+\e\tau
z,0)\otimes zdzd\tau
\bigg)
\Bigg\}.
\end{multline}
Thus,
\begin{multline}\label{limsupkkkjghgdt}
\limsup\limits_{\e\to
0^+}\int_0^{\e}\int_{\R^N}\frac{1}{\e}\Big|\nabla_x
Q_{\e}(x,t)\Big|^2\,dxdt\leq C\limsup\limits_{\e\to
0^+}\int_0^{\e}\int_{\R^N}\e\big|\nabla_x\chi_\e(x,t)\big|^2dxdt+\\
\limsup\limits_{\e\to 0^+}\frac{C}{\e}\int_0^{\e}\int_{\R^N}\bigg|
\nabla_x\Big(\int_0^t\chi_\e(x,s)ds\Big)\bigg|^2dxdt+
C\limsup\limits_{\e\to
0^+}\int_{\R^N}\bigg|\int_0^1\int_{\R^N}\omega(z)  u(x+\e \tau
z,0)\otimes z dzd\tau
\bigg|^2dx.
\end{multline}
Therefore,
\begin{multline}\label{limsupkkkjghgdjkltj}
\limsup\limits_{\e\to
0^+}\int_0^{\e}\int_{\R^N}\frac{1}{\e}\Big|\nabla_x
Q_{\e}(x,t)\Big|^2\,dxdt \leq C\Bigg(\limsup\limits_{\e\to
0^+}\frac{1}{2}\bigg(\int_{\R^N}\chi^2_\e(x,0)dx-\int_{\R^N}\chi^2_\e(x,\e)dx\bigg)+
\\ \limsup\limits_{\e\to
0^+}\frac{1}{\e}\int_0^{\e}\int_{\R^N}\bigg|
\nabla_x\Big(\int_0^t\chi_\e(x,s)ds\Big)\bigg|^2dxdt
+\limsup\limits_{\e\to 0^+}\int_{\R^N}\bigg|
\int_0^1\int_{\R^N}\omega(z) u(x+\e\tau z,0)\otimes z dzd\tau
\bigg|^2dx\Bigg)
\\
=\limsup\limits_{\e\to 0^+}\frac{C}{\e}\int_0^{\e}\int_{\R^N}\bigg|
\nabla_x\Big(\int_0^t\chi_\e(x,s)ds\Big)\bigg|^2dxdt
+C\limsup\limits_{\e\to
0^+}\int_{\R^N}\bigg|\int_0^1\int_{\R^N}\omega(z) u(x+\e \tau
z,0)\otimes zdzd\tau\bigg|^2dx.
\end{multline}
On the other hand
\begin{multline}\label{limsupkkkjghgdtsled}
\limsup\limits_{\e\to 0^+}\frac{1}{\e}\int_0^{\e}\int_{\R^N}\bigg|
\nabla_x\Big(\int_0^t\chi_\e(x,s)ds\Big)\bigg|^2dxdt=\limsup\limits_{\e\to
0^+}\int_0^{\e}\int_{\R^N}\frac{1}{\e}\bigg|\int_0^t\nabla_x
\chi_\e(x,s)\,ds\bigg|^2\,dxdt\leq\\ \limsup\limits_{\e\to
0^+}\int_0^{\e}\int_{\R^N}\frac{t}{\e}\bigg(\int_0^t\big|\nabla_x
\chi_\e(x,s)\big|^2\,ds\bigg)dxdt\leq\limsup\limits_{\e\to
0^+}\int_0^{\e}\int_{\R^N}\e\big|\nabla_x \chi_\e\big|^2\,dxdt=\\
\limsup\limits_{\e\to
0^+}\frac{1}{2}\bigg(\int_{\R^N}\chi^2_\e(x,0)dx-\int_{\R^N}\chi^2_\e(x,\e)dx\bigg)
=\limsup\limits_{\e\to
0^+}\frac{1}{2}\bigg(\int_{\R^N}\chi^2_1(x,0)dx-\int_{\R^N}\chi^2_1(x,\e^2)dx\bigg)=0.
\end{multline}
Then, by \er{limsupkkkjghgdjkltj} and \er{limsupkkkjghgdtsled} we
obtain
\begin{equation}
\label{limsupkkkjghgdtvthiim} \limsup\limits_{\e\to
0^+}\int_0^{\e}\int_{\R^N}\frac{1}{\e}\Big|\nabla_x
Q_{\e}(x,t)\Big|^2dxdt \leq  C\limsup\limits_{\e\to
0^+}\int_{\R^N}\bigg|\int_0^1\int_{\R^N}\omega(z)z\otimes  u(x+\e
\tau z,0)dzd\tau\bigg|^2dx=0.
\end{equation}
Thus by \er{limsupkkkjghggrad}, \er{limsupkkkjghgvtch} and
\er{limsupkkkjghgdtvthiim} we deduce \er{limsupkkkjghg}. Therefore,
by \er{glavnkjjghjbmainmainthfuyytufrthuihuiuuhykkk} and
\er{limsupkkkjghg}
we obtain
\begin{multline}\label{glavnkjjghjbmainmainthfuyytufrthuihuiuuhykkkguftyfift}
\lim_{\e\to
0^+}\int_0^T\int_{\R^N}\Bigg(\e\Big|\nabla_x\big\{\nabla_u
\eta\big(\bar u_{\e}(x,t)\big)\big\}\Big|^2+\frac{1}{\e}\Big|\nabla_x H_{F,\bar u_{\e}}(x,t)\Big|^2\Bigg)\,dxdt\leq\\
\leq\int_{J_{u}}\hat E_0\Big(u^+(x,t),
u^-(x,t),\vec\nu(x,t)\Big)\,\partial\mathcal{H}^N(x,t)+\delta\,,
\end{multline}
Therefore, taking
$\bar L_\e:=-\nabla_x H_{F,\bar u_{\e}}+F(\bar u_{\e})$, by
\er{glavnkjjghjbmainmainthfuyytufrthuihuiuuhykkkguftyfift} and the
fact that $\bar u_{\e}(x,T)=u_{\e}(x,T)=u^{(0)}_{\e}(x,T)\to u(x,T)$
in $L^2(\R^N,\R^k)$, we deduce $\partial_t \bar u_{\e}+\Div_x \bar
L_\e\equiv 0$ and
\begin{multline}\label{glavnkjjghjbmainmainth}
\lim_{\e\to
0^+}\Bigg\{\int_0^T\int_{\R^N}\Bigg(\e\Big|\nabla_x\big\{\nabla_u
\eta\big(\bar u_\e(x,t)\big)\big\}\Big|^2+\frac{1}{\e}\Big|\nabla_x
H_{F,\bar u_\e}(x,t)\Big|^2\Bigg)\,dxdt
+\int_{\R^N}\eta\big(\bar u_\e(x,T)\big)\,dx\Bigg\}=\\
\lim_{\e\to
0^+}\Bigg\{\int_0^T\int_{\R^N}\Bigg(\e\Big|\nabla_x\big\{\nabla_u
\eta\big(\bar u_\e(x,t)\big)\big\}\Big|^2+\frac{1}{\e}\Big| \bar
L_\e(x,t)-F\big(\bar
u_\e(x,t)\big)\Big|^2\Bigg)\,dxdt+\int_{\R^N}\eta\big(\bar
u_\e(x,T)\big)\,dx\Bigg\}\\= \lim_{\e\to
0^+}\Bigg\{\int_0^T\int_{\R^N}\Bigg(\e\bigg|\nabla_x
\Big\{\nabla_u\eta\big(\Div_x \bar
v_\e(x,t)\big)\Big\}\bigg|^2+\frac{1}{\e}\Big|
\partial_t \bar v_\e(x,t)+F\big(\Div_x
\bar v_\e(x,t)\big)\Big|^2\Bigg)\,dxdt
\\+\int_{\R^N}\eta\big(\Div_x\bar v_\e(x,T)\big)\,dx\Bigg\}
\leq\int_{J_{u}}\hat E_0\Big(u^+(x,t),
u^-(x,t),\vec\nu(x,t)\Big)\,\partial\mathcal{H}^N(x,t)+\int_{\R^N}\eta\big(u(x,T)\big)\,dx+\delta,
\end{multline}
(where as before we can define the corresponding function $\bar
v_\e$).

Finally, the result follows by letting $\delta\to0^+$ and using a
diagonal argument.
\end{proof}
\begin{proof}[Proof of Theorem
\ref{prcnlkkmainthhuiuiggg}]
The inequality \er{glavnkjjghjbmainmainthgjhguyhgggg}
follows by Theorem 2.1
in \cite{PII} or by Theorem
2.3 in \cite{PII}.
%
%
%
%
\end{proof}

\appendix
\section{Notations and basic results about $BV$-functions}
\label{sec:pre3}

\noindent$\bullet$ For given a real topological linear space $X$ we
denote by $X^*$ the dual space (the space of continuous linear
functionals from $X$ to $\R$).

\noindent$\bullet$ For given $h\in X$ and $x^*\in X^*$ we denote by
$\big<h,x^*\big>_{X\times X^*}$ the value in $\R$ of the functional
$x^*$ on the vector $h$.

%
%
%
%

\noindent$\bullet$ Given open set $G\subset\R^N$ we denote by
$\mathcal{D}(G,\R^d)$  the real topological linear space of
compactly supported $\R^d$-valued test functions i.e.
$C^\infty_c(G,\R^d)$ with the usual topology.

\noindent$\bullet$ Denote
$\mathcal{D}'(G,\R^d):=\big\{\mathcal{D}(G,\R^d)\big\}^*$ (the space
of $\R^d$-valued distributions in $G$).

\noindent$\bullet$ Given $h\in\mathcal{D}'(G,\R^d)$ and
$\delta\in\mathcal{D}(G,\R^d)$ we denote the value in $\R$ of the
distribution $h$ on the test function $\delta$ by
$<\delta,h>:=\big<\delta,h\big>_{\mathcal{D}(G,\R^d)\times
\mathcal{D}'(G,\R^d)}$.
%
%

\noindent$\bullet$ For a $p\times q$ matrix $A$ with $ij$-th entry
$a_{ij}$ and for a $q\times d$ matrix $B$ with $ij$-th entry
$b_{ij}$ we denote by $AB:=A\cdot B$ their product, i.e. the
$p\times d$ matrix, with $ij$-th entry
$\sum\limits_{k=1}^{q}a_{ik}b_{kj}$.

\noindent$\bullet$ We identify a $\vec u=(u_1,\ldots,u_q)\in\R^q$
with the $q\times 1$ matrix having $i1$-th entry $u_i$, so that for
a $p\times q$ matrix $A$ with $ij$-th entry $a_{ij}$ and for $\vec
v=(v_1,v_2,\ldots,v_q)\in\R^q$ we denote by $A\,\vec v :=A\cdot\vec
v$ the $p$-dimensional vector $\vec u=(u_1,\ldots,u_p)\in\R^p$,
given by $u_i=\sum\limits_{k=1}^{q}a_{ik}v_k$ for every $1\leq i\leq
p$.

\noindent$\bullet$ As usual $A^T$ denotes the transpose of the
matrix $A$.

\noindent$\bullet$ For $\vec u=(u_1,\ldots,u_p)\in\R^p$ and $\vec
v=(v_1,\ldots,v_p)\in\R^p$ we denote by $\vec u\vec v:=\vec
u\cdot\vec v:=\sum\limits_{k=1}^{p}u_k v_k$ the standard scalar
product. We also note that $\vec u\vec v=\vec u^T\vec v=\vec v^T\vec
u$ as products of matrices.

\noindent$\bullet$ For $\vec u=(u_1,\ldots,u_p)\in\R^p$ and $\vec
v=(v_1,\ldots,v_q)\in\R^q$ we denote by $\vec u\otimes\vec v$ the
$p\times q$ matrix with $ij$-th entry $u_i v_j$.

\noindent$\bullet$ For any $p\times q$ matrix $A$ with $ij$-th entry
$a_{ij}$ and $\vec v=(v_1,v_2,\ldots,v_d)\in\R^d$ we denote by
$A\otimes\vec v$ the $p\times q\times d$ tensor with $ijk$-th entry
$a_{ij}v_k$.

\noindent$\bullet$ Given a vector valued function
$f(x)=\big(f_1(x),\ldots,f_k(x)\big):\O\to\R^k$ ($\O\subset\R^N$) we
denote by $Df$ or by $\nabla_x f$ the $k\times N$ matrix with
$ij$-th entry $\frac{\partial f_i}{\partial x_j}$.

\noindent$\bullet$ Given a matrix valued function
$F(x):=\{F_{ij}(x)\}:\R^N\to\R^{k\times N}$ ($\O\subset\R^N$), we
denote $div\,F:=(l_1,\ldots,l_k)\in\R^k$, where
$l_i=\sum\limits_{j=1}^{N}\frac{\partial F_{ij}}{\partial x_j}$.

\noindent$\bullet$ Given a matrix valued function
$F(x)=\big\{f_{ij}(x)\big\}(1\leq i\leq p,\,1\leq j\leq
q):\O\to\R^{p\times q}$ ($\O\subset\R^N$) we denote by $DF$ or by
$\nabla_x F$ the $p\times q\times N$ tensor with $ijk$-th entry
$\frac{\partial f_{ij}}{\partial x_k}$.

\noindent$\bullet$ Given a vector measure $\mu=(\mu_1,\ldots,\mu_k)$
(where $\forall j=1,\ldots,k$ $\mu_j$ is a finite signed measure) we
denote by $\|\mu\|(E)$ the total variation of $\mu$ on the set $E$.

\noindent$\bullet$ For any $\mu$-measurable function $f$, we define
the product measure $f\cdot\mu$ by: $f\cdot\mu(E)=\int_E f\,d\mu$,
for every $\mu$-measurable set $E$.

 In what follows we present some known results on BV-spaces. We rely mainly on the book \cite{amb}
by Ambrosio, Fusco and Pallara.
%
%
%
%
%
%
\begin{definition}
Let $\Omega$ be a domain in $\R^N$ and let $f\in L^1(\Omega,\R^m)$.
We say that $f\in BV(\Omega,\R^m)$ if the following quantity is
finite:
\begin{equation*}
\int_\Omega|Df|:= \sup\bigg\{\int_\Omega f\cdot\Div\f \,dx :\,\f\in
C^1_c(\Omega,\R^{m\times N}),\;|\f(x)|\leq 1\;\forall x \bigg\}.
\end{equation*}
\end{definition}
\begin{definition}\label{defjac889878}
Let $\Omega$ be a domain in $\R^N$. Consider a function
$f\in L^1_{loc}(\Omega,\R^m)$ and a point $x\in\Omega$.\\
i) We say that $x$ is an {\em approximate continuity point} of $f$
if there exists $z\in\R^m$ such that
$$\lim\limits_{\rho\to
0^+}\frac{\int_{B_\rho(x)}|f(y)-z|\,dy} {\rho^N}=0.$$ In this case
we denote $z$ by $\tilde{f}(x)$. The set of  approximate continuity
points of
$f$ is denoted by $G_f$.\\
ii) We say that $x$ is an {\em approximate jump point} of $f$ if
there exist $a,b\in\R^m$ and $\vec\nu\in S^{N-1}$ such that $a\neq
b$ and
\begin{equation*}
\lim\limits_{\rho\to
0^+}\frac{\int_{B_\rho(x)}\big|\,f(y)-\chi(a,b,\vec\nu)(y)\,\big|\,dy}{\rho^N}=0,
\end{equation*}
where $\chi(a,b,\vec\nu)$ is defined by
\begin{equation*}
\chi(a,b,\vec\nu)(y):=
\begin{cases}
b\quad\text{if }\vec\nu\cdot y<0,\\
a\quad\text{if }\vec\nu\cdot y>0.
\end{cases}
\end{equation*}
The triple $(a,b,\vec\nu)$, uniquely determined, up to a permutation
of $(a,b)$ and a change of sign of $\vec\nu$, is denoted by
$(f^+(x),f^-(x),\vec\nu_f(x))$. We shall call $\vec\nu_f(x)$ the
{\em approximate jump vector} and we shall sometimes write simply
$\vec\nu(x)$ if the reference to the function $f$ is clear. The set
of approximate jump points is denoted by $J_f$. A choice of
$\vec\nu(x)$ for every $x\in J_f$ determines an orientation of
$J_f$. At an approximate continuity point $x$, we shall use the
convention $f^+(x)=f^-(x)=\tilde f(x)$.
\end{definition}
\begin{theorem}[Theorems 3.69 and 3.78 from \cite{amb}]\label{petTh}
Consider an open set $\Omega\subset\R^N$ and $f\in BV(\Omega,\R^m)$.
Then:\\
\noindent i) $\mathcal{H}^{N-1}$-a.e. point in
$\Omega\setminus J_f$ is a point of approximate continuity of $f$.\\
\noindent ii) The set $J_f$  is
$\sigma$-$\mathcal{H}^{N-1}$-rectifiable Borel set, oriented by
$\vec\nu(x)$. I.e. $J_f$ is $\sigma$-finite with respect to
$\mathcal{H}^{N-1}$, there exist countably many $C^1$ hypersurfaces
$\{S_k\}^{\infty}_{k=1}$ such that
$\mathcal{H}^{N-1}\Big(J_f\setminus\bigcup\limits_{k=1}^{\infty}S_k\Big)=0$,
and for $\mathcal{H}^{N-1}$-a.e. $x\in J_f\cap S_k$, the approximate
jump vector $\vec\nu(x)$ is  normal to $S_k$ at the point $x$.
\\ \noindent iii)
$\big[(f^+-f^-)\otimes\vec\nu_f\big](x)\in
L^1(J_f,d\mathcal{H}^{N-1})$.
\end{theorem}

\end{document}